\documentclass[12pt]{article}
\pdfoutput=1
\usepackage{geometry}
\usepackage[T1]{fontenc}
\usepackage[utf8]{inputenc}
\usepackage{lmodern}
\usepackage{microtype}

\usepackage{graphicx}
\usepackage[caption]{subfig}
\usepackage{mwe}
\usepackage{xcolor}
\usepackage{booktabs}
\usepackage{multirow}
\usepackage[export]{adjustbox}

\usepackage{amsmath,amssymb,amsfonts}
\usepackage{stmaryrd}
\usepackage{algorithmicx}
\usepackage{algpseudocode}
\usepackage{authblk}

\usepackage[
  backend=biber,
  style=numeric
]{biblatex}
\usepackage{hyperref}

\newenvironment{dedication}
  {\par\medskip\begin{center}\itshape}
  {\end{center}\par\medskip}

\newenvironment{funding}
  {\par\medskip\noindent\textbf{Funding:}\ }
  {\par\medskip}

\newenvironment{keywords}
  {\par\medskip\noindent\textbf{Keywords:}\ }
  {\par\medskip}

\newtheorem{algorithm}{Algorithm}
\newtheorem{proposition}{Proposition}

\newtheorem{remark}{Remark}
\newenvironment{proof}{\paragraph{Proof:}}{\hfill$\square$}

\newcommand{\R}{\mathbb{R}}
\newcommand{\J}{\mathcal{J}}

\newcommand{\bx}{\boldsymbol{x}}
\newcommand{\bp}{\boldsymbol{p}}

\newcommand{\epsstat}{\varepsilon_{\mathrm{stat}}}

\newcommand{\pmc}{\,\pm\,}

\newcommand\restr[2]{\ensuremath{\left.#1\right|_{#2}}}
\newcommand{\jump}[1]{\llbracket #1 \rrbracket}

\begin{document}

\title{Layerwise goal-oriented adaptivity for neural ODEs: an optimal control perspective}
\author[1,3]{Michael Hintermüller\thanks{Email: hintermueller@wias-berlin.de}}
\author[2]{Michael Hinze\thanks{Email: hinze@uni-koblenz.de}}
\author[3]{Denis Korolev\thanks{Email: korolev@wias-berlin.de}}
\affil[1]{\protect\raggedright Institute for Mathematics, Humboldt-Universität zu Berlin, Berlin, Germany}
\affil[2]{\protect\raggedright Mathematical Institute, Universität Koblenz, Koblenz, Germany}
\affil[3]{\protect\raggedright Weierstrass Institute for Applied Analysis and Stochastics, Berlin, Germany}
\date{}
\maketitle

\begin{dedication}
Dedicated to our esteemed colleague Ronald Hoppe, who passed away in February 2023.
\par
Ronald was an exceptional mentor and a good friend. He is truly missed.
\end{dedication}

\begin{abstract}
In this work, we propose a novel layerwise adaptive construction method for neural network architectures. Our approach is based on a goal--oriented dual-weighted residual technique for the optimal control of neural differential equations. This leads to an ordinary differential equation constrained optimization problem with controls acting as coefficients and a specific loss function. We implement our approach on the basis of a DG(0) Galerkin discretization of the neural ODE, leading to an explicit Euler time marching scheme. The resulting optimization problem is solved using the Adam algorithm and a BFGS method adapted to the $H^1$ topology induced by the regularization term. Finally, we apply our method to the construction of neural networks for the classification of data sets, where we present results for a selection of well known examples from the literature.
\end{abstract}


\section{Introduction}

In recent years, a new paradigm in mathematical machine learning has come into focus, which considers parameterized neural differential equations as continuum models for residual neural networks \cite{chen2018neural, dupont2019augmented, finlay2020train, haber2017stable, kidger2022neural, zhang2019anodev2, antil2023strong, antil2023optimal}. The hidden components of such a network are encoded in a state $\boldsymbol{x}$ ($=\boldsymbol{x}(t)$), which evolves continuously over the “time'' horizon\footnote{“Time'' is here merely artificial, rather than problem immanent.} $[0,T]$, $T>0$ and $t\in [0,T]$, according to a neural vector field $F(\boldsymbol{x}, \theta)$, thus transforming the (given) input data $\boldsymbol{x}_{\rm in} (=\boldsymbol{x}(0))$ into its feature representation $\boldsymbol{x}(T)$. The latter is then used to fit the output (label) data $\boldsymbol{y}$, which are given. This is typically achieved by learning the parameters $\theta$ ($=\theta(t)$) occurring in the neural differential equation via training data and the minimization of a suitable objective. Mathematically, this  leads to the solution of an optimization (or parameter identification) problem with ordinary  differential equation (ODE) constraints \cite{Benning2019JCD, haber2017stable, ruiz2023neural} reading
\begin{equation}\label{OCP intro}
\begin{aligned}
&\inf \  \mathcal{J}(\boldsymbol{x}, \theta):=J(\boldsymbol{x}) + \mathcal{R}(\theta), \ \   \text{over}  \ \ (\boldsymbol{x}, \theta),\\
&\text{subject to:} \quad  \dot{\boldsymbol{x}}=  F(\boldsymbol{x}, \theta), \quad \boldsymbol{x}(0) = \boldsymbol{x}_{\rm in},
\end{aligned}
\end{equation}
where $J$ typically represents the data fit term (thus depending also on $\boldsymbol{y}$) and $\mathcal{R}(\theta)$ is a suitable regularization ({\it prior}). In this setting, we refer to $\boldsymbol{x}$ as the {\it state} and $\theta$ as the {\it control} variable, respectively. 
This allows us to use analytical and computational tools developed within the realm of optimal control of differential equations.

More specifically, here we employ such tools to provide an answer to the question, “How the depth of a residual neural network should be adapted so that it satisfactorily solves a given task''. For this purpose, we rely on the concept of goal-oriented (mesh) adaptivity \cite{becker1999general, MR2421327,MR2745781,MR2642683, kraft2010dual} in our machine learning context. In fact, given an “ideal” continuous network $F(\, \cdot \,, \theta^{\star})$ that defines our {\it goal} (i.e., the target value for adaptation) $\mathcal{J}(\boldsymbol{x}^{\star}, \theta^{\star})$, we aim to construct a residual network, characterized by a discrete set of parameters $\theta^{\star}_{\tau}$, that fits this goal well. Accordingly, our adaptive approach is based on the error bound:
\begin{align}
\big|\mathcal{J}(\boldsymbol{x}^{\star},\theta^{\star})-\mathcal{J}(\boldsymbol{x}_{\tau}^{\star},\theta_{\tau}^{\star})\big|
\leq
\frac{1}{2}\sum_{k=1}^{K}\Big(
\big|\rho_{\mathrm{adj}}^{k}[\boldsymbol{x}^{\star}-\boldsymbol{x}_{\tau}^{\star}]\big|
+
\big|\rho_{\mathrm{grad}}^{k}[\theta^{\star}-\theta_{\tau}^{\star}]\big|
+
\big|\rho_{\mathrm{state}}^{k}[\boldsymbol{p}^{\star}-\boldsymbol{p}_{\tau}^{\star}]\big|
\Big)
+
|R|,
\end{align}
which contains the residual of the adjoint equation weighted by the corresponding error in the state, the residual of the gradient equation weighted by the error in the control, and the residual of the state equation weighted by the error in the adjoint. Here,  $R$ denotes higher-order terms which are neglected in the practical implementation, and $K\in\mathbb{N}$ is the number of time steps. Based on this error bound, we propose a computable indicator that drives our depth-adaptive neural network design. This generates a “time” grid, whose grid points are interpreted as layers of an “ideal” neural network.

The stationarity characterization in the form of the first-order optimality system for \eqref{OCP intro} and its corresponding discretization form the backbone of our approach. In this regard, it is crucial to ensure the existence of an optimal pair $(\boldsymbol{x}^{\star}, \theta^{\star})$ together with the associated Lagrange multiplier (adjoint state) $\boldsymbol{p}^{\star}$, with the triplet $(\boldsymbol{x}^{\star}, \theta^{\star}, \boldsymbol{p}^{\star})$ characterizing (one of) the ideal designs which we aim to meet. In order to achieve such a characterization in our functional-analytic setting, an appropriate regularization (or prior) is required. For this purpose, we stabilize \eqref{OCP intro} by adding an $H^{1}(0,T)$ regularization $\mathcal{R}$ to the objective $J$ that enforces smooth temporal variation of the network parameters and helps to show that an ideal design indeed exists. From an optimality perspective, this regularization induces the $H^{1}$ topology in which the gradient required for neural network optimization must be computed by solving a boundary value problem. On the discrete level, the latter is solved by combining a DG(0) Petrov–-Galerkin discretization for the state and adjoint equations with a continuous CG(1) Galerkin discretization for the gradient equation.

To the best of our knowledge, this approach to neural network design is new. Conceptually related approaches include the sensitivity-based insertion of layers into an existing neural network during training \cite{herberg2023sensitivity} and the adaptive successive approximation approach of \cite{aghili2024optimal}. In the latter, the optimality system resulting from the associated optimization problem is solved algorithmically with the aid of a fixed-point iteration (method of successive approximation), whereby in each iteration of the algorithm, the state and adjoint state are approximated on the basis of residual error estimators. The solution algorithm is therefore treated adaptively in its sub-steps. Our approach differs in that we derive the adaptivity directly from the optimality system. In this way, we take advantage of the fact that the optimality system mathematically represents a boundary value problem that can be treated adaptively using suitable methods.

The structure and contributions of the paper are as follows:
In Section \ref{sec 2: Neural ODEs}, the neural ODE is introduced and cast into its batch form. The corresponding optimal control problem together with a suitable functional-analytic framework is introduced in Section \ref{sec 3: OCP of neural ODEs}.  The first-order stationarity conditions are derived and analyzed in Section \ref{sec 4: optimality}.  In Section \ref{sec 5: discrete solver}, we propose our discretization of the optimality system. A dual-weighted a posteriori error representation for the training objective is derived in Section \ref{sec 6: DRW}, and the numerical realization of our layerwise algorithm is discussed in Section \ref{sec 7: Layerwise adaptive algorithm}. In Section \ref{sec 8: Numerical examples}, we report on numerical tests for our adaptive approach on binary and multiclass classification problems. The paper concludes with a summary and outlook.

\section{Neural ODEs}
\label{sec 2: Neural ODEs}
%Let $I=(0,T)$ with $T>0$. For any $q: I \rightarrow \mathbb{R}$, we sometimes write $q$ instead of $q(t)$ when the time argument is clear from context. We use the standard notation for Bochner spaces: 
%\begin{align*}
%H^{1}(I; \mathbb{R}^{d})= \{ w \in L^{2}(I; \mathbb{R}^{d}), \  \dot{w} \in   L^{2}(I; \mathbb{R}^{d}) \}.
%\end{align*}
%Define also notation for forms. 
The key idea of neural ODEs is to treat the depth of a neural network as a continuous variable, replacing discrete layers with a differential equation parameterized by weights and biases. Over a time interval $I=(0,T)$, with $T>0$, we consider the following time-dependent weight matrix and bias vector:
\begin{align}\label{coefficients}
W: I \rightarrow  \mathbb{R}^{d \times d}, \quad
b: I \rightarrow  \mathbb{R}^{d},
\end{align}
where $d\in\mathbb{N}$ denotes the dimension of the hidden state. We then collect these quantities into a time-dependent parameter vector $\theta: I \rightarrow \mathbb{R}^{n}$, defined by
\begin{align}\label{NF parameters}
\theta(t):=
\begin{pmatrix}
\operatorname{vec}W(t)\\ b(t)
\end{pmatrix},
\end{align}
where $n = d^2 + d$, and the operator $\mathrm{vec}: \mathbb{R}^{d \times d} \to \mathbb{R}^{d^2}$ stacks the columns of a matrix into a single vector. Let $f: \mathbb{R}^d \times \mathbb{R}^n \to \mathbb{R}^d$ denote the neural vector field
\begin{align}\label{neural field}
f(v, \theta(t)) := \boldsymbol{\sigma}(W(t)v + b(t)), \quad t\in [0,T], 
\end{align}
where $\boldsymbol{\sigma}: \mathbb{R}^{d} \to \mathbb{R}^{d}$ is the vector-valued activation function
\begin{align}
\boldsymbol{\sigma}(v) = (\sigma(v_{1}), \dots, \sigma(v_{d}) )^{\top},
\end{align}
and $\sigma : \mathbb{R} \rightarrow \mathbb{R}$ is a scalar activation function. Given an input $x_0 \in \mathbb{R}^{d_{\rm in}}$, $d_{\rm in}\in\mathbb{N}$, the neural ODE seeks a trajectory $x: [0,T] \to \mathbb{R}^d$ satisfying
\begin{align}\label{ResNet}
\begin{cases}
\dot{x}(t) = f(x(t), \theta(t)) \quad \text{for} \quad  t\in I,\\
x(0) = W_{\rm in}\, x_0,
\end{cases}
\end{align}
where $W_{\rm in} \in \mathbb{R}^{d \times d_{\rm in}}$ is the input projection matrix mapping the input data into the state space. Further, $\dot{x}$ denotes the time derivative of $x$ with respect to time $t$.

In supervised learning, a batch of $m$ training examples $\{(x_0^i, y^i)\}_{i=1}^m \subset \mathbb{R}^{d_{\rm in}} \times \mathbb{R}^{d_{\text{out}}}$, $d_{\rm out}\in\mathbb{N}$, is typically given.  
For both generality and practical implementation, it is convenient to express \eqref{ResNet} in a batch formulation.  For this, we define
\begin{align}\label{training data}
\boldsymbol{x}_{0} = (x^{1}_{0}, \ldots, x^{m}_{0})^{\top} \in  \mathbb{R}^{m d_{\text{in}}}, \quad \boldsymbol{y} = (y^{1}, \ldots, y^{m})^{\top} \in \mathbb{R}^{m d_{\text{out}}},
\end{align}
with $m\in\mathbb{N}$.
The batched neural vector field $F: \mathbb{R}^{md} \times \mathbb{R}^{n} \to \mathbb{R}^{md}$ acts component-wise as
\begin{align}\label{batched neural field}
F(\boldsymbol{x}, \theta)
:= \big(f(x^1, \theta), \ldots, f(x^m, \theta)\big)^{\top},
\quad
\boldsymbol{x} = (x^1, \ldots, x^m)^{\top} \in  \mathbb{R}^{md},
\end{align}
where each trajectory $x^i: [0,T] \to \mathbb{R}^d$ satisfies \eqref{ResNet} with $x^{i}(0)=W_{\text{in}}x_{0}^{i}$. The batched projected input is
\begin{align}\label{init cond}
\boldsymbol{x}_{\mathrm{in}} := (W_{\mathrm{in}} x_0^1, \ldots, W_{\mathrm{in}} x_0^m)^{\top} \in \mathbb{R}^{md}.
\end{align}
The $i$-th block of $\boldsymbol{x}$ is $x^i$; correspondingly,
$F^i(\boldsymbol{x},\theta)=f(x^i,\theta)$ and $\boldsymbol{x}_{\mathrm{in}}^i = W_{\mathrm{in}} x_0^i$. The resulting \emph{full-batch neural ODE} is then given by
\begin{align}\label{bResnet}
\begin{cases}
\dot{\boldsymbol{x}}(t)= F(\boldsymbol{x}(t), \theta(t)) \quad \text{for} \quad t\in I, \\
\boldsymbol{x}(0) = \boldsymbol{x}_{\rm in}.
\end{cases}
\end{align}
The system \eqref{bResnet} evolves the input $\boldsymbol{x}_{\text{in}}$ into a batch feature representation $\boldsymbol{x}(T)$. The output map $\boldsymbol{q}_{\mathrm{out}}: \mathbb{R}^{md} \to \mathbb{R}^{m d_{\text{out}}}$  transforms these features into predictions $\hat{\boldsymbol{y}} \in \mathbb{R}^{md_{\text{out}}}$, and is defined component-wise as
\begin{align}\label{output mapping}
\boldsymbol{q}_{\mathrm{out}}(\boldsymbol{x}(T))
= (\hat{y}^{1}, \ldots, \hat{y}^{m})^{\top}=:\hat{\boldsymbol{y}},
\end{align}
where $\hat{y}^{i}:=q_{\mathrm{out}}(W_{\mathrm{out}} x^i(T))$ for $i \in \{1,\ldots, m \}$, $q_{\mathrm{out}}: \mathbb{R}^{d_{\mathrm{out}}} \to \mathbb{R}^{d_{\mathrm{out}}}$ is a task-dependent continuous output function, and $W_{\mathrm{out}} \in \mathbb{R}^{d_{\mathrm{out}} \times d}$ is the output weight matrix.

A word on notation, before we continue: In this work, the Euclidean norm and dot product are denoted by $\lVert \, \cdot  \, \rVert$ and $(\cdot, \cdot)$, respectively. For batched vectors $\boldsymbol{x} \in \mathbb{R}^{md}$ as in \eqref{training data} or \eqref{batched neural field}, we use the batch norm associated with the product space $(\mathbb{R}^d)^{m}$, which is given by
\begin{align}\label{batch norm}
\lVert \boldsymbol{x} \rVert^{2} = \sum_{i=1}^{m} \lVert x^{i}\rVert^2.
\end{align}
For ease of notation, we assume that $W_{\rm in}$ and $W_{\rm out}$ are fixed, non-trainable parameters. This assumption does not affect the conclusions of our analysis, and extending the results to trainable $W_{\rm in}$ and $W_{\rm out}$ is straightforward. Nevertheless, where necessary, we provide computational formulas that allow these parameters to be included in the learning process.

\section{Optimal control of neural ODEs} 
\label{sec 3: OCP of neural ODEs}

The optimal control problem associated with \eqref{bResnet} consists in identifying optimal parameters $\theta$ based on the training data \eqref{training data} and a suitable learning objective. Here, we address this problem within the ``first optimize-then-discretize" approach, i.e., we first introduce an appropriate functional-analytic setting and derive stationarity conditions, which are then discretized for computational purposes.  We use the standard notation for Sobolev--Bochner spaces, such as $L^{2}(I; \mathbb{R}^{md})$ and $H^{1}(I; \mathbb{R}^{md})$; see, e.g., \cite{MR2500068}. Their respective norms are consistent with \eqref{batch norm} and are defined as follows: For $\boldsymbol{x} \in \mathbb{R}^{md}$ as in \eqref{batched neural field}, we have
\begin{align*}
\lVert \boldsymbol{x} \rVert_{L^{2}(I; \mathbb{R}^{md})}^{2} = \sum_{i=1}^{m} \int_{0}^{T} \lVert x^{i}(t)\rVert^{2} \, dt, \quad 
\lVert \boldsymbol{x} \rVert_{H^{1}(I; \mathbb{R}^{md})}^{2} = \lVert \dot{\boldsymbol{x}} \rVert_{L^{2}(I; \mathbb{R}^{md})}^{2} + \lVert \boldsymbol{x} \rVert_{L^{2}(I; \mathbb{R}^{md})}^{2}.
\end{align*}
Furthermore, $C([0,T]; \mathbb{R}^{d})$ denotes the space of continuous functions on $[0,T]$ with values in $\mathbb{R}^{d}$, equipped with the supremum norm $\lVert \boldsymbol{v} \rVert_{\infty} := \sup_{t \in [0,T]} \lVert \boldsymbol{v}(t) \rVert$. For brevity, we often omit the respective Euclidean space in the definition of function spaces in our proofs, specifically when the appropriate dimension is clear from the context.

We begin by discussing a suitable function space for the ODE parameters \eqref{NF parameters}. Regularization is naturally used to ensure that the quantities in \eqref{NF parameters} belong to a specific function class. For $\lambda>0$, the functional
\begin{align}\label{regularizer}
\mathcal{R}(\theta) = \frac{\lambda}{2} \int_0^T \left( \|\theta(t)\|^2 + \|\dot{\theta}(t)\|^2 \right) dt
\end{align}
ensures that the controls $\theta$ belong to $H^1(I)$. Indeed, without the derivative term in \eqref{regularizer}, one only has $\theta \in L^2(I)$, which admits discontinuities. We note that a regularization technique similar to \eqref{regularizer} was used in \cite{haber2017stable} to improve neural ODE stability within a ‘first discretize–then–optimize’ setting, where the derivative $\dot{\theta}(t)$ is replaced by its finite-difference approximation.

To guarantee the existence and uniqueness of solutions to \eqref{bResnet} in appropriate function spaces, we consider conditions on \eqref{NF parameters} that can be enforced by choosing a suitable regularization, as in \eqref{regularizer}. 
\begin{proposition}\label{Ex&Uniq}
Suppose that $\theta \in L^{2}(I; \mathbb{R}^{n})$ and that $\boldsymbol{\sigma} \in C(\mathbb{R}^{d})$ is Lipschitz continuous, i.e., there exists $L_\sigma > 0$ such that
\begin{align}\label{Lip act}
\lVert \boldsymbol{\sigma}(\boldsymbol{s}_2)- \boldsymbol{\sigma}(\boldsymbol{s}_1) \rVert \le L_\sigma \lVert \boldsymbol{s}_2 - \boldsymbol{s}_1 \rVert , \quad \forall \boldsymbol{s}_1, \boldsymbol{s}_2 \in \mathbb{R}^d.
\end{align}
Then the problem \eqref{bResnet} admits a unique solution $\boldsymbol{x} \in H^{1}(I; \mathbb{R}^{md})$. Moreover, there exist finite constants $C_1, C_2>0$ such that the following stability bound holds: 
\begin{align}\label{stability}
\lVert \boldsymbol{x} \rVert_{H^{1}(I; \mathbb{R}^{md})}\leq C_{1} \exp \big[ C_2 \lVert \theta \rVert_{L^{2}(I; \mathbb{R}^{n})}\big].
\end{align}
\end{proposition}
\begin{proof}
With a slight abuse of previous notation, but for the convenience of ODE theory, we define $f: \mathbb{R}^{d} \times [0,T] \rightarrow \mathbb{R}^{d}$ and $F : \mathbb{R}^{md} \times [0,T] \rightarrow \mathbb{R}^{md}$ as follows
\begin{align*}
f(x, t) := \boldsymbol{\sigma}(W(t)x + b(t)), \quad F(\boldsymbol{x}, t):= \big(f(x^1, t), \ldots, f(x^m, t)\big)^{\top}.
\end{align*}
Using the Lipschitz condition \eqref{Lip act}, for $\boldsymbol{v}, \boldsymbol{w} \in \mathbb{R}^{md}$, we estimate
\begin{equation}\label{Lip cond}
\begin{aligned}
\lVert F(\boldsymbol v,t)-F(\boldsymbol w,t) \rVert
=\Big[\sum_{i=1}^{m}\lVert f(\boldsymbol{v}^{i},t)-f(\boldsymbol{w}^{i},t)\lVert^{2}\Big]^{1/2}
\leq L_{\sigma}\lVert W(t) \lVert\,\lVert\boldsymbol v-\boldsymbol w\lVert.
\end{aligned}
\end{equation}
Using \eqref{Lip cond}, we get the following estimate
\begin{align*}
\lVert F(\boldsymbol{v},t) \rVert  &\leq   \lVert F(\boldsymbol{v},t) - F(0,t) \rVert +  \lVert F(0,t) \rVert  \leq  L_{\sigma}\lVert W(t) \lVert\,\lVert\boldsymbol v \lVert \, +  \, \lVert F(0,t) \rVert.
\end{align*}
Applying the Lipschitz condition \eqref{Lip act} to the second term above yields the following growth condition 
\begin{align}\label{growth condition}
\lVert F(\boldsymbol{v}, t) \rVert \leq  \sqrt{m} \, c(t) + a(t)  \lVert  \boldsymbol{v} \rVert,  
\end{align}
which holds for almost every $t \in [0,T]$. The coefficients in \eqref{growth condition} are given by
\begin{align}\label{coefficients growth}
a(t) = L_{\sigma} \lVert W(t) \rVert, \quad c(t) = L_{\sigma} \lVert b(t)\rVert + \lVert \boldsymbol{\sigma}(0) \rVert.
\end{align}
Since $\theta \in L^{2}(I)$, then $a, c \in L^{2}(I)$, and by the Cauchy–Schwarz inequality, we have $a, c \in L^{1}(I)$ as well. 

Note that $t \mapsto W(t)x+b(t)$ is not continuous, and $t \mapsto F(\boldsymbol{x},t)$ may fail to be continuous as well\footnote{Therefore, the famous Picard--Lindelöf theorem \cite[Theorem 3.1]{o1997existence} is not applicable. }. For this reason, we closely follow the approach of \cite[Theorem 3.4]{o1997existence}.  To begin with we define the norm 
\begin{align*}
\lVert \boldsymbol{v} \rVert_{A} = \lVert e^{-A} \boldsymbol{v}\rVert_{\infty}, \quad A(t) = \int_{0}^{t} a(s) \, ds. 
\end{align*}
Since $A(s) \leq  A(t)$ for $0\leq s \leq t$, we obtain the norm equivalence $e^{-A(T)} \lVert \boldsymbol{v} \rVert_{\infty} \leq \lVert \boldsymbol{v} \rVert_{A} \leq \lVert \boldsymbol{v} \rVert_{\infty}$, which makes $C([0,T])$ a Banach space when equipped with $\lVert \,\cdot \,\rVert_{A}$. For $\boldsymbol{v} \in C([0,T])$, we define an integral operator
\begin{align}\label{integral operator}
(\mathfrak{T}\boldsymbol{v})(t) = \boldsymbol{x}_{\mathrm{in}} + \int_{0}^{t} F(\boldsymbol{v}(s), s) \, ds.
\end{align}
From \eqref{growth condition}, we get $\|F(\boldsymbol{v}(s), s)\| \leq \sqrt{m} \, c(s) + a(s)M$, $M := \|\boldsymbol{v}\|_\infty$, ensuring that the integrand in \eqref{integral operator} is integrable and that $\mathfrak{T}$ is well-defined. For $0\leq t_1 < t_2$, we get
\begin{align*}
\|(\mathfrak{T}\boldsymbol{v})(t_2) - (\mathfrak{T}\boldsymbol{v})(t_1)\| \leq \int_{t_1}^{t_2}\lVert F(\boldsymbol{v}(s), s)\rVert \, ds  \leq  \int_{t_1}^{t_2} \big( \sqrt{m} \, c(s) + a(s)M \big)\, ds,
\end{align*}
where the right-hand side tends to $0$ as $t_2 \to t_1$ by the continuity of the Lebesgue integral. Hence, $\mathfrak{T}$ maps $C([0,T])$ into itself. We note that $A^{\prime}(t)= a(t)$ and apply \eqref{Lip cond} to deduce
\begin{align*}
\lVert \mathfrak{T}\boldsymbol{v} - \mathfrak{T}\boldsymbol{w}\rVert_A & \le \max_{t\in [0,T]}e^{-A(t)} \int_{0}^{t}\lVert F(\boldsymbol{v}(s),s) - F(\boldsymbol{w}(s),s)\rVert \, ds \leq \max_{t\in [0,T]} e^{-A(t)} \int_{0}^{t} a(s)\lVert \boldsymbol{v}(s) -\boldsymbol{w}(s)\rVert \,ds \\
&= \max_{t\in [0,T]} e^{-A(t)} \int_{0}^{t} A^{\prime}(s) e^{A(s)}e^{-A(s)}\lVert \boldsymbol{v}(s) -\boldsymbol{w}(s)\rVert \,ds \\
&  \leq \max_{t\in [0,T]} e^{-A(t)} \Big( \int_{0}^{t} A^{\prime}(s) e^{A(s)} \, ds \Big) \lVert \boldsymbol{v} -\boldsymbol{w} \rVert_{A} \\ 
& \leq \max_{t\in [0,T]} e^{-A(t)} \Big(e^{A(t)} - 1 \Big)\lVert \boldsymbol{v} -\boldsymbol{w} \rVert_{A} \leq \Big( 1 - e^{-A(T)}\Big)\lVert \boldsymbol{v} -\boldsymbol{w} \rVert_{A},
\end{align*}
Since $1 - e^{-A(T)} < 1$, the Banach fixed-point theorem guarantees a unique $\boldsymbol{x} \in C([0,T])$ such that $\mathfrak{T}\boldsymbol{x} = \boldsymbol{x}$. By the Lebesgue differentiation theorem, we get
\begin{align*}
\dot{\boldsymbol{x}}(t) = \frac{d}{dt} (\mathfrak{T}\boldsymbol{x})(t) = F(\boldsymbol{x}, t),
\end{align*}
which holds for almost every $t \in [0,T]$. The growth condition \eqref{growth condition} and the triangle inequality yield %concavity of $\sqrt{\cdot} \ $\footnote{$\sqrt{a+b} \leq\sqrt{2} \, (\sqrt{a}+\sqrt{b})$ for $a,b \in \mathbb{R}_{\geq 0}$.} then yield
\begin{equation}\label{der L2 est}
\begin{aligned}
\lVert \dot{\boldsymbol{x}} \rVert_{L^{2}(I)} &\leq  \big \lVert \sqrt{m} \,c + a \,  \lVert \boldsymbol{x} \rVert \big\rVert_{L^{2}(I)} \leq \sqrt{m} \, \lVert c \rVert_{L^{2}(I)} \, + \, \lVert \boldsymbol{x} \rVert_{\infty} \lVert a \rVert_{L^{2}(I)}. 
\end{aligned}
\end{equation}
Therefore, $\dot{\boldsymbol{x}} \in L^{2}(I)$. Hence, $\boldsymbol{x} \in H^{1}(I)$ and solves \eqref{bResnet}, as required. 

We proceed with the proof of the stability bound \eqref{stability}, beginning with the observation that 
\eqref{der L2 est} implies
\begin{align}\label{stab bound I}
\lVert \boldsymbol{x}\rVert_{H^{1}(I)} \leq  \lVert  \boldsymbol{x}\rVert_{L^{2}(I)} + \lVert \dot{\boldsymbol{x}} \rVert_{L^{2}(I)}  \leq  \big[\sqrt{T} +\lVert a \rVert_{L^2(I)} \big] \lVert \boldsymbol{x}\rVert_{\infty} + \sqrt{m} \, \lVert c \rVert_{L^2(I)}. 
\end{align}
Since $a(t) \geq 0$ and $c(t) \geq 0$ for almost every $t \in [0,T]$, the condition \eqref{growth condition} yields
\begin{align}\label{GR bound I}
\lVert \boldsymbol{x}(t) \rVert \leq \sqrt{m} \, \int_{0}^{t}c(s) \, ds + \int_{0}^{t} a(s) \lVert \boldsymbol{x}(s) \rVert \, ds
\end{align}
By the continuity of the Lebesgue integral, the first term on the right-hand side of \eqref{GR bound I} is a continuous, monotonically increasing function of $t$; applying Grönwall’s inequality \cite[Proposition 2.1]{emmrich1999discrete} then yields 
\begin{align*}
\lVert \boldsymbol{x}(t) \rVert \leq \sqrt{m} \,  \exp \Big(\int_{0}^{t} a(s) \, ds \Big) \, \int_{0}^{t}c(s) \, ds,  
\end{align*}
which holds for all $t\in[0,T]$. Using the Cauchy–Schwarz inequality and the fact that the integrands in the above estimate are non-decreasing, we obtain
%\mh{We could also go with the $L^1-$norms of $a$ and $c$. But this is not so important.}
\begin{align}\label{stab bound 2}
\lVert \boldsymbol{x} \rVert_{\infty} \leq \sqrt{T \,m} \, \exp \Big(\sqrt{T} \, \lVert a \rVert_{L^2{(I)}} \Big) \, \lVert c \rVert_{L^{2}(I)}.
\end{align}
Combining \eqref{stab bound 2} with \eqref{stab bound I} yields
\begin{align*}
\|\boldsymbol{x}\|_{H^{1}(I)}
\le \sqrt{m} \Big[
\sqrt{T} \, \big(\sqrt{T} + \|a\|_{L^2(I)}\big) \,
\exp\!\big(\sqrt{T} \, \|a\|_{L^2(I)}\big) + 1
\Big] \, \|c\|_{L^2(I)}.
\end{align*}
From this bound, it is clear that there exist constants $C_1, C_2>0$ depending only on $L_{\sigma}$, $m$, $\lVert \boldsymbol{\sigma}(0)\rVert$ and $T$ such that the stability bound \eqref{stability} holds. 
\end{proof}
\begin{remark}\label{remark on continuity}
If $\theta \in H^{1}(I)$, then by the continuous embedding $H^1(I) \hookrightarrow C([0,T])$ \cite{MR2424078}, the parameters in \eqref{NF parameters}, and thus the mapping $t \mapsto W(t)x + b(t)$, are continuous. Since $[0,T]$ is compact and $\theta \in C([0,T])$, the parameters in \eqref{NF parameters} are bounded on $[0,T]$. Hence, $C:= \sup_{t \in [0,T]}\|W(t)\|$ is finite and \eqref{Lip cond} becomes
\[
\begin{aligned}
\lVert F(\boldsymbol v,t)-F(\boldsymbol w,t) \rVert
\leq CL_{\sigma}\lVert\boldsymbol v-\boldsymbol w\lVert.
\end{aligned}
\]
independently of $t\in[0,T]$. Therefore, $F(\cdot,t)$ is Lipschitz continuous in $\boldsymbol{x}$, and the Lipschitz constant in \eqref{Lip cond} can be chosen independently of $t$. By the Picard--Lindelöf theorem, the problem \eqref{bResnet} then has a unique solution in $C^1([0,T])$.
\end{remark}
Proposition \ref{Ex&Uniq} yields the existence of the so-called {\it control-to-state map}
\begin{equation}\label{solution map}
\begin{aligned}
&\mathfrak{S}: L^{2}(I; \mathbb{R}^{n})  \rightarrow H^{1}(I; \mathbb{R}^{md}), \quad \theta \mapsto \boldsymbol{x}=:\mathfrak{S}(\theta). 
\end{aligned}
\end{equation}
This mapping is nonlinear, and therefore its boundedness and continuity must be proved separately. The boundedness of \eqref{solution map} readily follows from the stability bound \eqref{stability}. The next result establishes the continuous dependence of ODE solutions of \eqref{bResnet} on the control coefficients \eqref{NF parameters}, i.e., the continuity of \eqref{solution map}. 
%From the inclusion $H^{1}(I) \subset L^{2}(I)$, we readily infer that \eqref{solution map} is continuous as well. 
\begin{proposition}\label{continuity of S}
Suppose that $\boldsymbol{\sigma} \in C(\mathbb{R}^{d})$ and is Lipschitz continuous with some  Lipschitz constant $L_\sigma > 0$. Then the control-to-state mapping $\mathfrak{S}: L^{2}(I; \mathbb{R}^{n})  \rightarrow H^{1}(I; \mathbb{R}^{md})$ is continuous. 
\end{proposition}
\begin{proof}
Let $\theta_{1}, \theta_{2} \in L^{2}(I)$ be two distinct coefficients given by \eqref{NF parameters}, i.e.,
\begin{align*}
\theta_{1}(t)=
\begin{pmatrix}
\operatorname{vec}W_{1}(t)\\ b_{1}(t)
\end{pmatrix}, \quad \theta_{2}(t)=
\begin{pmatrix}
\operatorname{vec}W_{2}(t)\\ b_{2}(t)
\end{pmatrix}.
\end{align*}
According to Proposition \ref{Ex&Uniq}, $\theta_{1}$ and $\theta_{2}$ produce unique neural ODE solutions $\boldsymbol{x}_{1}, \boldsymbol{x}_{2} \in H^{1}(I)$ with the same initial value $\boldsymbol{x}_{\text{in}}$. Let $\delta \boldsymbol{x}(t)=\boldsymbol{x}_{1}(t) - \boldsymbol{x}_{2}(t)$. Observe that $\delta \boldsymbol{x} \in H^{1}(I)$ implies $\delta \boldsymbol{x} \in C([0,T])$ by the continuous embedding $H^{1}(I) \hookrightarrow C([0,T])$. The following inequality holds for almost every $t\in [0,T]$:  
\begin{align}\label{deco}
\lVert \delta \dot{\boldsymbol{x}}(t) \rVert \leq \big[\lVert  F(\boldsymbol{x}_{1}, \theta_1) - F(\boldsymbol{x}_{2}, \theta_1) \rVert    + \lVert  F(\boldsymbol{x}_{2}, \theta_1) - F(\boldsymbol{x}_{2}, \theta_2) \rVert\big](t).
\end{align}
Similar to \eqref{Lip cond}, from the Lipschitz continuity of $\boldsymbol{\sigma}$, we obtain the following bound for almost every $t \in [0,T]$: 
\begin{align}\label{t-Lip bound}
\big[\lVert  F(\boldsymbol{x}_{1}, \theta_1) - F(\boldsymbol{x}_{2}, \theta_1) \rVert \big](t) \leq a(t) \lVert \delta \boldsymbol{x} (t)\rVert,  
\end{align}
where $a \in L^{1}(I)$ is given by \eqref{coefficients growth}. By integrating inequality \eqref{deco}, using $\delta \boldsymbol{x}(0)=0$ and \eqref{t-Lip bound}, we get
\begin{align}\label{GR}
\lVert \delta \boldsymbol{x}(t) \rVert \leq  \int_{0}^{t}a(s) \lVert \delta \boldsymbol{x}(s) \rVert \, ds + \int_{0}^{t}\big[\lVert  F(\boldsymbol{x}_{2}, \theta_{1}) - F(\boldsymbol{x}_{2}, \theta_{2}) \rVert\big](s) \, ds.
\end{align}
By the continuity of the Lebesgue integral, the first term on the right-hand side of \eqref{GR} is a continuous, monotonically increasing function of $t$. Additionally, $a(t) \geq 0$ for almost every $t \in [0,T]$. Therefore, Grönwall’s inequality \cite[Proposition 2.1]{emmrich1999discrete} applied to \eqref{GR} yields
\begin{align*}
\lVert \delta \boldsymbol{x}(t) \rVert 
& \leq  \exp \Big(\int_{0}^{t} a(s) \, ds \Big)\int_0^t \big[\lVert F(\boldsymbol{x}_{2}, \theta_{1}) - F(\boldsymbol{x}_{2}, \theta_{2}) \rVert\big](s)  \, ds, 
\end{align*}
By the Cauchy--Schwarz inequality and monotonicity of the above integrals, we get
\begin{align*}
\lVert \delta \boldsymbol{x}(t) \rVert
 \leq \sqrt{T} \,C_{a}  \lVert F(\boldsymbol{x}_{2}, \theta_{1}) - F(\boldsymbol{x}_{2}, \theta_{2})\rVert_{L^{2}(I)}, 
\end{align*}
where $C_{a} = \exp \big(\int_{0}^{T} a(t) \, dt \big) < \infty$. We further use the Lipschitz condition on $\boldsymbol{\sigma}$, Young's inequality and the elementary inequality $\sqrt{z_{1}+z_{2}} \leq \sqrt{z_{1}}+\sqrt{z_{2}}$ for $z_{1}, z_{2} \in \mathbb{R}_{\geq 0}$, and obtain
\begin{align*}
\lVert  \delta \boldsymbol{x}(t) \rVert & \leq  \sqrt{T} \, C_{a}  \, L_{\sigma}  \Big[\sum_{i=1}^{m}\int_0^T  \big\lVert \big(W_{1}(t)- W_{2}(t)\big)\boldsymbol{x}_{2}^{i}(t) + b_{1}(t)- b_{2}(t) \big\lVert^{2} \, dt\Big]^{1/2} \\ 
& \leq \sqrt{2\, T} \, C_{a} \, L_{\sigma} \Big[\int_0^T \big\lVert W_{1}(t) - W_{2}(t) \big\rVert^{2}   \big\lVert \boldsymbol{x}_{2}(t)  \big\rVert^{2}  + m\int_0^T \big\rVert b_{1}(t)- b_{2}(t) \big\lVert^{2} \, dt\Big]^{1/2}  \\ 
& \leq \sqrt{2\, T} \, C_{a} \, L_{\sigma} \, \lVert \boldsymbol{x}_2\rVert_{\infty} \rVert W_{1} - W_{2} \rVert_{L^{2}(I)} + \sqrt{2 \, T \, m} \, C_{a} \, L_{\sigma}\,\rVert b_{1} - b_{2} \rVert_{L^{2}(I)},
\end{align*}
where we also used that $\lVert \boldsymbol{x}_{2}\rVert_{\infty} < \infty$. Since the above bound is uniform in $t$, it follows that  $\lVert  \delta \boldsymbol{x} \rVert_{\infty} \rightarrow 0$ as $\lVert \theta_{1} - \theta_{2} \lVert_{L^2} \rightarrow 0$. Therefore, $\lVert  \delta \boldsymbol{x} \rVert_{L^{2}(I)} \rightarrow 0$ as $\lVert \theta_{1} - \theta_{2} \lVert_{L^2} \rightarrow 0$ as well.

To complete the proof, it remains to show that $\lVert  \delta \dot{\boldsymbol{x}} \rVert_{L^{2}(I)} \rightarrow 0$ as $\lVert \theta_{1} - \theta_{2} \lVert_{L^2} \rightarrow 0$. We start by applying Young's inequality to \eqref{deco}, and obtain for almost every $t\in [0,T]$ the following
\begin{align*}
\lVert \delta \dot{\boldsymbol{x}}(t) \rVert^{2} \leq 2  \lVert  \big[F(\boldsymbol{x}_{1}, \theta_1) - F(\boldsymbol{x}_{2}, \theta_1) \big](t)\rVert^{2}    + 2\lVert  \big[F(\boldsymbol{x}_{2}, \theta_1) - F(\boldsymbol{x}_{2}, \theta_2)\big](t) \rVert^{2}.
\end{align*}
Integrating and applying \eqref{t-Lip bound} together with the previous estimates, we get
\begin{align*}
\lVert \delta \dot{\boldsymbol{x}} \rVert_{L^{2}(I)} \leq \sqrt{2} \lVert a\rVert_{L^{2}(I)} \lVert \delta \boldsymbol{x}\rVert_{\infty} + 2 L_{\sigma} \, \Big[\lVert \boldsymbol{x}_2\rVert_{\infty} \rVert W_{1} - W_{2} \rVert_{L^{2}(I)} + \sqrt{m} \,\rVert b_{1} - b_{2} \rVert_{L^{2}(I)} \Big].  
\end{align*}
Since $\lVert  \delta \boldsymbol{x} \rVert_{\infty} \rightarrow 0$ as $\lVert \theta_{1} - \theta_{2} \lVert_{L^2} \rightarrow 0$, the proof is complete. 
\end{proof}

After these preparatory results, we are ready to formulate an optimal control problem with neural ODEs. We introduce the following spaces: 
\begin{align}\label{spaces}
\mathcal{W}: = H^{1}(I; \mathbb{R}^{md}), \quad \mathcal{Q}:=L^{2}(I; \mathbb{R}^{md}), \quad  \mathcal{V}: = \mathcal{Q} \times \mathbb{R}^{md},  \quad \mathcal{U}:=H^{1}(I; \mathbb{R}^{n}).
\end{align}
In our problem, we aim to determine the parameters $\theta \in \mathcal{U}$ of \eqref{bResnet} from the training data \eqref{training data} by optimizing the following loss function 
\begin{align}\label{Cross-entropy function}
J: \mathcal{W} \rightarrow \mathbb{R}_{\geq 0}, \quad J(\boldsymbol{x}) =l(\boldsymbol{x}(T)),
\end{align}
where $l: \mathbb{R}^{md} \rightarrow \mathbb{R}$ is a continuous, problem-dependent function. We define the superposition operator
\begin{align}\label{Nemytskii op}
N_{F}: \mathcal{W} \times \mathcal{U} \rightarrow \mathcal{Q}, \quad (\boldsymbol{x}, \theta) \mapsto F(\boldsymbol{x}(\cdot), \theta(\cdot)),  
\end{align} 
and consider the optimal control (or parameter identification) problem
\begin{equation}\label{OCP std}
\begin{aligned}
&\inf \  \mathcal{J}(\boldsymbol{x}, \theta):=J(\boldsymbol{x}) + \mathcal{R}(\theta), \ \   \text{over}  \ \ (\boldsymbol{x}, \theta) \in \mathcal{W}\times \mathcal{U},\\
&\text{subject to:} \quad e(\boldsymbol{x}, \theta) = 0,
\end{aligned}
\end{equation}
where $\mathcal{J}: \mathcal{W} \times \mathcal{U} \rightarrow \mathbb{R}$ includes the regularizer \eqref{regularizer}, and the constraint is defined by
\begin{align*}
e: \mathcal{W} \times \mathcal{U} \rightarrow \mathcal{V}, \quad (\boldsymbol{x}, \theta) \mapsto e(\boldsymbol{x}, \theta): = \begin{pmatrix}
\dot{\boldsymbol{x} } - N_{F}(\boldsymbol{x}, \theta), \ \boldsymbol{x}(0) - \boldsymbol{x}_{\mathrm{in}} 
\end{pmatrix}^{\top} ,
\end{align*}
where $\dot{\boldsymbol{x}}$ is now understood in the generalized sense. Since $\mathcal{W} \hookrightarrow C([0,T]; \mathbb{R}^{md})$,  the initial value $\boldsymbol{x}(0)$ in the above equation is well-defined. Altogether, \eqref{OCP std} is an equality-constrained optimization problem, with the neural ODE defining the constraint. 

Next we establish the existence of minimizers for the problem \eqref{OCP std}, i.e., optimal pairs $(\boldsymbol{x}^{\star}, \theta^{\star}) \in \mathcal{W} \times \mathcal{U}$ satisfying $\mathcal{J}(\boldsymbol{x}^{\star},\theta^{\star}) \le \mathcal{J}(\boldsymbol{x}, \theta)$ for all $(\boldsymbol{x}, \theta) \in \mathcal{W} \times \mathcal{U}$. 
\begin{proposition}[Existence of optimal controls]\label{Existence of minimizers}
Suppose that $\theta \in \mathcal{U}$ with $\mathcal{U}:=H^{1}(I; \mathbb{R}^{n})$. Then the problem \eqref{OCP std} has an optimal solution. 
\end{proposition}
\begin{proof}
We define $F_{\mathrm{ad}} = \{ (\boldsymbol{x}, \theta) \in \mathcal{W} \times \mathcal{U} \, : \, \boldsymbol{x} = \mathfrak{S} (\theta) \}$. Proposition \ref{Ex&Uniq} implies that $F_{\mathrm{ad}} \neq \emptyset$. Therefore, $\inf \, \{\mathcal{J}(\boldsymbol{x},\theta) :  (\boldsymbol{x},\theta)  \in F_{\mathrm{ad}} \}=:j \in \mathbb{R}_{\geq 0}$. By properties of the infimum, we can pick a (infimizing) sequence $(\boldsymbol{x}_{k},\theta_{k})_{k\in\mathbb{N}}\subset F_{\text{ad}}$, such that $\mathcal{J}(\boldsymbol{x}_{k}, \theta_{k}) \rightarrow  j$ as $k \rightarrow \infty$. Furthermore, we have
\begin{align*}
\lim_{\lVert \theta \rVert_{\mathcal{U}} \to \infty} \mathcal{J}(\boldsymbol{x}, \theta)
\;\ge\;
\lim_{\lVert \theta \rVert_{\mathcal{U}} \to \infty} \mathcal{R}(\theta)
= +\infty,
\end{align*}
i.e., $\mathcal{J}$ is radially unbounded with respect to $\theta$. It then follows that $(\theta_{k})_{k\in\mathbb{N}}$ is bounded in $\mathcal{U}$; otherwise, if  $\lVert \theta_{k} \rVert_{ \mathcal{U}} \rightarrow \infty$ as $k\rightarrow \infty$, then necessarily $\mathcal{J}(\boldsymbol{x}_{k}, \theta_{k}) \rightarrow \infty$, contradicting the notion of an infimizing sequence. The stability bound \eqref{stability} yields that $(\boldsymbol{x}_{k})_{k\in\mathbb{N}}$ is bounded as well. Since $\mathcal{W}$ and $\mathcal{U}$ are Hilbert spaces and thus reflexive, the Banach--Alaoglu theorem implies that the sequence $(\boldsymbol{x}_{k}, \theta_{k})_{k\in\mathbb{N}}$ has a weakly convergent subsequence in $\mathcal{W} \times \mathcal {U}$, which we denote by $(\boldsymbol{x}_{k^{\prime}}, \theta_{k^{\prime}})_{k^{\prime}\in\mathbb{N}}$, with weak limit $(\boldsymbol{x}^{\star}, \theta^{\star}) \in \mathcal{W} \times \mathcal{U}$. From the compact embedding $\mathcal{U} \hookrightarrow C([0,T])$, we obtain that  $\theta_{k^{\prime}} \rightarrow \theta^{\star}$ strongly in $C([0,T])$. Since $C([0,T]) \subset L^{2}(I)$, Proposition \ref{continuity of S} implies that the control-to-state mapping $\mathfrak{S}: C([0,T]) \rightarrow \mathcal{W}$ is continuous, which in turn implies that $\boldsymbol{x}_{k^{\prime}} \rightarrow \boldsymbol{x}^{\star}$ strongly in $\mathcal{W}$. Since $\boldsymbol{x}^{\star}=\mathfrak{S}(\theta^{\star})$, we deduce that $(\boldsymbol{x}^{\star}, \theta^{\star}) \in F_{\mathrm{ad}}$. Since $J$ is continuous and $\mathcal{R}$ is weakly lower semicontinuous (because the norm is convex and continuous), we obtain
\begin{align*}
\mathcal{J}(\boldsymbol{x}^{\star}, \theta^{\star})  \leq  \,  \underset{k^{\prime} \rightarrow \infty}\liminf \, \mathcal{J}(\boldsymbol{x}_{k'}, \theta_{k'}) = \underset{k^{\prime} \rightarrow \infty}{\lim } \, \mathcal{J}(\boldsymbol{x}_{k'}, \theta_{k'})  = \underset{k \rightarrow \infty}{\lim } \, \mathcal{J}(\boldsymbol{x}_{k}, \theta_{k})=j,
\end{align*}
where the above subsequence of real numbers converges since the entire sequence converges, implying that the $\liminf$ and the limit indeed coincide. Therefore, $\mathcal{J}(\boldsymbol{x}^{\star}, \theta^{\star}) = j$, so that $(\boldsymbol{x}^{\star}, \theta^{\star})$ is a minimizer of $\mathcal{J}$.
\end{proof}
We note that upon elimination of the state variable in \eqref{OCP std} via \eqref{solution map}, the reduced objective $\widehat{\mathcal{J}}(\theta): = \mathcal{J}(\boldsymbol{x}(\theta), \theta)$ is generally non-convex. Therefore, minimizers of the reduced optimization problem need not be unique. We also mention that Proposition~\ref{Existence of minimizers} infers the strong continuity of \eqref{solution map} along minimizing sequences arising from $H^{1}$-regular controls. However, under weaker assumptions on the controls, such as $\mathcal{U}=L^{2}(I; \mathbb{R}^{n})$, it is somehow difficult to argue the existence of optimal controls via the direct method of calculus of variations.

\section{Stationarity conditions}
\label{sec 4: optimality}
In this section, we derive first-order optimality (or stationarity) conditions for the problem \eqref{OCP std}. 
We begin by considering the corresponding constraint qualification, which ensures the existence of a unique Lagrange multiplier (adjoint state) needed to characterize optimal solutions via the derivative of the Lagrangian. For this, we study the properties of the adjoint equation. In fact, the increased regularity of the adjoint variable allows us to reformulate the adjoint equation in a form suitable for our discretization. The function spaces needed in this section are defined in \eqref{spaces}. 

Let $D_{1}N_{F}(\boldsymbol{x}, \theta) : \mathcal{W} \rightarrow \mathcal{Q}$ denote the derivative of \eqref{Nemytskii op} with respect to its first argument: 
\begin{align*}
D_1 N_F(\boldsymbol{x},\theta) 
= \operatorname{diag}\Big(
D_1 f(x^1,\theta),\, 
D_1 f(x^2,\theta),\, 
\ldots,\, 
D_1 f(x^m,\theta) 
\Big).
\end{align*}
Since the batched vector field \eqref{batched neural field} acts component-wise across the batch, $D_1 N_F(\boldsymbol{x}, \theta)$ is block-diagonal. Each block stands for the Jacobian of \eqref{neural field} and is given by 
\begin{equation}\label{N-component}
\begin{aligned}
D_1 f(x^i(t), \theta(t)) = \operatorname{diag}\big(\boldsymbol{\sigma}'(W(t)x^i(t) + b(t))\big) \, W(t), 
\end{aligned}
\end{equation}
where $\operatorname{diag}\big(\boldsymbol{\sigma}'(v)\big) \in \mathbb{R}^{d\times d}$ is the diagonal matrix with $\boldsymbol{\sigma}'(v)\in \mathbb{R}^d$ along the main diagonal. The derivative $D_{2}N_{F}(\boldsymbol{x}, \theta) : \mathcal{U} \rightarrow \mathcal{Q}$ of \eqref{Nemytskii op} with respect to its second argument is given by 
\begin{align*}
D_{2}N_{F}(\boldsymbol{x},\theta)
=\begin{pmatrix}
D_2f(x^1,\theta)\\
\vdots\\
D_2f(x^m,\theta)
\end{pmatrix}.
\end{align*}
The next result shows that $D_1 N_F(\boldsymbol{x},\theta)$ and $D_2 N_F(\boldsymbol{x},\theta)$ are in fact Fr\'echet derivatives and that the constraint qualification for the equality-constrained optimization problem \eqref{OCP std} holds. Here, $L(\mathcal{U}, \mathcal{Q})$ denotes the space of linear and continuous operators from $\mathcal{U}$ to $\mathcal{Q}$, etc. 
\begin{proposition}[Constraint qualification]\label{constraint qualification}
Suppose that $\boldsymbol{\sigma}' \in L^{\infty}(\mathbb{R}^{d})$ and the selection \eqref{spaces} holds. Then we have the following: 
\begin{enumerate}
    \item $D_1 N_F(\boldsymbol{x}, \theta) \in L(\mathcal{Q}, \mathcal{Q})$ and $D_2 N_F(\boldsymbol{x}, \theta) \in L(\mathcal{U}, \mathcal{Q})$.
    \item The Jacobian $e^{\prime}(\boldsymbol{x}, \theta) \in L(\mathcal{W}\times \mathcal{U},  \mathcal{V})$ is surjective. 
\end{enumerate}
\end{proposition}
\begin{proof}
The linearity of $D_1 N_F(\boldsymbol{x}, \theta)$ and $D_2 N_F(\boldsymbol{x}, \theta)$ follows directly from their componentwise representations. Since $D_1 N_F(\boldsymbol{x}, \theta)$ is block-diagonal, for $\boldsymbol{v}\in\mathcal{Q}$ we obtain
\begin{align*}
\lVert D_1 N_F(\boldsymbol{x}, \theta)\boldsymbol{v}\rVert_{\mathcal{Q}}^2
= \sum_{i=1}^m \lVert D_1 f(x^i,\theta)v^i\rVert_{L^2(I)}^2 
\leq \sum_{i=1}^m \lVert D_1 f(x^i,\theta)\rVert_{L^\infty(I)}^2
\lVert v^i\rVert_{L^2(I)}^2
\leq C_1^2\lVert\boldsymbol{v}\rVert_{\mathcal{Q}}^2,
\end{align*}
where $C_1:=\lVert\boldsymbol{\sigma}'\rVert_{L^\infty(\mathbb{R}^d)}\lVert W\rVert_\infty$. The constant $C_1$ is finite since $\theta\in\mathcal{U}$, the embedding $\mathcal{U}\hookrightarrow C([0,T])$ is continuous, and $[0,T]$ is compact. Hence $D_1 N_F(\boldsymbol{x}, \theta)\in L(\mathcal{Q},\mathcal{Q})$. For $\xi\in\mathcal{U}$, we similarly have
\begin{align*}
\lVert D_2 N_F(\boldsymbol{x}, \theta)\xi\rVert_{\mathcal{Q}}^2
=\sum_{i=1}^m\lVert D_2 f(x^i,\theta)\xi\rVert_{L^2(I)}^2 
\leq \sum_{i=1}^m\lVert D_2 f(x^i,\theta)\rVert_{L^\infty(I)}^2
\lVert\xi\rVert_{L^2(I)}^2
\leq C_2^2\lVert\xi\rVert_{\mathcal{U}}^2,
\end{align*}
where $C_2:=\big[\sum_{i=1}^m\lVert D_2 f(x^i,\theta)\rVert_{L^\infty(I)}^2\big]^{1/2}$.
In view of \eqref{NF parameters}, we get
\begin{align*}
D_2 f(x^i,\theta)
=
\operatorname{diag}\big(\boldsymbol{\sigma}'(Wx^i+b)\big)
\begin{pmatrix}
(x^i)^\top\otimes I_d & I_d
\end{pmatrix}.
\end{align*}
Therefore, $\lVert D_2 f(x^i,\theta)\rVert_{L^\infty(I)}\leq
\lVert\boldsymbol{\sigma}'\rVert_{L^\infty(\mathbb{R}^d)}
\big(1+\lVert x^i\rVert_\infty^2\big)^{1/2}$. Since $\boldsymbol{x}\in\mathcal{W}$, the continuous embedding $\mathcal{W}\hookrightarrow C([0,T])$ and the compactness of $[0,T]$ imply that $C_2<\infty$. Thus $D_2 N_F(\boldsymbol{x}, \theta)\in L(\mathcal{U},\mathcal{Q})$, which proves the first claim.

Next, we show that the Jacobian $e^{\prime}(\boldsymbol{x}, \theta) \in L(\mathcal{W}\times \mathcal{U},  \mathcal{V})$ is surjective. It is given by
\begin{align*}
e^{\prime}(\boldsymbol{x}, \theta) = \begin{pmatrix}
\partial_{t} - D_1 N_F(\boldsymbol{x}, \theta) &  - D_2 N_F(\boldsymbol{x}, \theta)  \\ 
\delta_{0} & 0 
\end{pmatrix},
\end{align*}
where $\delta_{0} \in \mathcal{W}^{\ast}$ is the Dirac delta distribution concentrated at $0$. Equivalently, we need to show that for arbitrary $q = (q_{1}, q_{0}) \in \mathcal{V}$, the equation  $e^{\prime}(\boldsymbol{x}, \theta) \delta \boldsymbol{h} = q$ admits a solution $\delta \boldsymbol{h} = (\delta \boldsymbol{x}, \delta \theta) \in \mathcal{W} \times \mathcal{U}$. This equation can be written as 
\begin{equation}\label{evolution problem}
\begin{aligned}
\delta \dot{\boldsymbol{x}}   = D_{1}N_{F}(\boldsymbol{x}, \theta)\delta \boldsymbol{x} + g_1,  \quad \delta \boldsymbol{x}(0)
 = q_{0},
\end{aligned}
\end{equation}
where $g_{1}:=q_{1} + D_{2}N_{F}(\boldsymbol{x}, \theta)\delta \theta$ with $g_{1} \in \mathcal{Q}$ . Hence, it now suffices to show that, for any given $\delta \theta \in \mathcal{U}$, there exists a unique solution $\delta \boldsymbol{x} \in \mathcal{W}$ to the linear non-autonomous inhomogeneous ODE \eqref{evolution problem}. Since the mapping $t \mapsto D_{1}N_{F}(\boldsymbol{x}(t), \theta(t))$ belongs to $L^{\infty}(I; \mathbb{R}^{md \times md})$, one can construct a unique solution in $W^{1, \infty}(I; \mathbb{R}^{md})$ for the homogeneous counterpart of \eqref{evolution problem}  (the case $g_{1} = 0$) and further use a variation-of-constants approach to obtain a unique solution in $\mathcal{W}$ for \eqref{evolution problem}; see \cite[Chapter 5]{pazy2012semigroups} for the details.
\end{proof}

For our purposes, specifically for the subsequent Petrov–-Galerkin discretization, we require a variational form of \eqref{bResnet}, which reads: For $\theta \in \mathcal{U}$, find $\boldsymbol{x}  \in \mathcal{W}$ such that 
\begin{align}\label{varResnet}
\mathcal{F}(\boldsymbol{x} , \theta; \varphi) = 0, \quad \text{for all} \ \  \varphi \in \mathcal{V}, 
\end{align}
where the corresponding form $\mathcal{F}: \mathcal{W} \times \mathcal{U} \times \mathcal{V} \rightarrow \mathbb{R}$ is given by
\begin{align}\label{state form}
\mathcal{F}(\boldsymbol{x} , \theta; \varphi) : = \int_{0}^{T}\big(\dot{\boldsymbol{x} } - N_{F}(\boldsymbol{x} , \theta), \ \varphi_{1} \big) \ dt  + \big(\boldsymbol{x}(0) - \boldsymbol{x} _{\text{in}}, \varphi_{0}\big)
\end{align}
with the test functions $\varphi = (\varphi_{1},\varphi_{0}) \in \mathcal{V}$.
\begin{remark}
We closely follow the notation introduced in \cite{kraft2010dual}, where the functionals, as the one given in \eqref{state form}, depend nonlinearly on the arguments before the semicolon and linearly on those after it.
\end{remark}
To characterize minimizers of \eqref{OCP std}, we derive an optimality system by considering the derivative of the corresponding Lagrangian $\mathcal{L}: \mathcal{W} \times \mathcal{U} \times  \mathcal{V} \rightarrow \mathbb{R}$, which is given by
\begin{align}
\mathcal{L}(\boldsymbol{x}, \theta; \boldsymbol{p}): = \mathcal{J}(\boldsymbol{x}, \theta) - \mathcal{F}(\boldsymbol{x}, \theta; \boldsymbol{p}), 
\end{align}
where we used the variational form \eqref{state form} of the constraint $e(\boldsymbol{x}, \theta)$. Proposition~\ref{Existence of minimizers} guarantees the existence of an optimal pair $(\boldsymbol{x}^{\star}, \theta^{\star}) \in \mathcal{W} \times \mathcal{U}$. In addition, Proposition~\ref{constraint qualification} ensures that $e^{\prime}(\boldsymbol{x}^{\star}, \theta^{\star})$ is surjective; therefore, there exists a unique Lagrange multiplier (adjoint state) $\boldsymbol{p}^{\star}:=(\boldsymbol{z}^{\star}, \boldsymbol{q}^{\star}) \in \mathcal{V}$, where $\boldsymbol{z}^{\star} \in \mathcal{Q}$  and $\boldsymbol{q}^{\star} \in \mathbb{R}^{md}$, such that the following first-order optimality condition holds: 
\begin{align}
\mathcal{L}^{\prime}(\boldsymbol{x}^{\star}, \theta^{\star};  \boldsymbol{p}^{\star}, \varphi) = 0, \quad \forall{\varphi} \in \mathcal{W} \times \mathcal{U} \times \mathcal{V},   
\end{align}
or, equivalently, in component form: 
\begin{equation}\label{opt system abstract}
\begin{aligned}
&(\mathrm{adjoint \  equation}) &&
D_{1}\mathcal{J}(\boldsymbol{x}^{\star}, \theta^{\star}; \varphi_x)
- D_{1}\mathcal{F}(\boldsymbol{x}^{\star}, \theta^{\star}; \boldsymbol{p}^{\star}, \varphi_x)
&&= 0, && \forall \varphi_x \in \mathcal{W}, \\[0.25ex]
&(\mathrm{gradient\ equation}) &&
D_{2}\mathcal{J}(\boldsymbol{x}^{\star}, \theta^{\star}; \varphi_\theta)
- D_{2}\mathcal{F}(\boldsymbol{x}^{\star}, \theta^{\star}; \boldsymbol{p}^{\star}, \varphi_\theta)
&&= 0, && \forall \varphi_\theta \in \mathcal{U}, \\[0.25ex]
& 
(\mathrm{state\ equation}) &&
\mathcal{F}(\boldsymbol{x}^{\star}, \theta^{\star}; \varphi_p)
&&= 0, && \forall \varphi_p \in \mathcal{V}.
\end{aligned}
\end{equation}
We now study the optimality system \eqref{opt system abstract}. From the first equation in \eqref{opt system abstract}, we derive a variational form of the adjoint equation: find $\boldsymbol{p}=(\boldsymbol{z}, 
\boldsymbol{q})\in \mathcal{V}$ such that 
\begin{equation}\label{adjoint equation 1}
\begin{aligned}
\int_{0}^{T} \big (\dot{\boldsymbol{w}}, \boldsymbol{z}) \ dt -\int_{0}^{T} \big(D_{1}N_{F}(\boldsymbol{x}, \theta) \boldsymbol{w}, \boldsymbol{z} \big) \ dt &+ \big( \boldsymbol{w}(0), \boldsymbol{q} \big)  =  \langle J^{\prime}(\boldsymbol{x}), \boldsymbol{w} \rangle_{\mathcal{W}^{\ast}, \mathcal{W}}, 
\end{aligned}
\end{equation}
for all $\boldsymbol{w} \in \mathcal{W}$. Due to the chosen loss function \eqref{Cross-entropy function}, which depends only on the evaluation of the state $\boldsymbol{x} \in \mathcal{W}$ at the terminal time $T$, it holds that
\begin{align}\label{adjoint RHS}
\langle J^{\prime}(\boldsymbol{x}), \boldsymbol{w} \rangle_{\mathcal{W}^{\ast}, \mathcal{W}} 
= \langle l^{\prime}(\boldsymbol{x}(T)) \, \delta_{T}, \boldsymbol{w} \rangle_{\mathcal{W}^{\ast}, \mathcal{W}} 
= \big( l^{\prime}(\boldsymbol{x}(T)), \boldsymbol{w}(T) \big),
\end{align}
where $\delta_{T} \in H^{1}(I; \mathbb{R}^{md})^{\ast}$ is the Dirac delta distribution concentrated at $T$ and $ l^{\prime}(\boldsymbol{x}(T)) \in \mathbb{R}^{md}$. The following holds true for the adjoint variable $\boldsymbol{p}=(\boldsymbol{z},\boldsymbol{q})$ in the adjoint equation \eqref{adjoint equation 1}, where we omit the $\star$ for better readability.
\begin{proposition}\label{Adjoint regularity}
Suppose that $\boldsymbol{\sigma}' \in L^{\infty}(\mathbb{R}^{d})$ and the selection \eqref{spaces} holds. Then we have the following: 
\begin{enumerate}
\item The multiplier $\boldsymbol{z}$ satisfies $\boldsymbol{z} \in \mathcal{W}$.
\item The multiplier $\boldsymbol{q}$ is equal to $\boldsymbol{z}(0)$ and the adjoint equation \eqref{adjoint equation 1} can be equivalently written as follows: find $\boldsymbol{z} \in \mathcal{W}$ such that 
\begin{equation}\label{adjoint problem 2}
\begin{aligned}
-\int_{0}^{T}  \big(\dot{\boldsymbol{z}}, \varphi_{1}\big) \, dt  &= \int_{0}^{T}  \big( D_{1}N_{F}(\boldsymbol{x}, \theta)^{\ast} \boldsymbol{z}, \varphi_{1} \big) \, dt, \\ 
(\boldsymbol{z}(T),\varphi_{0})  &= (l^{\prime}(\boldsymbol{x}(T)), \varphi_{0})
\end{aligned}
\end{equation}
for all $\varphi = (\varphi_{1},\varphi_{0}) \in \mathcal{V}$. 
\end{enumerate}
\end{proposition} 
\begin{proof}
We test \eqref{adjoint equation 1} with $ \boldsymbol{w}\in \mathcal{D}(I)$, where $\mathcal{D}(I):=C^{\infty}_{c}(I)$ is the space of smooth functions with compact support in $I$. Due to the compact support, the boundary term in \eqref{adjoint equation 1} and the derivative \eqref{adjoint RHS} vanish, yielding
\begin{align}
\int_{0}^{T} (\dot{\boldsymbol{w}}, \boldsymbol{z})\, dt
= -\langle \dot{\boldsymbol{z}}, \boldsymbol{w} \rangle_{\mathcal{D}', \mathcal{D}} = \int_{0}^{T} \big( D_{1}N_{F}(\boldsymbol{x}, \theta)^{\ast} \boldsymbol{z}, \boldsymbol{w} \big)\, dt,
\end{align}
where $\dot{\boldsymbol{z}} \in \mathcal{D}'(I)$ denotes the distributional time derivative of $\boldsymbol{z}$. However, $\dot{\boldsymbol{z}}$ is represented by
$-D_{1}N_{F}(\boldsymbol{x}, \theta)^{\ast} \boldsymbol{z}$. 
Since $D_{1}N_{F}(\boldsymbol{x}, \theta) \in L(\mathcal{Q})$ implies 
$D_{1}N_{F}(\boldsymbol{x}, \theta)^{\ast} \in L(\mathcal{Q})$, and $\boldsymbol{z} \in \mathcal{Q}$, we conclude that $\dot{\boldsymbol{z}} \in \mathcal{Q}$, which proves the claim.

For the second claim, we first apply generalized integration by parts to the first term in \eqref{adjoint equation 1} with $\boldsymbol{w} \in C^{\infty}([0,T])$, which yields
\begin{align}\label{int by parts}
\int_{0}^{T} \big(\boldsymbol{z}, \dot{\boldsymbol{w}} \big) \, dt = \big( \boldsymbol{z}(T), \boldsymbol{w}(T)  \big) - \big ( \boldsymbol{z}(0) , \boldsymbol{w}(0)\big) - \int_{0}^{T} \big (\dot{\boldsymbol{z}}, \boldsymbol{w} \big) \, dt.
\end{align}
The traces of $\boldsymbol{z} \in H^{1}(I)$ in \eqref{int by parts} are well-defined due to the continuity of the embedding $H^{1}(I) \hookrightarrow C([0,T])$, and the last term in \eqref{int by parts} is meaningful, as the first claim implies that $\dot{\boldsymbol{z}}$ is also the weak derivative of $\boldsymbol{z}$.  Combining  \eqref{adjoint equation 1}, \eqref{adjoint RHS} and \eqref{int by parts}, we get 
\begin{align*}
- \int_{0}^{T} \big ( \dot{\boldsymbol{z}}, \boldsymbol{w} \big) \, dt - \int_{0}^{T} \big(D_{1}N_{F}(\boldsymbol{x}, \theta)^{\ast} \boldsymbol{z}, \boldsymbol{w} \big) \, dt &+ \big(\boldsymbol{q} - \boldsymbol{z}(0), \boldsymbol{w}(0) \big) + \big(\boldsymbol{z}(T)-l^{\prime}(\boldsymbol{x}(T)), \boldsymbol{w}(T) \big)=0.
\end{align*} 
The fundamental lemma of the calculus of variations yields $\boldsymbol{q} = \boldsymbol{\boldsymbol{z}}(0)$.  Finally, using the density of $C^{\infty}([0,T])$ in $\mathcal{Q}$, we obtain \eqref{adjoint problem 2}.
\end{proof}
%\begin{remark}
%The assumptions also hold for...\textcolor{red}{no need for bias bounded} We define a uniform parameter box for the neural ODE parameters \eqref{NF parameters} as
%\begin{align}
%\mathbb{B}_R := \Big\{ \theta : I \to \mathbb{R}^{n} \;:\; 
%\|\theta(t)\|_\infty \le R, \ \forall t \in I \Big\} \subset L^{\infty}(I; \mathbb{R}^{n}), 
%\end{align}
%where $\|\theta(t)\|_\infty := \max\Big\{|W_{ij}(t)|, \ |b_j(t)| : 1 \le i,j \le d \Big\}.$ By ensuring that  
%\(|W_{ij}(t)|, |b_j(t)| \le R\) for all \(t \in I\), we guarantee that $\theta \in \mathbb{B}_R$. 
%\end{remark}

\begin{remark}
Inspecting Proposition \ref{Adjoint regularity} reveals that the increased $H^1$ regularity of $\boldsymbol{z}$ depends not only on the regularity of $\theta$, but also on the choice of the objective function. In our case, the objective \eqref{Cross-entropy function} depends only on the terminal state $\boldsymbol{x}(T)$. However, if the objective involves the evaluation of  $\boldsymbol{x}(\widetilde{T})$ at $\widetilde{T} < T$, one obtains
\begin{align*}
\dot{\boldsymbol{z}} = - \big(D_1 N_F(\boldsymbol{x}, \theta)^{\ast} \boldsymbol{z} + l'(\boldsymbol{x}(\widetilde{T})) \, \delta_{\widetilde{T}} \big),
\end{align*}
where $\delta_{\widetilde{T}} \in H^{1}(I; \mathbb{R}^{md})^{\ast}$ is the Dirac delta distribution concentrated at $\widetilde{T}$. Consequently, we expect that $\boldsymbol{z} \in L^2(I; \mathbb{R}^{md})$ only. We note that such objectives appear naturally in time-series modeling with neural ODEs; see \cite{kidger2022neural} and references therein.
\end{remark}
%\mh{The Riesz representer of $D_{2}\mathcal{L}(\boldsymbol{x}, \theta; \boldsymbol{p}, \cdot) \in \mathcal{U}^\ast$
%is the gradient of the reduced cost functional
%$\hat{\mathcal{J}}(\theta) := \mathcal{J}(\boldsymbol{x}(\theta), \theta)$,
%where $\boldsymbol{p} = \boldsymbol{p}(\boldsymbol{x}(\theta))$; cf.~\cite{hinze2008optimization}.}
We henceforth identify the adjoint with its trajectory component and write $\boldsymbol{p}:=\boldsymbol{z}$ and proceed by examining the structure of the gradient equation in \eqref{opt system abstract}. We recall that the Riesz representer of $D_{2}\mathcal{L}(\boldsymbol{x}, \theta; \boldsymbol{p}, \cdot) \in \mathcal{U}^\ast$
is the gradient $g$ of the reduced cost functional
$\widehat{\mathcal{J}}(\theta) = \mathcal{J}(\boldsymbol{x}(\theta), \theta)$,
where $\boldsymbol{p} = \boldsymbol{p}(\boldsymbol{x}(\theta))$; cf. \cite{hinze2008optimization}. That is, the gradient $g \in \mathcal{U}$ is the unique element satisfying
\begin{align}\label{gradient equation}
\langle g, \varphi_{\theta} \rangle_{\mathcal{U}}
= \int_{0}^{T} \lambda\Big[ (\theta, \varphi_{\theta}) +  (\dot{\theta}, \dot{\varphi}_{\theta}) \Big]\, + \big( D_{2}N_{F}(\boldsymbol{x},\theta)^{\ast} \boldsymbol{p}, \varphi_{\theta} \big)\, dt
, \quad \forall \varphi_{\theta} \in \mathcal{U}.
\end{align}
For convenience, we provide the strong form of \eqref{opt system abstract}, where we make use of \eqref{adjoint problem 2} for the adjoint problem: 
\begin{equation}\label{strong opt system}
\begin{aligned}
-\dot{\boldsymbol{p}}^{\star}& =  D_{1}N_{F}(\boldsymbol{x}^{\star}, \theta^{\star})^{\ast} \boldsymbol{p}^{\star}, \quad \boldsymbol{p}^{\star}(T)  =  l^{\prime}(\boldsymbol{x}^{\star}(T)), \\[0.5ex] 
\lambda \big(\theta^{\star} - \ddot{\theta}^{\star} \big)  &= -D_{2}N_{F}(\boldsymbol{x}^{\star},\theta^{\star})^{\ast}\boldsymbol{p}^{\star}, \quad \dot{\theta}^{\star}(0) = \dot{\theta}^{\star}(T)=0, \\[0.5ex] 
\dot{\boldsymbol{x}}^{\star} &= N_F(\boldsymbol{x}^{\star}, \theta^{\star}), \quad \boldsymbol{x}^{\star}(0) = \boldsymbol{x}_{\text{in}}.
\end{aligned}
\end{equation}
For a given $\theta$, the gradient $g$ in \eqref{gradient equation} is then obtained from the two-point boundary value problem
\begin{equation}\label{grad-strong}
g - \ddot{g} = \lambda \big(\theta - \ddot{\theta} \big) + D_{2}N_{F}(\boldsymbol{x},\theta)^{\ast}\boldsymbol{p}, \quad 
\dot{g}(0) = \dot{g}(T)=0,
\end{equation}
where $\boldsymbol{x}$ and $\boldsymbol{p}$ are the state and adjoint state trajectories associated to $\theta$.
The weak form of the problem \eqref{grad-strong} is then given in \eqref{gradient equation}. This boundary value problem structure stems from the regularizer \eqref{regularizer} inducing the $H^{1}(I)$ topology in which the gradient must now be computed.

\begin{remark}
To include $W_{\text{in}}$ and $W_{\text{out}}$ as trainable parameters, we need the  derivatives $D_{3}\mathcal{J}(\boldsymbol{x}, \theta, W_{\text{in}}, W_{\text{out}})$ and $D_{4}\mathcal{J}(\boldsymbol{x}, \theta, W_{\text{in}}, W_{\text{out}})$. The former is given by
\begin{equation*}
\begin{aligned}
D_{3}\mathcal{J}(\boldsymbol{x}, \theta, W_{\text{in}}, W_{\text{out}}) = \sum_{i=1}^{m} \boldsymbol{p}^{i}(0) (\boldsymbol{x}_{0}^{i})^{\top}.  
\end{aligned}
\end{equation*}
$D_{4}\mathcal{J}(\boldsymbol{x}, \theta, W_{\text{in}}, W_{\text{out}})$ is independent of the neural vector field in \eqref{ResNet} and depends only on the problem-specific loss \eqref{Cross-entropy function} and the output mapping $\boldsymbol{q}_{\mathrm{out}}$; it can therefore be computed using standard calculus.
\end{remark}

\section{Petrov--Galerkin discretization and numerical solver}
\label{sec 5: discrete solver}
In this section, we consider the discretization of the optimality system \eqref{opt system abstract}. We employ a discontinuous Petrov–-Galerkin discretization for both the state and the adjoint state, taking advantage of the higher regularity of the adjoint variable, as established in Proposition \ref{Adjoint regularity}. Both the state and adjoint discretizations lead to DG(0) schemes, which can be interpreted as the forward and backward Euler time-stepping schemes, respectively. For the discretization of the gradient equation \eqref{gradient equation}, we use a CG(1) Galerkin approximation, where the piecewise-linear in time parametrization of the control parameters is consistent with the regularization scheme \eqref{regularizer}.
%ensuring that the derivative $\dot{\theta}(t)$ remains non-vanishing on the discrete level. 

We start by defining a temporal partition of the time interval $[0,T]$ as follows:
\begin{align*}
\mathcal{T}_{K} = \{ 0= t_{0} < t_{1} <  \cdots < t_{K-1} < t_{K} = T\},
\end{align*}
and denote $I_{k}:=(t_{k-1}, t_{k})$, and the midpoint of $I_k$ by $t_{k-\frac{1}{2}}$. Furthermore, let $h_{k} = t_{k} - t_{k-1}$, for all ${k}=1,\ldots, K$, and $h: = \max \ h_{k}$. The subscript $\tau$ indicates dependence on the generally nonuniform temporal mesh $\{h_k\}_{k=1}^{K}$ and will subsequently be used for the corresponding discrete functions and spaces. We define the jump of a function $v$ at instant $t_{k}$ as $\jump{v}^{k}: = v(t_{k}^{+}) - v(t_{k}^{-})$, where
\begin{align*}
v(t_{k}^{+}): = \underset{s \rightarrow 0^{+}}{\lim} \ v(t_{k} + s), \quad v(t_{k}^{-}): = \underset{s \rightarrow 0^{+}}{\lim} \ v(t_{k} - s). 
\end{align*}
For the discretization of the trial space $\mathcal{W}$, we select the following discrete space: 
\begin{align*}
\mathcal{W}_{\tau}:  = \{\boldsymbol{w}_{\tau} \in \mathcal{Q} :  \restr{\boldsymbol{w}_{\tau}}{[t_{k-1}, t_{k})} \in \mathbb{P}_{r}([t_{k-1}, t_{k})  ;  \mathbb{R}^{md}), \ \forall{k}=1, \ldots ,K, \  \boldsymbol{w}_{\tau}(T) \in \mathbb{R}^{md} \}, 
\end{align*}
where $\mathbb{P}_{r}(I_{k}; \mathbb{R}^{md})$ is the space of polynomials on the interval $I_{k}$ taking values in $\mathbb{R}^{md}$, with polynomial degree less than or equal to $r$. Note that $\mathcal{W}_{\tau} \not\subset \mathcal{W}$; that is, the above discretization is non-conforming for the trial space. However, we have $\mathcal{W}_{\tau} \subset \mathcal{V}$.  For the discretization of the test space $\mathcal{V}$, we choose the following subspace of $\mathcal{V}$:  
\begin{align*}
\mathcal{V}_{\tau}:  = \{\boldsymbol{v}_{\tau} \in \mathcal{Q} :  \restr{\boldsymbol{v}_{\tau}}{(t_{k-1}, t_{k}]} \in \mathbb{P}_{r}((t_{k-1}, t_{k}]  ;  \mathbb{R}^{md}), \ \forall{k}=1, \ldots ,K, \ \boldsymbol{v}_{\tau}(0) \in \mathbb{R}^{md} \}.
\end{align*}
Then, we have $\mathcal{V}_{\tau} \subset \mathcal{V}$; that is, the discretization is conforming in the test space. Note that $\mathcal{W}_{\tau}$ and $\mathcal{V}_{\tau}$ both have dimension $((r+1)K + 1)md$. For the discretization of the gradient equation \eqref{gradient equation} we use the conforming discrete space $\mathcal{U}_{\tau} \subset \mathcal{U}$, which we specify later in the section.

We define the discrete state equation with the discontinuous Galerkin form $\mathcal{F}_{\text{DG}} :\mathcal{W}_{\tau} \times \mathcal{U}_{\tau} \times \mathcal{V}_{\tau} \rightarrow \mathbb{R}$ (see Appendix~\ref{sec:App A} for the derivation):  for $\theta_{\tau} \in \mathcal{U}_{\tau}$, find $\boldsymbol{x}_{\tau} \in \mathcal{W}_{\tau}$ such that  
\begin{equation}\label{VF state}
\begin{aligned}
\mathcal{F}_{\text{DG}}(\boldsymbol{x}_{\tau}, \theta_{\tau}; \varphi):& = \sum_{k=1}^{K} \bigg[ \int_{I_{k}} \big(\dot{\boldsymbol{x}}_{\tau} - N_F(\boldsymbol{x}_{\tau}, \theta_\tau),  \varphi_{1} \big) \, dt  + \big(\jump{\boldsymbol{x}_{\tau}}^{k}, \varphi_{1}(t_{k}^{-})\big) \bigg]  + \big(\boldsymbol{x}_{\tau}(0) - \boldsymbol{x}_{\text{in}},   \varphi_0\big)=0
\end{aligned}
\end{equation}
for all $\varphi=(\varphi_{1}, \varphi_{0}) \in \mathcal{V}_{\tau}$. The formulation \eqref{VF state} yields $((r+1)K + 1) md$ nonlinear algebraic equations in the same number of unknowns. Indeed, the solvability of this system follows from the Lipschitz continuity of $F(\cdot, \theta)$, yielding the discrete control-to-state map 
\begin{align}\label{eq: discrete control-to-state map}
\mathfrak{S}_{\tau} : \mathcal{U}_{\tau} \rightarrow \mathcal{W}_{\tau}, 
\quad \theta_{\tau} \mapsto \boldsymbol{x}_{\tau} =: \mathfrak{S}_{\tau}(\theta_{\tau}),
\end{align}
which is continuous due to the continuity of the map $t\mapsto F(\cdot, \theta(t))$, see Remark~\ref{remark on continuity}. Using \eqref{eq: discrete control-to-state map}, we define the discrete reduced objective by
\begin{align}\label{eq: discrete reduced objective}
\widehat{\mathcal{J}}_{\tau}(\theta_{\tau})
:=
\mathcal{J}
\bigl(
\mathfrak{S}_{\tau}(\theta_{\tau}),
\theta_{\tau}
\bigr).
\end{align} 
The existence of a discrete optimal control pair $(\boldsymbol{x}_{\tau}^{\star}, \theta_{\tau}^{\star}) \in \mathcal{W}_{\tau} \times \mathcal{U}_{\tau}$ can then be argued along the lines of Proposition~\ref{Existence of minimizers}, which is further simplified by the equivalence of weak and strong convergence in finite-dimensional spaces. By linearizing \eqref{VF state} at $(\boldsymbol{w}_{\tau}, \theta_{\tau})$ as in Proposition~\ref{constraint qualification}, we obtain $((r+1)K + 1) md$ linear algebraic equations, whose solvability yields the surjectivity of the Jacobian arising from the linearization. Then, a discrete stationary point $(\boldsymbol{x}_{\tau}^{\star}, \theta_{\tau}^{\star}; \boldsymbol{p}_{\tau}^{\star}) \in \mathcal{W}_{\tau} \times \mathcal{U}_{\tau} \times \mathcal{V}_{\tau}$ exists and satisfies
the discrete necessary optimality condition 
\begin{align}\label{discrete opt cond}
\mathcal{L}^{\prime}_{\mathrm{DG}}(\boldsymbol{x}_{\tau}^{\star}, \theta_{\tau}^{\star}; \boldsymbol{p}_{\tau}^{\star}, \varphi) = 0, \quad \forall{\varphi} \in \mathcal{W}_\tau \times \mathcal{U}_\tau \times \mathcal{V}_\tau,   
\end{align}
where $\mathcal{L}_{\mathrm{DG}}(\boldsymbol{x}_{\tau}, \theta_{\tau}; \boldsymbol{p}_{\tau}):= \mathcal J(\boldsymbol{x}_{\tau},\theta_\tau) -\mathcal{F}_{\text{DG}}(\boldsymbol{x}_{\tau}, \theta_\tau; \boldsymbol{p}_{\tau})$.
In component form we then have
\begin{equation}\label{opt system discrete}
\begin{aligned}
D_{1}\mathcal{J}(\boldsymbol{x}_{\tau}^{\star}, \theta_{\tau}^{\star}; \varphi_x) 
- D_{1}\mathcal{F}_{\text{DG}}(\boldsymbol{x}_{\tau}^{\star}, \theta_{\tau}^{\star}; \boldsymbol{p}_{\tau}^{\star}, \varphi_x) 
&= 0, 
&& \forall\, \varphi_x \in \mathcal{W}_\tau, \\[0.5ex]
D_{2}\mathcal{J}(\boldsymbol{x}_{\tau}^{\star}, \theta_\tau^{\star}; \varphi_\theta) 
- D_{2}\mathcal{F}_{\text{DG}}(\boldsymbol{x}_{\tau}^{\star}, \theta_\tau^{\star}; \boldsymbol{p}_{\tau}^{\star}, \varphi_\theta) 
&= 0, 
&& \forall\, \varphi_\theta \in \mathcal{U}_\tau, \\[0.5ex]
\mathcal{F}_{\text{DG}}(\boldsymbol{x}_{\tau}^{\star}, \theta_\tau^{\star}; \varphi_p) 
&= 0, 
&& \forall\, \varphi_p \in \mathcal{V}_\tau.
\end{aligned}
\end{equation}
The discrete adjoint equation for the optimality system \eqref{opt system discrete} reads: find $\boldsymbol{p}_{\tau} \in \mathcal{V}_{\tau}$ such that
\begin{equation}\label{adjoint equation}
\begin{aligned}
\sum_{k=1}^{K} \int_{I_{k}} \big( - \dot{\boldsymbol{p}}_{\tau} - D_{1}N_{F}(\boldsymbol{x}_{\tau}, \theta_\tau)^{\ast}\boldsymbol{p}, \varphi_{x} \big) \ dt   &- \sum_{k=1}^{K}\big (\jump{\boldsymbol{p}_{\tau}}^{k-1}, \varphi_x (t_{k-1}^{+}))  = 0 \\
\big( \boldsymbol{p}_{\tau}(T), \varphi_{x}(T) \big) & =  \big(l^{\prime}(\boldsymbol{x}_{\tau}(T)), \varphi_{x}(T) \big),
\end{aligned}
\end{equation} 
for all $\varphi_x \in \mathcal{W}_{\tau}$. The derivation of \eqref{adjoint equation} is provided in Appendix~\ref{sec:App A}.  We then consider a discretization that yields an explicit Euler scheme forward in time for the discrete state equation in \eqref{opt system discrete} and backward in time for the discrete adjoint equation \eqref{adjoint equation}, which we state later in the section. Specifically, we use piecewise-constant in time polynomial functions, corresponding to $r=0$ in $\mathcal{W}_{\tau}$ and $\mathcal{V}_{\tau}$, cf. also \cite{munoz2019explicit}. We refer to Appendix~\ref{sec:App B} for derivation details of this time-marching interpretation of our DG(0) scheme.

Let us now proceed with the numerical discretization of the gradient equation \eqref{gradient equation}. Let $\mathcal{U}_{\tau} \subset \mathcal{U}$ denote the finite-dimensional space of piecewise-linear in time finite-element functions with values in $\mathbb{R}^{n}$: 
\begin{align*}
\mathcal{U}_\tau
:=
\Big\{
u_{\tau} \in \mathcal{U} :
\ u_{\tau}|_{[t_{k-1},t_k]} \in \mathcal{P}_1([t_{k-1},t_k];\mathbb{R}^n),\ 
k=1,\dots,K
\Big\} \cong \mathcal{Q}_{\tau} \otimes \mathbb{R}^{n}, 
\end{align*}
where $\mathcal{Q}_{\tau} \subset H^{1}(I)$ consists of piecewise-linear, globally continuous $\mathbb{R}$-valued functions, and is spanned by the well-known nodal hat functions $\{ \phi^{i} \}_{i=0}^{K}$ \cite{MR2322235}.  Both boundary nodes $i=0$ and $i=K$ are included due to the Neumann boundary conditions in \eqref{grad-strong}. In view of the above, $g_{\tau} \in \mathcal{U}_{\tau}$ and $\theta_\tau \in \mathcal{U}_{\tau}$ are given by
\begin{align}\label{vectorial ansatz functions}
g_{\tau}(t)= \sum_{k=0}^{K} g_{\tau}^{k} \  \phi^{k}(t), \quad \theta_\tau(t) = \sum_{k=0}^{K}\theta_{\tau}^{k} \  \phi^{k}(t).
\end{align}
Consistently with the DG(0) state–adjoint discretization and the midpoint evaluation of the control in the neural field of  our time-marching interpretation, we define the piecewise-constant control pullback on $I_k$ by
\begin{align}\label{eq:control pullbacks}
q_{\tau}(t)
=
q_{\tau}^{k}
:=
D_{2}N_{F}
\bigl(
\boldsymbol{x}_{\tau}^{k-1},
\theta_{\tau}^{k-\frac12}
\bigr)^{\ast}
\boldsymbol{p}_{\tau}^{k}
\in\mathbb{R}^{n} \quad
\mathrm{for} \quad t\in I_k.
\end{align}
The above term can be computed using the backward mode automatic differentiation. The corresponding term in the continuous gradient equation \eqref{gradient equation} is then approximated by 
\begin{align}\label{control pullback approximation}
\left\langle
D_{2}N_{F}
(\boldsymbol{x}_{\tau},\theta_{\tau})^{\ast}
\boldsymbol{p}_{\tau}, \,
\varphi_{\theta}\right\rangle_{L^{2}(I;\mathbb{R}^{n})} & \approx
\sum_{k=1}^{K}
\int_{I_k}\big(q_{\tau}^{k}, \varphi_{\theta}(t)\big)\,dt.
\end{align}
The nodal vectors $g_{\tau}^{k}, \theta_{\tau}^{k} \in \mathbb{R}^{n}$ in \eqref{vectorial ansatz functions} and $\boldsymbol{q}_{\tau}^{k} \in \mathbb{R}^{n}$  in \eqref{eq:control pullbacks} are stored row-wise:
\begin{align}\label{eq:mat coefficients}
\mathbf G_\tau
=
\begin{pmatrix}
(g_\tau^0)^\top\\[-0.2em]
\vdots\\[-0.2em]
(g_\tau^K)^\top
\end{pmatrix}
\in\R^{(K+1)\times n},
\quad
\mathbf\Theta_\tau
=
\begin{pmatrix}
(\theta_\tau^0)^\top\\[-0.2em]
\vdots\\[-0.2em]
(\theta_\tau^K)^\top
\end{pmatrix}
\in\R^{(K+1)\times n}, \quad 
\mathbf{Q}_{\tau}
= \begin{pmatrix}
(q_\tau^1)^\top\\[-0.2em]
\vdots\\[-0.2em]
(q_\tau^K)^\top 
\end{pmatrix}\in \mathbb{R}^{K \times n}.  
\end{align}
Furthermore, $\mathbf{B}^{\tau} = \mathbf{A}^{\tau} + \mathbf{M}^{\tau}$, where $\mathbf{A}^{\tau}, \mathbf{M}^{\tau}\in \mathbb{R}^{(K+1) \times (K+1)}$ are the stiffness and the mass matrices of $\mathcal{Q}_{\tau}$-basis functions with the respective entries $\mathbf{A}^{\tau}_{k,l}=(\dot{\phi}^{k}, \dot{\phi}^{l})_{L^{2}(I)}$ and $\mathbf{M}^{\tau}_{k,l}=(\phi^{k}, \phi^{l})_{L^{2}(I)}$. For coefficient matrices $\mathbf{V}_{\tau}, \mathbf{W}_{\tau} \in \mathbb{R}^{(K+1) \times n}$ representing $v_{\tau}, w_{\tau} \in \mathcal{U}_{\tau}$, the discrete $H^{1}$ inner product is given by
\begin{align}\label{discrete H1 inner prod}
\langle v_{\tau},w_{\tau}\rangle_{H^{1}(I;\mathbb{R}^{n})}
=
\operatorname{trace} \,(\mathbf{V}_{\tau}^{\top}\,
\mathbf{B}^{\tau}\,
\mathbf{W}_{\tau}).
\end{align}

We use the approximation \eqref{control pullback approximation}, $g_{\tau}$ as the discrete trial function and $\theta_\tau$ as data in the continuous gradient equation \eqref{gradient equation}, and test the latter with $\varphi_{\theta}= \phi^{l}e_{j}$, where $l=0,\ldots, K$,  $j=1,\ldots, n$ and $e_j$ is the $j$-th canonical basis vector of $\mathbb{R}^n$. Assembling the resulting equations yields the following algebraic form of the discrete gradient equation:
\begin{align}\label{dg eq}
\mathbf{B}^{\tau}\,\mathbf{G}_{\tau} = \mathbf{R}_{\tau}, 
\end{align}
where the right-hand side $\mathbf{R}_{\tau} \in \mathbb{R}^{(K+1)\times n}$ is given by
\begin{align}\label{eq: gradient RHS}
\mathbf{R}_{\tau}:=\lambda \, \mathbf{B}^{\tau} \, \mathbf{\Theta}_{\tau} + 
\mathbf{C}^{\tau}\mathbf{Q}_{\tau}.
\end{align}
Here, $\mathbf{C}^{\tau}\in\mathbb{R}^{(K+1)\times K}$ is the mixed DG(0)--CG(1) matrix with entries
\begin{align*}
\mathbf{C}_{lk}^{\tau}
:=
\int_{I_k}\phi^{l}(t)\,dt
=
\frac{h_k}{2}
\bigl(
\delta_{l,k-1}
+
\delta_{l,k}
\bigr),
\qquad
0\leq l\leq K,
\quad
1\leq k\leq K,
\end{align*}
where $\delta_{l,k}$ is the Kronecker symbol. Since $\mathbf{B}^{\tau}$ is symmetric positive definite and tridiagonal, the $n$ columns of $\mathbf{G}_{\tau}$ are computed by the same tridiagonal solve with $n$ right-hand sides corresponding to the columns of $\mathbf{R}_{\tau}$, which can be done efficiently in parallel.

In summary, on a temporal grid $\mathcal{T}_{K}$ we find an approximation of the discrete stationary point $(\boldsymbol{x}_{\tau}^{\star}, \theta_{\tau}^{\star}, \boldsymbol{p}_{\tau}^{\star})$ of \eqref{opt system discrete} by considering the following system: 
\begin{equation}\label{full discrete optimality sys}
\begin{alignedat}{2}
\boldsymbol{x}_\tau^0
&=
\boldsymbol{x}_{\mathrm{in}},
\\[0.4ex]
\theta_\tau^{k-\frac12}
&=
\frac12
\left(
\theta_\tau^{k-1}+\theta_\tau^k
\right),
&\qquad&
k=1,\ldots,K,
\\[0.4ex]
\boldsymbol{x}_\tau^k
&=
\boldsymbol{x}_\tau^{k-1}
+
h_k
F\left(
\boldsymbol{x}_\tau^{k-1},
\theta_\tau^{k-\frac12}
\right),
&&
k=1,\ldots,K,
\\[0.4ex]
\boldsymbol{p}_\tau^K
&=
l'\left(
\boldsymbol{x}_\tau^K
\right),
\\[0.4ex]
\boldsymbol{p}_\tau^{k-1}
&=
\boldsymbol{p}_\tau^k
+
h_k
D_{1} F \left(
\boldsymbol{x}_\tau^{k-1},
\theta_\tau^{k-\frac12}
\right)^\ast
\boldsymbol{p}_\tau^k,
&&
k=K,\ldots,1,
\\[0.4ex]
\mathbf B^\tau\mathbf G_\tau
&=\mathbf{R}_{\tau}.
\end{alignedat}
\end{equation}
The above system is evaluated at a sequence of nodal coefficient matrices $\{ \mathbf \Theta_{\tau}^{(\mathrm{it})}\}_{\mathrm{it}\geq 0}$, where $\mathrm{it} \in \mathbb{Z}_{\geq0}$ denotes the iteration count. At a current iterate $\mathbf \Theta_{\tau}^{(\mathrm{it})}$, the forward Euler scheme (forward pass) is used to obtain $\boldsymbol{x}_{\tau}^{K}$, which is further used to compute $l^{\prime}(\boldsymbol{x}^{K}_{\tau})$ to initialize the backward Euler scheme for the adjoint approximation (backward pass). The latter is then used to assemble the right-hand side \eqref{eq: gradient RHS} of the discrete gradient equation \eqref{dg eq}, which is subsequently solved for $\mathbf{G}_{\tau}^{(\mathrm{it})}$. This drives an optimizer-dependent update
\begin{align}\label{optimization update}
\mathbf{\Theta}_{\tau}^{(\mathrm{it}+1)} = \mathbf\Theta_{\tau}^{(\mathrm{it})} + \alpha_{\mathrm{it}} \mathbf{D}_{\tau}^{(\mathrm{it})},
\end{align}
where $\mathbf{G}_{\tau}$ is used to construct an update direction $\mathbf{D}_{\tau}^{(\mathrm{it})}$, and $\alpha_{\mathrm{it}}>0$ is the step size. We want to point out early on that $\mathbf{D}_{\tau}^{(\mathrm{it})}$ may not be a descent direction. The updated coefficient matrix defines the next piecewise-linear control through the ansatz \eqref{vectorial ansatz functions}. This
procedure is repeated until the stationarity measure
\begin{align}\label{stationarity measure}
\mathfrak{s}_{\tau}(\mathbf{\Theta}_{\tau}): = \lVert g_{\tau} \rVert_{H^{1}(I; \mathbb{R}^{n})} = \big[\operatorname{trace} ( \mathbf{G}_{\tau}^{\top} \mathbf{B}^{\tau}\mathbf{G}_{\tau})\big]^{1/2}
\end{align}
falls below the prescribed tolerance $\epsstat >0$ or the solver iteration budget is exhausted.

We further comment on the numerical optimization required for the update \eqref{optimization update}. Identifying the discrete control $\theta_{\tau}$ in the objective \eqref{eq: discrete reduced objective} with its coefficient matrix $\mathbf\Theta_{\tau}$ gives the fully discrete reduced objective
\begin{align}\label{fully discrete objective}
\mathfrak{J}_{\tau}(\mathbf\Theta_{\tau}) = l(\boldsymbol{x}_{\tau}^{K}( \mathbf{\Theta}_{\tau})) + \frac{\lambda}{2} \, \operatorname{trace} ( \mathbf{\Theta}_{\tau}^{\top} \mathbf{B}^{\tau}\mathbf \Theta_{\tau}),
\end{align}
which is used to monitor the performance. For our optimization, we use the Adam algorithm \cite{kingma2014adam} as a warm start, followed by BFGS. One step of Adam requires computing the moments
\begin{align}
\mathrm{\mathbf{m}}^{(\mathrm{it})} = \beta_{1}  \mathrm{\mathbf{m}}^{(\mathrm{it-1})} + (1- \beta_{1}) \mathbf{G}_{\tau}^{(\mathrm{it})}, \quad \mathrm{\mathbf{v}}^{(\mathrm{it})} = \beta_{2}  \mathrm{\mathbf{v}}^{(\mathrm{it-1})} + (1- \beta_{2}) \mathbf{G}_{\tau}^{(\mathrm{it})} \odot \mathbf{G}_{\tau}^{(\mathrm{it})},
\end{align}
where $\odot$ is the Hadamard product, $\beta_1,\beta_2\in[0,1)$ and $\mathrm{\mathbf m}^{(-1)}=\mathrm{\mathbf v}^{(-1)}= \mathbf 0 \in
\mathbb{R}^{(K+1)\times n}$. The corresponding bias-corrected moments are given by
\begin{align}\label{eq: Adam moments}
\widehat{\mathrm{\mathbf m}}^{(\mathrm{it})}
=
\frac{
\mathrm{\mathbf m}^{(\mathrm{it})}
}{
1-\beta_1^{\mathrm{it}+1}
}, \qquad 
\widehat{\mathrm{\mathbf v}}^{(\mathrm{it})}
=
\frac{
\mathrm{\mathbf v}^{(\mathrm{it})}
}{1-\beta_2^{\mathrm{it}+1}}.
\end{align}
The Adam update direction for \eqref{optimization update} is then defined by
\begin{align}\label{eq:adam-direction}
\mathbf D_\tau^{(\mathrm{it})}
=
-
\frac{
\widehat{\mathrm{\mathbf m}}^{(\mathrm{it})}
}{
\sqrt{
\widehat{\mathrm{\mathbf v}}^{(\mathrm{it})}
}
+
\varepsilon_{\mathrm{Ad}}
},
\end{align}
where the square root and division are applied component-wise and $\varepsilon_{\mathrm{Ad}}>0$ prevents division by zero.  In our case, Adam provides an efficient warm start. Near stationarity, its momentum and adaptive scaling may cause oscillations that prevent the stationarity measure \eqref{stationarity measure} from crossing the prescribed tolerance $\epsstat$.

To increase the likelihood of reaching the stationarity tolerance $\epsstat$, we switch from Adam to BFGS once the Adam iterate enters the switching region $\mathfrak{s}_{\tau}(\mathbf{\Theta}_{\tau}^{(\mathrm{it})}) \leq \kappa_{s} \, \epsstat$, where $\kappa_{s}>1$. We note that standard finite-dimensional BFGS implementations operate on a vectorization of $\mathbf{\Theta}_\tau$ and use the Euclidean inner product, which corresponds to the Frobenius inner product in the present matrix notation. However, the natural inner product for $\mathbf{\Theta}_\tau$ is the discrete $H^1$ inner product \eqref{discrete H1 inner prod}. A metric-aware BFGS method could incorporate this inner product directly into its line search and quasi-Newton update \cite{kelley1987quasi}. To use a standard Euclidean implementation, we instead introduce the change of coordinates
\begin{align}\label{Cholesky transform}
\widehat{\mathbf{\Theta}}_{\tau}
=(\mathbf{L}^{\tau})^{\top}\mathbf{\Theta}_{\tau}, \qquad \mathbf{\Theta}_{\tau}
=(\mathbf{L}^{\tau})^{-\top} \, \widehat{\mathbf{\Theta}}_{\tau},
\end{align}
where $\mathbf{L}^{\tau}$ comes from the Cholesky factorization $\mathbf{B}^{\tau}=\mathbf{L}^{\tau}(\mathbf{L}^{\tau})^{\top}$ of the discrete Riesz matrix. For $\mathbf{V}_{\tau}, \mathbf{W}_{\tau} \in \mathbb{R}^{(K+1) \times n}$ representing $v_{\tau}\in \mathcal{U}_{\tau}$ and $w_{\tau}\in \mathcal{U}_{\tau}$, we define $\widehat{\mathbf{V}}_{\tau}
=(\mathbf{L}^{\tau})^{\top}\mathbf{V}_{\tau}$,  $\widehat{\mathbf{W}}_{\tau}
=(\mathbf{L}^{\tau})^{\top}\mathbf{W}_{\tau}$, and obtain 
\begin{align}\label{isometry property}
(\widehat{\mathbf{V}}_{\tau},
 \widehat{\mathbf{W}}_{\tau}
)_{\mathbf{F}} = \operatorname{trace}\, (\widehat{\mathbf{V}}_{\tau}^{\top} \, \widehat{\mathbf{W}}_{\tau})
= \operatorname{trace}\, (\mathbf{V}_{\tau}^{\top}
\mathbf{B}^{\tau}
\mathbf{W}_{\tau})
=
\langle 
v_{\tau}, w_{\tau}
\rangle_{H^{1}(I; \mathbb{R}^{n})}.
\end{align}
Hence, \eqref{Cholesky transform} defines an isometry from the discrete control space $\mathcal{U}_{\tau}$ to the coefficient space equipped with $(\cdot, \cdot)_{\mathbf{F}}$. In particular, $\lVert \widehat{\mathbf{\Theta}}_{\tau}\rVert_{\mathbf{F}} = \lVert \theta_{\tau}\rVert_{H^{1}(I; \mathbb{R}^{n})}$. We define the transformed objective by
\begin{align}\label{modified reduced objective}
\widetilde{\mathfrak{J}}_{\tau}(\widehat{\mathbf\Theta}_{\tau}): = \mathfrak{J}_{\tau}((\mathbf{L}^{\tau})^{- \top} \, \widehat{\mathbf{\Theta}}_{\tau}).
\end{align}
The chain rule then yields the Euclidean gradient of the transformed objective: 
\begin{align*}
\widehat{\mathbf{G}}_{\tau}: = \nabla_{\widehat{\Theta}_{\tau}} \widetilde{\mathfrak{J}}_{\tau}(\widehat{\Theta}_{\tau}) = (\mathbf{L}^{\tau})^{-1} \mathbf{R}_{\tau} = (\mathbf{L}^{\tau})^{\top} \mathbf{G}_{\tau}.
\end{align*}
Note that $\lVert\widehat{\mathbf{G}}_{\tau}\rVert_{\mathbf{F}} = \mathfrak{s}_{\tau}(\mathbf{\Theta}_{\tau})$. For a given iterate $\mathrm{it}\in \mathbb{Z}_{\geq0}$, we set
\begin{align*}
\widehat{\boldsymbol{\theta}}_{\tau}^{(\mathrm{it})}
=
\operatorname{vec}_{\tau}
\bigl(\widehat{\mathbf{\Theta}}_{\tau}^{(\, \mathrm{it})}\bigr), \quad 
\widehat{\boldsymbol{g}}_{\tau}^{(\mathrm{it})}=
\operatorname{vec}_{\tau}
\bigl(\widehat{\mathbf{G}}_{\tau}^{(\mathrm{it})}\bigr).
\end{align*} 
where $\operatorname{vec}_{\tau}: \mathbb{R}^{(K+1)\times n} \rightarrow \mathbb{R}^{(K+1)n}$ denotes vectorization, and
$\operatorname{mat}_{\tau}=\operatorname{vec}_{\tau}^{-1}$ is its inverse. The BFGS update is then determined from the equation
\begin{align}\label{vectorized BFGS direction}
\widehat{\mathbf{H}}_{\tau}^{\,(\mathrm{it})}
\widehat{\boldsymbol{d}}_{\tau}^{\,(\mathrm{it})}
= - \widehat{\boldsymbol{g}}_{\tau}^{\,(\mathrm{it})}.
\end{align}
where $\widehat{\mathbf{H}}_{\tau}^{\,(\mathrm{it})}$ is the positive-definite Hessian approximation. The update for \eqref{optimization update} is then given by
\begin{align}\label{BFGS direction transformation}
\mathbf{D}_{\tau}^{(\mathrm{it})}
=
(\mathbf{L}^{\tau})^{-\top}\, \widehat{\mathbf{D}}_{\tau}^{(\mathrm{it})},
\end{align}
where $\widehat{\mathbf{D}}_{\tau}^{(\mathrm{it})}=\operatorname{mat}_{\tau} (\widehat{\boldsymbol{d}}_{\tau}^{\,(\mathrm{it})})$. The secant vectors for the BFGS update are given by
\begin{align*}
\widehat{\boldsymbol{s}}_{\tau}^{(\mathrm{it})}=
\widehat{\boldsymbol{\theta}}_{\tau}^{(\mathrm{it}+1)}
-
\widehat{\boldsymbol{\theta}}_{\tau}^{(\mathrm{it})}
=
\alpha_{\mathrm{it}}
\widehat{\boldsymbol{d}}_{\tau}^{(\mathrm{it})}, \quad 
\widehat{\boldsymbol{\gamma}}_{\tau}^{(\mathrm{it})}
=
\widehat{\boldsymbol{g}}_{\tau}^{(\mathrm{it}+1)}
-
\widehat{\boldsymbol{g}}_{\tau}^{(\mathrm{it})}.
\end{align*}
The BFGS Hessian approximation is then updated according to
\begin{align}\label{BFGS Hessian update}
\widehat{\mathbf{H}}_{\tau}^{(\mathrm{it}+1)}
= 
\widehat{\mathbf{H}}_{\tau}^{(\mathrm{it})}
-
\frac{
\widehat{\mathbf{H}}_{\tau}^{(\mathrm{it})}
\widehat{\boldsymbol{s}}_{\tau}^{(\mathrm{it})}
\bigl(
\widehat{\boldsymbol{s}}_{\tau}^{(\mathrm{it})}
\bigr)^{\top}
\widehat{\mathbf{H}}_{\tau}^{(\mathrm{it})}
}{
\bigl(
\widehat{\boldsymbol{s}}_{\tau}^{(\mathrm{it})}
\bigr)^{\top}
\widehat{\mathbf{H}}_{\tau}^{(\mathrm{it})}
\widehat{\boldsymbol{s}}_{\tau}^{(\mathrm{it})}
}+
\frac{
\widehat{\boldsymbol{\gamma}}_{\tau}^{(\mathrm{it})}
\bigl(
\widehat{\boldsymbol{\gamma}}_{\tau}^{(\mathrm{it})}
\bigr)^{\top}
}{
\bigl(
\widehat{\boldsymbol{\gamma}}_{\tau}^{(\mathrm{it})}
\bigr)^{\top}
\widehat{\boldsymbol{s}}_{\tau}^{(\mathrm{it})}
},
\end{align}
where for initialization (upon termination of the Adam phase) we use the identity $\mathbb{I}_{(K+1)n}$. Observe that the isometry \eqref{isometry property} preserves the curvature products in the BFGS update 
\begin{align}
\bigl(
\widehat{\boldsymbol{\gamma}}_{\tau}^{(\mathrm{it})}
\bigr)^{\top}
\widehat{\boldsymbol{s}}_{\tau}^{(\mathrm{it})}
=
\bigl(
(\mathbf{L}^{\tau})^{\top}\mathbf{\Gamma}_{\tau}^{(\mathrm{it})},
(\mathbf{L}^{\tau})^{\top}\mathbf{S}_{\tau}^{(\mathrm{it})}
\bigr)_{\mathbf{F}}=
\operatorname{trace}
(
(\mathbf{\Gamma}_{\tau}^{(\mathrm{it})})^{\top}
\mathbf{B}^{\tau}
\mathbf{S}_{\tau}^{(\mathrm{it})})
=\langle
\gamma_{\tau}^{(\mathrm{it})},
s_{\tau}^{(\mathrm{it})}
\rangle_{H^{1}(I;\mathbb{R}^{n})}.
\end{align}
Consequently, Euclidean BFGS in the transformed coordinates respects the discrete $H^{1}$ geometry. 

Since $\widehat{\mathbf{H}}_{\tau}^{(\mathrm{it})}$ is positive definite, \eqref{vectorized BFGS direction} yields
\begin{align*}
\bigl(
\widehat{\boldsymbol{g}}_{\tau}^{(\mathrm{it})}
\bigr)^{\top}
\widehat{\boldsymbol{d}}_{\tau}^{(\mathrm{it})}
=
-
\bigl(
\widehat{\boldsymbol{g}}_{\tau}^{(\mathrm{it})}
\bigr)^{\top}
\bigl(
\widehat{\mathbf{H}}_{\tau}^{(\mathrm{it})}
\bigr)^{-1}
\widehat{\boldsymbol{g}}_{\tau}^{(\mathrm{it})}
<0,
\end{align*}
provided that $\widehat{\boldsymbol{g}}_{\tau}^{(\mathrm{it})}\neq 0$, or, equivalently, $\mathfrak{s}_{\tau}(\mathbf{\Theta}_{\tau}^{(\mathrm{it})})>0$. Thus, unless the current iterate is already stationary, $\widehat{\boldsymbol{d}}_{\tau}^{(\mathrm{it})}$ is a descent direction. The step size $\alpha_{\mathrm{it}}>0$ is chosen to satisfy the Wolfe conditions \cite{wright1999numerical}:
\begin{align}
\widetilde{\mathfrak{J}}_{\tau}
\left(
\widehat{\mathbf{\Theta}}_{\tau}^{(\mathrm{it})}
+
\alpha_{\mathrm{it}}
\widehat{\mathbf{D}}_{\tau}^{(\mathrm{it})}
\right)
&\leq
\widetilde{\mathfrak{J}}_{\tau}
\left(
\widehat{\mathbf{\Theta}}_{\tau}^{(\mathrm{it})}
\right)
+
c_{1}\alpha_{\mathrm{it}}
\left(
\widehat{\boldsymbol{g}}_{\tau}^{(\mathrm{it})}
\right)^{\top}
\widehat{\boldsymbol{d}}_{\tau}^{(\mathrm{it})},
\label{eq:Wolfe sufficient decrease}
\\
\left(
\widehat{\boldsymbol{g}}_{\tau}^{(\mathrm{it}+1)}
\right)^{\top}
\widehat{\boldsymbol{d}}_{\tau}^{(\mathrm{it})}
&\geq
c_{2}
\left(
\widehat{\boldsymbol{g}}_{\tau}^{(\mathrm{it})}
\right)^{\top}
\widehat{\boldsymbol{d}}_{\tau}^{(\mathrm{it})}.
\label{eq:Wolfe curvature}
\end{align}
Here, $c_{1}, c_{2} \in (0,1)$ are prescribed line search parameters; we use
$c_{1}=10^{-4}$ and $c_{2}=0.9$. Thus, the first Wolfe condition ensures sufficient decrease in the objective, and the second Wolfe condition guarantees $(
\widehat{\boldsymbol{\gamma}}_{\tau}^{(\mathrm{it})}
)^{\top}
\widehat{\boldsymbol{s}}_{\tau}^{(\mathrm{it})}>0$, preserving positive definiteness of the sequence \eqref{BFGS Hessian update}. Checking \eqref{eq:Wolfe sufficient decrease} requires evaluating the objective, and thus, the state at trial step sizes. Checking \eqref{eq:Wolfe curvature} additionally requires evaluating the gradient of the objective \eqref{modified reduced objective}. Consequently, the Wolfe line search triggers additional solves of \eqref{full discrete optimality sys}. The complete numerical solver  is then given in Algorithm \ref{alg:fixed-depth-optimality-system}. We remark here that beyond BFGS, other quasi-Newton methods, such as that proposed in \cite{newman2021train}, may offer attractive alternatives within our framework. We also point to multilevel optimization techniques such as those developed in \cite{baraldi2025multilevelproximaltrustregionmethod,MR4486512}.

\begin{remark}
Since the stationarity measure $\mathfrak{s}_{\tau}(\mathbf{\Theta}_{\tau})$ need not decrease monotonically in Algorithm \ref{alg:fixed-depth-optimality-system}, one may maintain a checkpoint $(\mathbf\Theta_\tau^{\mathrm{best}},
\mathfrak s_\tau^{\mathrm{best}})$, overwritten whenever
$\mathfrak s_\tau^{(j)}<\mathfrak s_\tau^{\mathrm{best}}$.
This requires only one additional control array storage. If the checkpoint is returned, its state, adjoint, and gradient are recomputed once.
\end{remark}

\begin{algorithm}[Fixed-depth discretized neural ODE solver]
\label{alg:fixed-depth-optimality-system}
\small
\begin{algorithmic}[1]
\Require Grid $\mathcal T_K$, initial nodal control
$\mathbf\Theta_\tau^{(0)}$, training data, regularization parameter
$\lambda$, stationarity tolerance $\varepsilon_{\mathrm{stat}}$,
switching factor $\kappa_{\mathrm{s}}>1$, and Adam, BFGS, and
line-search parameters
\Ensure Control $\mathbf\Theta_\tau$, state $\boldsymbol x_\tau$,
adjoint $\boldsymbol p_\tau$, $H^1$ gradient $\mathbf G_\tau$, and
stationarity measure $\mathfrak s_\tau(\mathbf\Theta_\tau)$

\State Assemble $\mathbf M^\tau$, $\mathbf A^\tau$, and
$\mathbf C^\tau$, set
$\mathbf B^\tau=\mathbf M^\tau+\mathbf A^\tau$, and compute
$\mathbf B^\tau=\mathbf L^\tau(\mathbf L^\tau)^\top$

\State Set $\mathrm{it}\gets0$, initialize the Adam moments
by \eqref{eq: Adam moments}, and set
$\mathbf\Theta_\tau^{(0)}$ as the initial iterate

\State Compute $\boldsymbol x_\tau^{(0)}$, the objective, and
$\boldsymbol p_\tau^{(0)}$; assemble $\mathbf Q_\tau^{(0)}$; solve
\[
\mathbf B^\tau\mathbf G_\tau^{(0)}
=
\lambda\mathbf B^\tau\mathbf\Theta_\tau^{(0)}
+
\mathbf C^\tau\mathbf Q_\tau^{(0)};
\qquad
\mathfrak s_\tau^{(0)}
=
\left[
\operatorname{trace}
\left(
(\mathbf G_\tau^{(0)})^\top
\mathbf B^\tau\mathbf G_\tau^{(0)}
\right)
\right]^{1/2}
\]

\Statex \textit{Adam phase}
\While{
$\mathfrak s_\tau^{(\mathrm{it})}>
\kappa_{\mathrm{s}}\varepsilon_{\mathrm{stat}}$
and the Adam iteration limit has not been reached
}
    \State Update the Adam moments, form
    $\mathbf D_\tau^{(\mathrm{it})}$ by
    \eqref{eq:adam-direction}, and set
    \[
    \mathbf\Theta_\tau^{(\mathrm{it}+1)}
    \gets
    \mathbf\Theta_\tau^{(\mathrm{it})}
    +
    \alpha_{\mathrm{it}}\mathbf D_\tau^{(\mathrm{it})}
    \]

    \State Compute the state
    $\boldsymbol x_\tau^{(\mathrm{it}+1)}$, objective
    $\mathfrak J_\tau(\mathbf\Theta_\tau^{(\mathrm{it}+1)})$, and
    adjoint $\boldsymbol p_\tau^{(\mathrm{it}+1)}$; assemble
    $\mathbf Q_\tau^{(\mathrm{it}+1)}$ and solve
    \[
    \mathbf B^\tau\mathbf G_\tau^{(\mathrm{it}+1)}
    =
    \lambda\mathbf B^\tau\mathbf\Theta_\tau^{(\mathrm{it}+1)}
    +
    \mathbf C^\tau\mathbf Q_\tau^{(\mathrm{it}+1)}
    \]
    according to \eqref{full discrete optimality sys}. Compute $\mathfrak s_\tau^{(\mathrm{it}+1)}$ according to \eqref{stationarity measure}.

    \State Set $\mathrm{it}\gets\mathrm{it}+1$
\EndWhile

\If{
$\varepsilon_{\mathrm{stat}}
<
\mathfrak s_\tau^{(\mathrm{it})}
\leq
\kappa_{\mathrm{s}}\varepsilon_{\mathrm{stat}}$
and BFGS is enabled
}
    \State Set
    \[
    \widehat{\boldsymbol\theta}_\tau^{(\mathrm{it})}
    \gets
    \operatorname{vec}_\tau
    \left(
    (\mathbf L^\tau)^\top\mathbf\Theta_\tau^{(\mathrm{it})}
    \right),
    \qquad
    \widehat{\boldsymbol g}_\tau^{(\mathrm{it})}
    \gets
    \operatorname{vec}_\tau
    \left(
    (\mathbf L^\tau)^\top\mathbf G_\tau^{(\mathrm{it})}
    \right),
    \qquad
    \widehat{\mathbf H}_\tau^{(\mathrm{it})}
    \gets\mathbb{I}_{(K+1)n}
    \]

    \Statex \textit{BFGS phase}
    \While{
    $\mathfrak s_\tau^{(\mathrm{it})}>
    \varepsilon_{\mathrm{stat}}$
    and the BFGS iteration limit has not been reached
    }
        \State Solve
        \[
        \widehat{\mathbf H}_\tau^{(\mathrm{it})}
        \widehat{\boldsymbol d}_\tau^{(\mathrm{it})}
        =
        -\widehat{\boldsymbol g}_\tau^{(\mathrm{it})}
        \]
        according to \eqref{vectorized BFGS direction}, and choose
        $\alpha_{\mathrm{it}}>0$ satisfying the Wolfe conditions
        \eqref{eq:Wolfe sufficient decrease} and \eqref{eq:Wolfe curvature}.

        \State Set
        \[
        \widehat{\boldsymbol\theta}_\tau^{(\mathrm{it}+1)}
        \gets
        \widehat{\boldsymbol\theta}_\tau^{(\mathrm{it})}
        +
        \alpha_{\mathrm{it}}
        \widehat{\boldsymbol d}_\tau^{(\mathrm{it})}
        \]
        and recover $\mathbf\Theta_\tau^{(\mathrm{it}+1)}$ from the triangular solve
        \[
        (\mathbf L^\tau)^\top
        \mathbf\Theta_\tau^{(\mathrm{it}+1)}
        =
        \operatorname{mat}_\tau
        \left(
        \widehat{\boldsymbol\theta}_\tau^{(\mathrm{it}+1)}
        \right)
        \]

        \State Compute 
        $\boldsymbol x_\tau^{(\mathrm{it}+1)}$, 
        $\mathfrak J_\tau(\mathbf\Theta_\tau^{(\mathrm{it}+1)})$,  $\boldsymbol p_\tau^{(\mathrm{it}+1)}$; assemble
        $\mathbf Q_\tau^{(\mathrm{it}+1)}$ and solve $\mathbf{G}_\tau^{(\mathrm{it}+1)}$
        according to \eqref{full discrete optimality sys}. 
        \State Compute $\mathfrak s_\tau^{(\mathrm{it}+1)}$ according to \eqref{stationarity measure} and set
        \[
        \widehat{\boldsymbol g}_\tau^{(\mathrm{it}+1)}
        \gets
        \operatorname{vec}_\tau
        \left(
        (\mathbf L^\tau)^\top
        \mathbf G_\tau^{(\mathrm{it}+1)}
        \right)
        \]

        \State Update
        $\widehat{\mathbf H}_\tau^{(\mathrm{it}+1)}$
        by \eqref{BFGS Hessian update} and set
        $\mathrm{it}\gets\mathrm{it}+1$
    \EndWhile
\EndIf

\State Set
\[
(\mathbf\Theta_\tau,\boldsymbol x_\tau,\boldsymbol p_\tau,
\mathbf G_\tau)
\gets
(\mathbf\Theta_\tau^{(\mathrm{it})},
\boldsymbol x_\tau^{(\mathrm{it})},
\boldsymbol p_\tau^{(\mathrm{it})},
\mathbf G_\tau^{(\mathrm{it})})
\]
\State \Return $\mathbf\Theta_\tau$, $\boldsymbol x_\tau$,
$\boldsymbol p_\tau$, $\mathbf G_\tau$, and
$\mathfrak s_\tau(\mathbf\Theta_\tau)
=\mathfrak s_\tau^{(\mathrm{it})}$.
\end{algorithmic}
\end{algorithm}

\section{The dual-weighted residual error estimation}
\label{sec 6: DRW}
In this section, we derive an a posteriori error bound based on the analysis of the continuous and discrete optimality systems, which will be used for the layerwise adaptive construction of neural network architectures derived from the discretization of neural ODEs. 

First, we note that in view of Proposition \ref{Adjoint regularity}, every solution $(\boldsymbol{x}^\star,\theta^\star, \boldsymbol{p}^\star)$ of the optimality system \eqref{opt system abstract} also solves 
\begin{align}\label{stat cond}
\mathcal{L}^{\prime}_{\mathrm{DG}}(\boldsymbol{x}^{\star}, \theta^{\star};  \boldsymbol{p}^{\star}, \varphi) = 0, \quad \forall{\varphi} \in \mathcal{V} \times \mathcal{U} \times \mathcal{V}.    
\end{align}
Component-wise, \eqref{stat cond} is given by: 
\begin{equation}\label{opt system abstract II}
\begin{aligned}
&(\mathrm{adjoint \  equation}) &&
D_{1}\mathcal{J}(\boldsymbol{x}^{\star}, \theta^{\star}; \varphi_x)
- D_{1}\mathcal{F}_{\text{DG}}(\boldsymbol{x}^{\star}, \theta^{\star}; \boldsymbol{p}^{\star}, \varphi_x)
&&= 0, && \forall \varphi_x \in \mathcal{V}, \\[0.25ex]
&(\mathrm{gradient\ equation}) &&
D_{2}\mathcal{J}(\boldsymbol{x}^{\star}, \theta^{\star}; \varphi_\theta)
- D_{2}\mathcal{F}_{\text{DG}}(\boldsymbol{x}^{\star}, \theta^{\star}; \boldsymbol{p}^{\star}, \varphi_\theta)
&&= 0, && \forall \varphi_\theta \in \mathcal{U}, \\[0.25ex]
& 
(\mathrm{state\ equation}) &&
\mathcal{F}_{\text{DG}}(\boldsymbol{x}^{\star}, \theta^{\star}; \varphi_p)
&&= 0, && \forall \varphi_p \in \mathcal{V}.
\end{aligned}
\end{equation}
Indeed, because $\boldsymbol{x}^{\ast}$ and  $\boldsymbol{p}^{\ast}$ are continuous in time, the jump terms in the above forms vanish. Thanks to this formulation, $\boldsymbol{x}_{\tau}^{\ast} \in \mathcal{W}_{\tau}$ and $\boldsymbol{p}_{\tau}^{\ast} \in \mathcal{V}_{\tau}$ become admissible trial functions in the state and the adjoint equation of \eqref{opt system abstract II}, respectively, even though $\mathcal{W}_{\tau} \not\subset \mathcal{W}$ and $\mathcal{V}_{\tau} \not\subset \mathcal{W}$. 

We further observe that due to the inclusions $\mathcal{W}_{\tau} \subset \mathcal{V}, \ \mathcal{U}_{\tau} \subset \mathcal{U}$ and $\mathcal{V}_{\tau} \subset \mathcal{V}$, the error 
$e = (\boldsymbol{x}^{\star} - \boldsymbol{x}_{\tau}^{\star},\, \theta^{\star} - \theta_{\tau}^{\star},\, \boldsymbol{p}^{\star} - \boldsymbol{p}_{\tau}^{\star}) 
\in \mathcal{V} \times \mathcal{U} \times \mathcal{V}$ 
is a feasible test function in \eqref{opt system abstract II}. 
Substituting it into \eqref{opt system abstract II} yields the Galerkin orthogonality condition
\begin{align}\label{MH:GO}
\mathcal{L}^{\prime}_{\mathrm{DG}}(\boldsymbol{x}^{\star}, \theta^{\star}; \boldsymbol{p}^{\star}, e) = 0.
\end{align}
With the continuous optimality system \eqref{opt system abstract II} and its discrete counterpart \eqref{opt system discrete} at our disposal, we adapt \cite[Theorem 2.2]{kraft2010dual} to the neural ODE setting and provide an a posteriori representation formula for the error in the functional $\mathcal{J}$.  
\begin{proposition}\label{Prop:error representation formula} Suppose that $\boldsymbol{\sigma}^{\prime \prime \prime} \in L^{\infty}(\mathbb{R}^{d})$.  Let $(\boldsymbol{x}^{\star}, \theta^{\star}, \boldsymbol{p}^{\star}) \in \mathcal{W} \times \mathcal{U} \times \mathcal{W}$ and $(\boldsymbol{x}_{\tau}^{\star}, \theta_\tau^{\star}, \boldsymbol{p}_{\tau}^{\star})\in \mathcal{W}_\tau \times \mathcal{U}_\tau \times \mathcal{V}_\tau$ be solutions of \eqref{opt system abstract II} and \eqref{opt system discrete}, respectively. Then 
\begin{align}\label{error representation formula}
\mathcal{J}(\boldsymbol{x}^{\star}, \theta^{\star}) - \mathcal{J}(\boldsymbol{x}_\tau^{\star}, \theta_\tau^{\star}) = \frac{1}{2} \rho_{\mathrm{adj}}[ \boldsymbol{x}^{\star} - \tilde{\boldsymbol{x}}_{\tau}] + \frac{1}{2}\rho_{\mathrm{grad}}[ \theta^{\star} - \tilde{\theta}_{\tau}] +  \frac{1}{2}  \rho_{\mathrm{state}}[ \boldsymbol{p}^{\star} - \tilde{\boldsymbol{p}}_{\tau}]   + R.
\end{align}
The residuals $\rho_{\mathrm{adj}} \in \mathcal{V}^{\ast}$, $\rho_{\mathrm{grad}} \in \mathcal{U}^{\ast}$ and $\rho_{\mathrm{state}} \in \mathcal{V}^{\ast}$ are given by
\begin{equation}\label{residual formulas}
\begin{aligned}
\rho_{\mathrm{adj}}[ \boldsymbol{x}^{\star} - \tilde{\boldsymbol{x}}_{\tau}]
&= \sum_{k=1}^{K} \int_{I_{k}} 
    \big( 
        - \dot{\boldsymbol{p}}_{\tau}^{\star} 
        - D_{1}N_{F}(\boldsymbol{x}_{\tau}^{\star}, \theta_\tau^{\star})^{\ast}\boldsymbol{p}_{\tau}^{\star},\ 
        \boldsymbol{x}^{\star} - \tilde{\boldsymbol{x}}_{\tau} 
    \big) \, dt - \sum_{k=1}^{K} 
    \big( 
        \jump{\boldsymbol{p}_{\tau}^{\star}}^{k-1},\ 
        \boldsymbol{x}^{\star}(t_{k-1}^{+}) - \tilde{\boldsymbol{x}}_{\tau}(t_{k-1}^{+}) 
    \big), \\[0.5em]
\rho_{\mathrm{grad}}[ \theta^{\star} - \tilde{\theta}_{\tau}]
&= \sum_{k=1}^{K} \int_{I_{k}}  
        \lambda \Big[
            \big(\theta_{\tau}^{\star},  \theta^{\star} - \tilde{\theta}_{\tau}\big) 
            + \big(\dot{\theta}_{\tau}^{\star}, \dot{\theta}^{\star} - \dot{\tilde{\theta}}_{\tau}\big)
        \Big]
        + \big(
            D_{2}N_{F}(\boldsymbol{x}_{\tau}^{\star}, \theta_\tau^{\star})^{\ast}\boldsymbol{p}_{\tau}^{\star},\ 
            \theta^{\star} - \tilde{\theta}_{\tau} 
        \big)
     \, dt, \\[0.5em]
\rho_{\mathrm{state}}[ \boldsymbol{p}^{\star} - \tilde{\boldsymbol{p}}_{\tau}] 
&= \sum_{k=1}^{K} \int_{I_{k}} 
    \big( 
        \dot{\boldsymbol{x}}_{\tau}^{\star} - N_{F}(\boldsymbol{x}_{\tau}^{\star}, \theta_{\tau}^{\star}),\ 
        \boldsymbol{p}^{\star} - \tilde{\boldsymbol{p}}_{\tau}
    \big) \, dt  
    + \sum_{k=1}^{K} 
    \big( 
        \jump{\boldsymbol{x}_{\tau}^{\star}}^{k},\ 
        \boldsymbol{p}^{\star}(t_{k}^{-}) - \tilde{\boldsymbol{p}}_{\tau}(t_{k}^{-}) 
    \big), 
\end{aligned}
\end{equation}
where $\tilde{\boldsymbol{x}}_{\tau} \in \mathcal{W}_{\tau}$, $\tilde{\boldsymbol{p}}_{\tau}\in \mathcal{V}_{\tau}$ and $\tilde{\theta}_{\tau} \in \mathcal{U}_{\tau}$ are arbitrary, and $R$ is the remainder
\begin{equation}\label{remainder}
\begin{aligned}
R = \frac{1}{2}& \int_{0}^{1} \Big(\mathcal{J}^{\prime \prime \prime }(\boldsymbol{x}_{\tau}^{\star} + se_{\boldsymbol{x}}, \theta_{\tau}^{\star} + se_{\theta};  e, e, e)  - \mathcal{F}^{\prime \prime \prime }_{\mathrm{DG}}(\boldsymbol{x}_{\tau}^{\star} + se_{\boldsymbol{x}}, \theta_{\tau}^{\star} + se_{\theta}; \boldsymbol{p}_{\tau}^{\star} + se_{\boldsymbol{p}}, e, e, e) \Big)s(s-1) \ ds,
\end{aligned}
\end{equation}
which is cubic in the error $e = (\boldsymbol{x}^{\star} - \boldsymbol{x}_{\tau}^{\star},\, \theta^{\star} - \theta_{\tau}^{\star},\, \boldsymbol{p}^{\star} - \boldsymbol{p}_{\tau}^{\star}) 
\in \mathcal{V} \times \mathcal{U} \times \mathcal{V}$.
\end{proposition}
\begin{proof} Using the third equation in \eqref{opt system abstract II} and the third equation in \eqref{opt system discrete}, we obtain
\begin{align*}
\mathcal{J}(\boldsymbol{x}^{\star}, \theta^{\star}) - \mathcal{J}(\boldsymbol{x}_\tau^{\star}, \theta_\tau^{\star}) & =  \mathcal{L}_{\mathrm{DG}}(\boldsymbol{x}^{\star}, \theta^{\star}; \boldsymbol{p}^{\star}) + \mathcal{F}_{\text{DG}}(\boldsymbol{x}^{\star}, \theta^{\star}; \boldsymbol{p}^{\star}) - \mathcal{L}_{\mathrm{DG}}(\boldsymbol{x}_\tau^{\star}, \theta_\tau^{\star}; \boldsymbol{p}_\tau^{\star}) - \mathcal{F}_{\text{DG}}(\boldsymbol{x}_\tau^{\star}, \theta_\tau^{\star}; \boldsymbol{p}_\tau^{\star}) \\[0.5em]
& = \mathcal{L}_{\mathrm{DG}}(\boldsymbol{x}^{\star}, \theta^{\star}; \boldsymbol{p}^{\star}) - \mathcal{L}_{\mathrm{DG}}(\boldsymbol{x}_\tau^{\star}, \theta_\tau^{\star}; \boldsymbol{p}_\tau^{\star}) \\[0.5em] 
& = \int_{0}^{1} \mathcal{L}^{\prime}_{\mathrm{DG}}(\boldsymbol{x}_{\tau}^{\star} + se_{\boldsymbol{x}}, \theta_{\tau}^{\star} + se_{\theta}; \boldsymbol{p}_{\tau}^{\star} + se_{\boldsymbol{p}}, e) \ ds + \frac{1}{2} \mathcal{L}^{\prime}_{\mathrm{DG}}(\boldsymbol{x}_{\tau}^{\star}, \theta_{\tau}^{\star}; \boldsymbol{p}_{\tau}^{\star}, e) \\[0.5em]
& -  \frac{1}{2} \mathcal{L}^{\prime}_{\mathrm{DG}}(\boldsymbol{x}_{\tau}^{\star}, \theta_{\tau}^{\star}; \boldsymbol{p}_{\tau}^{\star}, e) - \frac{1}{2} \mathcal{L}^{\prime}_{\mathrm{DG}}(\boldsymbol{x}^{\star}, \theta^{\star}; \boldsymbol{p}^{\star}, e), 
\end{align*}
where the last term is zero in view of \eqref{MH:GO}, and the first and the last two terms in the above expression form the remainder of the trapezoidal quadrature rule:
\begin{align*}
R = \int_{0}^{1} \mathcal{L}^{\prime}_{\mathrm{DG}}(\boldsymbol{x}_{\tau}^{\star} + se_{\boldsymbol{x}}, \theta_{\tau}^{\star} + se_{\theta}; \boldsymbol{p}_{\tau}^{\star} + se_{\boldsymbol{p}}, e) \ ds - \frac{1}{2}\Big[\mathcal{L}^{\prime}_{\mathrm{DG}}(\boldsymbol{x}_{\tau}^{\star}, \theta_{\tau}^{\star}; \boldsymbol{p}_{\tau}^{\star}, e) +  \mathcal{L}^{\prime}_{\mathrm{DG}}(\boldsymbol{x}^{\star}, \theta^{\star}; \boldsymbol{p}^{\star}, e) \Big]. 
\end{align*}
Therefore, we get the following equality
\begin{equation}\label{Difference quad}
\begin{aligned}
\mathcal{J}(\boldsymbol{x}^{\star}, \theta^{\star}) - \mathcal{J}(\boldsymbol{x}_\tau^{\star}, \theta_\tau^{\star}) & = \frac{1}{2} \mathcal{L}^{\prime}_{\mathrm{DG}}(\boldsymbol{x}_{\tau}^{\star}, \theta_{\tau}^{\star}; \boldsymbol{p}_{\tau}^{\star}, \boldsymbol{x}^{\star}- \boldsymbol{x}_\tau^{\star}, \theta^{\star} - \theta_\tau^{\star}, \boldsymbol{p}^{\star} - \boldsymbol{p}_\tau^{\star}) + R \\[0.5em] 
& = \frac{1}{2} \mathcal{L}^{\prime}_{\mathrm{DG}}(\boldsymbol{x}_{\tau}^{\star}, \theta_{\tau}^{\star}; \boldsymbol{p}_{\tau}^{\star}, \boldsymbol{x}^{\star} - \tilde{\boldsymbol{x}}_\tau, \theta^{\star} - \tilde{\theta}_\tau, \boldsymbol{p}^{\star} - \tilde{\boldsymbol{p}}_\tau) + R, 
\end{aligned}
\end{equation}
where we replaced $(\boldsymbol{x}_\tau^{\star}, \theta_\tau^{\star}, \boldsymbol{p}_\tau^{\star})$ by an arbitrary $(\tilde{\boldsymbol{x}}_\tau, \tilde{\theta}_\tau, \tilde{\boldsymbol{p}}_\tau) \in \mathcal{W}_{\tau} \times \mathcal{U}_{\tau} \times \mathcal{V}_{\tau}$ in view of \eqref{discrete opt cond}. Writing \eqref{Difference quad} in terms of individual components yields 
\begin{align*}
\mathcal{J}(\boldsymbol{x}^{\star}, \theta^{\star}) - \mathcal{J}(\boldsymbol{x}_\tau^{\star}, \theta_\tau^{\star}) & = \frac{1}{2} \Big[ D_{1}\mathcal{J}(\boldsymbol{x}_{\tau}^{\star}, \theta_{\tau}^{\star}; \boldsymbol{x}^{\star} - \tilde{\boldsymbol{x}}_\tau) 
- D_{1}\mathcal{F}_{\text{DG}}(\boldsymbol{x}_{\tau}^{\star}, \theta_{\tau}^{\star}; \boldsymbol{p}_{\tau}^{\star}, \boldsymbol{x}^{\star} - \tilde{\boldsymbol{x}}_\tau) 
 \Big] \\[0.5ex]
& + \frac{1}{2} \Big[D_{2}\mathcal{J}(\boldsymbol{x}_{\tau}^{\star}, \theta_\tau^{\star}; \theta^{\star} - \tilde{\theta}_\tau) 
- D_{2}\mathcal{F}_{\text{DG}}(\boldsymbol{x}_{\tau}^{\star}, \theta_\tau^{\star}; \boldsymbol{p}_{\tau}^{\star}, \theta^{\star} - \tilde{\theta}_\tau) \Big]  \\[0.5ex]
& + \frac{1}{2} \mathcal{F}_{\text{DG}}(\boldsymbol{x}_{\tau}^{\star}, \theta_\tau^{\star}; \boldsymbol{p}^{\star} - \tilde{\boldsymbol{p}}_\tau) + R  \\
&= \frac{1}{2} \rho_{\mathrm{adj}}[ \boldsymbol{x}^{\star} - \tilde{\boldsymbol{x}}_{\tau}]   + \frac{1}{2} \rho_{\mathrm{grad}}[ \theta^{\star} - \tilde{\theta}_{\tau}]  + \frac{1}{2} \rho_{\mathrm{state}}[ \boldsymbol{p}^{\star} - \tilde{\boldsymbol{p}}_{\tau}]  +R.    
\end{align*}
Expressing the remainder $R$ using the Peano kernel formula for the trapezoidal quadrature rule yields
\begin{align*}
R= \frac{1}{2} \int_{0}^{1}\mathcal{L}^{\prime \prime \prime }_{\mathrm{DG}}(\boldsymbol{x}_{\tau}^{\star} + se_{\boldsymbol{x}}, \theta_{\tau}^{\star} + se_{\theta}; \boldsymbol{p}_{\tau}^{\star} + se_{\boldsymbol{p}}, e, e, e)  s(s-1)  \ ds.
\end{align*}
Expanding the trilinear form $\mathcal{L}^{\prime \prime \prime }_{\mathrm{DG}}(\boldsymbol{x}_{\tau}^{\star} + se_{\boldsymbol{x}}, \theta_{\tau}^{\star} + se_{\theta}; \boldsymbol{p}_{\tau}^{\star} + se_{\boldsymbol{p}}, \cdot, \cdot, \cdot)$ in the above integrand in terms of its components yields the remainder \eqref{remainder}. The assumption $\boldsymbol{\sigma}^{\prime \prime \prime} \in L^{\infty}(\mathbb{R}^{d})$ then guarantees that the remainder remains bounded. 
\end{proof}
%\begin{remark}
%Note that it is in principle possible to not discretize the control space in %\eqref{opt system discrete}; cf. \cite{hinze2005variational}. 
%\end{remark}
%For our purposes, we use the following interval norm 
%\begin{align}\label{interval norm}
%\lVert v \rVert_{I_{k}} = \underset{t \in I_{k}}{\sup} \ \lVert v(t) \rVert.
%\end{align}
We note that Proposition~\ref{Prop:error representation formula} requires sufficiently regular activation functions. For example, we use the smooth hyperbolic tangent in our numerical examples. Less regular activation functions can, in principle, be approximated (or smoothed) to meet the required assumptions; see, e.g., \cite{dong2024descent}.

\section{Layerwise adaptive algorithm} 
\label{sec 7: Layerwise adaptive algorithm}

In this section, we present a depth-adaptive algorithm for training residual neural networks (discretized neural ODEs) and provide details of its numerical realization. The overall methodology follows the paradigm
\begin{equation*}
\mathrm{OPTIMIZE} \quad \Longrightarrow \quad \mathrm{ESTIMATE} \quad \Longrightarrow \quad \mathrm{MARK} \quad \Longrightarrow \quad \mathrm{REFINE}.
\end{equation*}
Starting from an optimized fixed-depth architecture, we refine it by inserting layers at grid locations selected using indicators derived from the error estimate \eqref{DRW error bound}, followed by re-optimization. This procedure is repeated until a prescribed target or stopping criterion is satisfied, or the layer budget is exhausted.

We first describe the construction of a computable error indicator which guides the adaptive process. The error representation formula \eqref{error representation formula} and the triangle inequality give us
\begin{align}\label{DRW error bound}
|\mathcal{J}(\boldsymbol{x}^{\star},\theta^{\star})-\mathcal{J}(\boldsymbol{x}_{\tau}^{\star},\theta_{\tau}^{\star})|\leq\frac{1}{2}\sum_{k=1}^{K}\bigg[\big|\rho_{\mathrm{adj}}^{k}[\boldsymbol{x}^{\star}-\tilde{\boldsymbol{x}}_{\tau}] \big|+\big|\rho_{\mathrm{grad}}^{k}[\theta^{\star}-\tilde{\theta}_{\tau}]\big|+\big|\rho_{\mathrm{state}}^{k}[\boldsymbol{p}^{\star}-\tilde{\boldsymbol{p}}_{\tau}]\big|\bigg]+|R|.
\end{align}
where the interval-wise residuals are given by
\begin{equation}\label{interval-wise residuals}
\begin{aligned}
\rho_{\mathrm{adj}}^{k}[\boldsymbol{x}^{\star}-\tilde{\boldsymbol{x}}_{\tau}] & \; =   \int_{I_{k}} 
    \big( 
        - \dot{\boldsymbol{p}}_{\tau}^{\star} 
        - D_{1}N_{F}(\boldsymbol{x}_{\tau}^{\star}, \theta_\tau^{\star})^{\ast}\boldsymbol{p}_{\tau}^{\star},\ 
        \boldsymbol{x}^{\star} - \tilde{\boldsymbol{x}}_{\tau} 
    \big) \, dt -  
    \big( 
        \jump{\boldsymbol{p}_{\tau}^{\star}}^{k-1},\ 
        \boldsymbol{x}^{\star}(t_{k-1}^{+}) - \tilde{\boldsymbol{x}}_{\tau}(t_{k-1}^{+}) 
    \big), \\[1.5ex] 
\rho_{\mathrm{grad}}^{k}[\theta^{\star}-\tilde{\theta}_{\tau}] & \; =\int_{I_{k}}  
        \lambda \Big[
            \big(\theta_{\tau}^{\star},  \theta^{\star} - \tilde{\theta}_{\tau}\big) 
            + \big(\dot{\theta}_{\tau}^{\star}, \dot{\theta}^{\star} - \dot{\tilde{\theta}}_{\tau}\big)
        \Big]
        + \big(
            D_{2}N_{F}(\boldsymbol{x}_{\tau}^{\star}, \theta_\tau^{\star})^{\ast}\boldsymbol{p}_{\tau}^{\star}, 
            \theta^{\star} - \tilde{\theta}_{\tau} 
        \big)
     \, dt,   \\[1.5ex]
\rho_{\mathrm{state}}^{k}[\boldsymbol{p}^{\star}-\tilde{\boldsymbol{p}}_{\tau}] & \; =  \int_{I_{k}} 
    \big( 
        \dot{\boldsymbol{x}}_{\tau}^{\star} - N_{F}(\boldsymbol{x}_{\tau}^{\star}, \theta_{\tau}^{\star}),\ 
        \boldsymbol{p}^{\star} - \tilde{\boldsymbol{p}}_{\tau}
    \big) \, dt  
    + 
    \big( 
        \jump{\boldsymbol{x}_{\tau}^{\star}}^{k},\ 
        \boldsymbol{p}^{\star}(t_{k}^{-}) - \tilde{\boldsymbol{p}}_{\tau}(t_{k}^{-}) 
        \big). 
\end{aligned}
\end{equation}
Neglecting the remainder in \eqref{DRW error bound}, the residuals \eqref{interval-wise residuals} furnish an error indicator
\begin{align}\label{error indicator}
\triangle = \frac{1}{2}\sum_{k=1}^{K} \Big( \big|\rho_{\mathrm{adj}}^{k}[\boldsymbol{x}^{\star}-\tilde{\boldsymbol{x}}_{\tau}]\big| + \big|\rho_{\mathrm{grad}}^{k}[\theta^{\star}-\tilde{\theta}_{\tau}]\big| + \big|\rho_{\mathrm{state}}^{k}[\boldsymbol{p}^{\star}-\tilde{\boldsymbol{p}}_{\tau}]\big| \Big) =: \frac{1}{2}\sum_{k=1}^{K} \triangle_{k}.
\end{align}
The above indicator is not yet computable: it involves the exact triplet $(\boldsymbol{x}^{\star}, \theta^{\star}, \boldsymbol{p}^{\star})$ and the discrete stationary point $(\boldsymbol{x}_{\tau}^{\star}, \theta^{\star}_{\tau}, \boldsymbol{p}^{\star}_\tau)$ of \eqref{opt system discrete}. We now specify the discretization that turns each interval-wise residual into a computable expression.

In what follows, $(\boldsymbol{x}_{\tau}^{\star}, \theta^{\star}_{\tau}, \boldsymbol{p}^{\star}_\tau)$ is replaced by its approximation $(\boldsymbol{x}_{\tau}, \theta_{\tau}, \boldsymbol{p}_{\tau})$ obtained from \eqref{full discrete optimality sys} and $(\tilde{\boldsymbol{x}}_{\tau}, \tilde{\theta}_{\tau}, \tilde{\boldsymbol{p}}_{\tau}):=(\boldsymbol{x}_{\tau}, \theta_{\tau}, \boldsymbol{p}_{\tau})$. The triplet $(\boldsymbol{x}^{\star}, \theta^{\star}, \boldsymbol{p}^{\star})$ is not available and must be replaced for computing purposes by a suitable reconstruction obtained from $(\boldsymbol{x}_{\tau}, \theta_{\tau}, \boldsymbol{p}_{\tau})$.  We recall that $\boldsymbol{x}_{\tau}(t) = \boldsymbol{x}_{\tau}^{k-1}$ for $t \in [t_{k-1}, t_{k})$ and $\boldsymbol{p}_{\tau}(t) = \boldsymbol{p}_{\tau}^{k}$ for $t \in (t_{k-1}, t_{k}]$. We use $s = \frac{t-t_{k-1}}{h_{k}} \in [0,1]$ as the reference coordinate on $[t_{k-1}, t_{k}]$, and we abbreviate
\begin{align*}
\Delta \boldsymbol{x}_{\tau}^{k} = \boldsymbol{x}_{\tau}^{k} - \boldsymbol{x}_{\tau}^{k-1}, \quad \Delta \boldsymbol{p}_{\tau}^{k} = \boldsymbol{p}_{\tau}^{k} - \boldsymbol{p}_{\tau}^{k-1},  \quad \Delta \theta_{\tau}^{k} = \theta_{\tau}^{k} - \theta_{\tau}^{k-1},
\end{align*}
The discrete control is available and is piecewise-linear: $\restr{\theta_{\tau}}{I_{k}}(s)= \theta_{\tau}^{k-1} + s \Delta \theta_{\tau}^{k}$. Then, we replace the unavailable $\boldsymbol{x}^{\star}$ and $\boldsymbol{p}^{\star}$ by the continuous, piecewise-linear reconstructions $\mathcal{Q}_{{\tau}}^{\boldsymbol{x}^{\star}}[\boldsymbol{x}_{\tau}] \in \mathcal{W}$ and $\mathcal{Q}_{{\tau}}^{\boldsymbol{p}^{\star}}[\boldsymbol{p}_{\tau}]\in \mathcal{W}$, which are defined interval-wise as follows: 
\begin{align}\label{reconstructions}
\mathcal{Q}_{{\tau}}^{\boldsymbol{x}^{\star}}[\boldsymbol{x}_{\tau}](s) = \boldsymbol{x}_{\tau}^{k-1} + s \Delta \boldsymbol{x}_{\tau}^{k} \quad \text{and} \quad \mathcal{Q}_{{\tau}}^{\boldsymbol{p}^{\star}}[\boldsymbol{p}_{\tau}](s) = \boldsymbol{p}_{\tau}^{k-1} + s \Delta \boldsymbol{p}_{\tau}^{k}.
\end{align}
Observe that $\mathcal{Q}_{{\tau}}^{\boldsymbol{x}^{\star}}[\boldsymbol{x}_{\tau}](s) - \boldsymbol{x}_{\tau}^{k-1} = s \Delta \boldsymbol{x}_{\tau}^{k}$ and $\mathcal{Q}_{{\tau}}^{\boldsymbol{p}^{\star}}[\boldsymbol{p}_{\tau}](s) - \boldsymbol{p}_{\tau}^{k} = -(1-s) \Delta \boldsymbol{p}_{\tau}^{k}$ on $[t_{k-1}, t_{k}]$. The first difference vanishes at $t_{k-1}^{+}$ and the second at $t_{k}^{-}$; hence, the corresponding jump pairings in \eqref{interval-wise residuals} vanish. Plugging \eqref{reconstructions} into the state and the adjoint residuals in \eqref{interval-wise residuals} and transforming $I_{k}$ to $[0,1]$, we get
\begin{equation}
\begin{aligned}
\widetilde{\rho}_{\mathrm{adj}}^{\,k}\big[\mathcal{Q}_{{\tau}}^{\boldsymbol{x}^{\star}}[\boldsymbol{x}_{\tau}]-\boldsymbol{x}_{\tau}\big] \; & =   -h_{k}\int_{0}^{1} 
    \big( 
        D_{1}N_{F}(\boldsymbol{x}_{\tau}^{k-1}, \theta_\tau(s))^{\ast}\boldsymbol{p}_{\tau}^{k},\ 
        s \Delta \boldsymbol{x}_{\tau}^{k} 
    \big) \, ds , \\ 
\widetilde{\rho}_{\mathrm{state}}^{\,k}\big[\mathcal{Q}_{{\tau}}^{\boldsymbol{p}^{\star}}[\boldsymbol{p}_{\tau}]-\boldsymbol{p}_{\tau}\big] \; & = h_{k} \int_{0}^{1} 
    \big( N_{F}(\boldsymbol{x}_{\tau}^{k-1}, \theta_{\tau}(s)),\ (1-s) \Delta \boldsymbol{p}_{\tau}^{k} 
    \big) \, ds.   
\end{aligned}
\end{equation}
The above integrals are evaluated using Gauss--Legendre quadrature on $[0,1]$, starting with order $q=8$. The order is doubled until all local state and
adjoint residuals agree within the combined absolute and relative tolerances
$10^{-9}$ and $10^{-7}$, respectively, or until $q=128$ is reached.

For the control residual, $\theta^{\star}$ is replaced by the globally continuous piecewise-quadratic reconstruction $\mathcal{Q}_{\tau}^{\theta^{\star}}[\theta_{\tau}] \in \mathcal{U}$, which is defined in the reference coordinates on each interval $[t_{k-1}, t_k]$ by
\begin{align*}
\mathcal{Q}_{\tau}^{\theta^{\star}}[\theta_{\tau}](s)=A_{k}s^{2} +B_{k}s + C_{k}.
\end{align*}
The coefficients $A_k, B_k, C_k \in \mathbb{R}^n$ are determined by the conditions $\mathcal{Q}_{\tau}^{\theta^{\star}}[\theta_{\tau}](0) = \theta_\tau^{k-1}$, $
\mathcal{Q}_{\tau}^{\theta^{\star}}[\theta_{\tau}](1) = \theta_\tau^{k}$ and
$\dot{\mathcal{Q}}_{\tau}^{\theta^{\star}}[\theta_{\tau}](0) = S_{k-1}$, where $S_{k-1}$ is the left-sided slope, which is given by $S_0 = 0$ and $S_k = \frac{\Delta \theta_{\tau}^{k}}{h_k}$. This gives explicitly $A_k = \Delta \theta_{\tau}^{k} - B_k$, $B_k = S_{k-1} \, h_k,$ and $C_k = \theta_\tau^{k-1}$. 

\begin{figure}[h]
    \centering
    \includegraphics[width=1\linewidth]{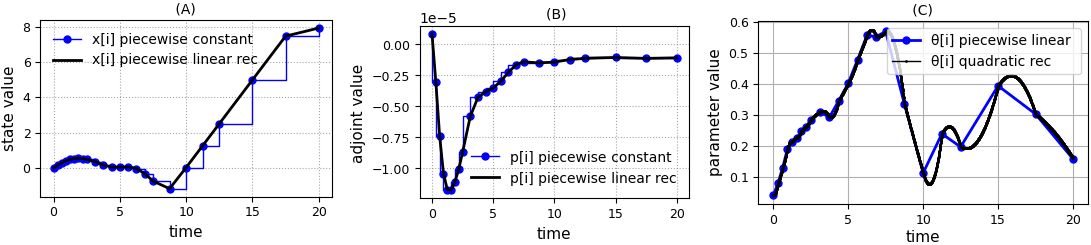}
    \caption{Discrete components of the state, control and the adjoint, and their respective reconstructions. (\textbf{A}): $x_{\tau,j}^{i}$ ($j$-th component of the $i$-th equation) of the discrete state $\boldsymbol{x}_{\tau}$ and its piecewise linear reconstruction. (\textbf{B}): component $p_{\tau, j}^{i}$ of the discrete adjoint $\boldsymbol{p}_\tau$ and its piecewise linear reconstruction. (\textbf{C}): component $\theta^{i}$ of the discrete control $\theta_{\tau}$ and its piecewise quadratic reconstruction }
    \label{fig:reconstruction}
\end{figure}

Figure~\ref{fig:reconstruction} shows piecewise linear reconstructions from the (piecewise constant) state approximation $\boldsymbol{x}_{\tau}$ and adjoint state approximation $\boldsymbol{p}_{\tau}$, as well as a piecewise quadratic reconstruction from the (piecewise linear) $\theta_{\tau}$, taken from our Swiss Roll  example with $T=20$ (see the next section).
To obtain a computable control residual, we replace $\theta^{\star}$ by
$\mathcal{Q}_{\tau}^{\theta^{\star}}[\theta_{\tau}]$, $\theta_{\tau}^{\star}$ by $\theta_\tau$ and use $\tilde{\theta}_{\tau}:=\theta_{\tau}$ in
\eqref{interval-wise residuals}. Transforming $I_k$ to $[0,1]$ gives
\begin{align}\label{grad res approx}
\widetilde{\rho}_{\mathrm{grad}}^{\,k}\big[\mathcal{Q}_{\tau}^{\theta^{\star}}[\theta_{\tau}]-\theta_\tau\big]
=
h_k\int_0^1
\lambda\Big[
\big(\theta_\tau,\mathcal{Q}_{\tau}^{\theta^{\star}}[\theta_{\tau}]-\theta_\tau\big)
+
\big(\dot{\theta}_\tau,
\dot{\mathcal{Q}}_{\tau}^{\theta^{\star}}[\theta_{\tau}]-\dot{\theta}_\tau\big)
\Big]
+
\big(
\tilde q_\tau^{k},
\mathcal{Q}_{\tau}^{\theta^{\star}}[\theta_{\tau}]-\theta_\tau
\big)
\,ds.
\end{align}
where $\tilde q_\tau^{k}(s)= D_2N_F\bigl(
\boldsymbol{x}_\tau^{k-1},
\theta_\tau(s)
\bigr)^\ast
\boldsymbol{p}_\tau^k$. The integral containing $\theta_\tau$ and $\mathcal{Q}_{\tau}^{\theta^{\star}}[\theta_{\tau}]$ is computed as follows: 
\begin{align*}
\int_{I_k} \big( \theta_\tau, \mathcal{Q}_{\tau}^{\theta^{\star}}[\theta_{\tau}] - \theta_{\tau} \big) \, dt
 = h_k  
\big(\theta^{k-1}_\tau, \alpha^k \big) +h_{k} \big(\theta^k_\tau, \beta^k \big) - \begin{pmatrix}\theta_\tau^{k-1} \\ \theta_\tau^{k}\end{pmatrix}^{\!\!\top}
(\mathbf{M}_k^{\tau} \otimes \mathbb{I}_{n})
\begin{pmatrix}\theta_\tau^{k-1} \\ \theta_\tau^{k}\end{pmatrix},
\end{align*}
where the local mass matrix $M_k^{\tau}$, and the coefficients $\alpha^{k}, \beta^{k} \in \mathbb{R}^{n}$ are given by
\begin{align*}
 \mathbf{M}_k^{\tau} =
\begin{pmatrix}
\frac{h_{k}}{3} &  \frac{h_{k}}{6} \\
 \frac{h_{k}}{6} & \frac{h_{k}}{3}
\end{pmatrix} , \quad \alpha^{k} = \frac{A_k}{12} + \frac{B_k}{6} + \frac{C_k}{2}, 
\quad
\beta^{k} = \frac{A_k}{4} + \frac{B_k}{3} + \frac{C_k}{2}.
\end{align*} 
For the derivatives $\dot{\theta}_\tau$ and $\dot{\mathcal{Q}}_{\tau}^{\theta^{\star}}[\theta_{\tau}]$, the respective integral is given by
\begin{align*}
\int_{I_k}
\bigl(
\dot{\theta}_\tau,
\dot{\mathcal{Q}}_{\tau}^{\theta^{\star}}[\theta_{\tau}]-\dot{\theta}_\tau
\bigr)\,dt
=
\frac{1}{h_k}
\int_0^1
\bigl(
\Delta\theta_\tau^k,
2A_ks+B_k-\Delta\theta_\tau^k
\bigr)\,ds
=
\frac{1}{h_k}
\bigl(
\Delta\theta_\tau^k,
A_k+B_k-\Delta\theta_\tau^k
\bigr)
=0.
\end{align*}
The control pullback term in \eqref{grad res approx} is evaluated using Gauss--Legendre quadrature on $[0,1]$, as described above for the state and the adjoint residuals. Therefore, \eqref{grad res approx} can now be evaluated. We define
\begin{align}\label{error indicator approximation}
\eta_{k}: = \big|\widetilde{\rho}_{\mathrm{adj}}^{\, k}\big[\mathcal{Q}_{\tau}^{\boldsymbol{x}^{\star}}[\boldsymbol{x}_{\tau}]-\boldsymbol{x}_{\tau}\big]\big| + \big|\widetilde{\rho}_{\mathrm{grad}}^{\, k}\big[\mathcal{Q}_{\tau}^{\theta^{\star}}[\theta_{\tau}]-\theta_\tau\big]\big| + \big|\widetilde{\rho}_{\mathrm{state}}^{\, k}\big[\mathcal{Q}_{\tau}^{\boldsymbol{p}^{\star}}[\boldsymbol{p}_{\tau}]-\boldsymbol{p}_{\tau}\big]\big|,
\end{align}
which provides a computable approximation of the local indicator
$\triangle_k$ in \eqref{error indicator}. We further define $\eta_K:=\frac{1}{2}\sum_{k=1}^{K}\eta_k$ as an approximation of $\triangle$.

Equipped with the computable indicators \eqref{error indicator approximation}, we define a marking strategy $\operatorname{Mark}(\mathcal{T}_K,\{\eta_k\}_{k=1}^{K})$ that selects a set of intervals $\mathcal{M}_{K}\subset\{1,\ldots,K\}$ for layer insertion. We begin with the conservative strategy which marks the interval $I_{k^\star}$ with the largest local indicator:
\begin{align}\label{conservative refinement}
k^\star = \operatorname*{arg\,max}_{1\leq k\leq K} \, \eta_k.
\end{align}
Hence, $\mathcal{M}_{K}=\{k^\star\}$ with $|\mathcal{M}_{K}|=1$. However, several intervals may be marked simultaneously. Under D\"orfler-type marking \cite{doerfler1996convergent}, we select the smallest marked set $\mathcal{M}_{K}$ satisfying
\begin{align}
\sum_{k\in \mathcal{M}_{K}}\eta_{k}
\geq
\vartheta_{\mathcal{D}}\sum_{k=1}^{K}\eta_{k},
\end{align}
where $0<\vartheta_{\mathcal{D}} < 1$ is the fraction of the total indicator mass to be marked. The achieved bulk fraction
\begin{align}\label{empirical bulk}
\beta_{K}^{\vartheta_{\mathcal{D}}} = \frac{\sum_{k \in\mathcal{M}_{K} }\eta_{k}}{\sum_{k=1}^{K} \eta_{k}}
\end{align}
satisfies $\beta_{K}^{\vartheta_{\mathcal{D}}} \geq \vartheta_{\mathcal{D}}$ and may be strictly larger as individual indicator contributions cannot be subdivided. If $|\mathcal{M}_{K}| > K_{\max} - K$, the marked set is truncated to respect the remaining layer budget. The refined partition is then defined by
\begin{align}
\mathcal{T}_{K^{+}} = \operatorname{Refine} (\mathcal{T}_{K}, \mathcal{M}_{K}), \qquad K^{+} = K + |\mathcal{M}_{K}|,
\end{align}
where each marked interval $I_{k}$, $k \in \mathcal{M}_{K}$, is bisected at its midpoint $t_{k-\frac{1}{2}}$.

It remains to specify the layer-insertion procedure. Let $\tau$ and $\tau^+$ denote the discretizations associated with $\mathcal{T}_{K}$ and $\mathcal{T}_{K^+}$, respectively. The optimized control $\theta_{\tau}^\star$ on $\mathcal{T}_{K}$ is transferred to the refined grid $\mathcal{T}_{K^+}$ by the midpoint prolongation $\theta_{\tau^+}=\mathcal{P}_{\tau\to\tau^+}\theta_{\tau}^\star$, which preserves all existing nodal values and assigns
\begin{align*}
\bigl(
\mathcal{P}_{\tau\to\tau^+}\theta_\tau^\star
\bigr)
(t_{k-\frac12})
=
\frac{1}{2}
\left(
\theta_{\tau}^\star(t_{k-1})
+
\theta_{\tau}^\star(t_k)
\right),
\qquad
k\in\mathcal{M}_{K}.
\end{align*}
We then set $\theta_{\tau^+}^{(0)}:=\theta_{\tau^+}$ as the initial guess for the optimization on $\mathcal{T}_{K^+}$. Using the midpoint prolongation, one can also transfer the Adam moments \eqref{eq: Adam moments} from the coarser network to a finer optimization stage. Algorithm \ref{alg:stationarity-driven-layer-insertion} summarizes our adaptive approach.
\begin{algorithm}[Goal-oriented adaptive training]
\label{alg:stationarity-driven-layer-insertion}
\small
\begin{algorithmic}[1]
\Require Initial partition $\mathcal{T}_{K_0}$ and associated
discretization $\tau_0$; initial coefficient matrix
$\Theta_{\tau_0}^{(0)}$; depth budget $K_{\max}$; stationarity
tolerance $\varepsilon_{\mathrm{stat}}$; marking rule
$\operatorname{Mark}$; Adam state transfer flag
\Ensure Final adaptive partition, optimized coefficient matrix,
and phase history
%\State Freeze the calibrated input and output maps
\State Set $K\gets K_0$, $\mathcal{T}_K\gets\mathcal{T}_{K_0}$,
$\tau\gets\tau_0$, and
$\Theta_\tau^{(0)}\gets\Theta_{\tau_0}^{(0)}$

\While{the outer stopping criterion is not satisfied}
    \State Apply Algorithm~\ref{alg:fixed-depth-optimality-system}
    on $\mathcal{T}_K$ to obtain an approximation of
    $(\boldsymbol{x}_\tau^{\star},\Theta_\tau^\star,
    \boldsymbol{p}_\tau^{\star})$

    \If{$\mathfrak{s}_\tau\leq\varepsilon_{\mathrm{stat}}$}
        \State Record the phase as stationary

    \ElsIf{continuation without stationarity is enabled}
        \State Record the phase as nonstationary and permit
        exploratory refinement
        \State Exclude it from estimator and rate diagnostics
        requiring stationarity

    \Else
        \State \textbf{break}
    \EndIf

    \State Compute the local indicators
    $\eta_k$, $k=1,\ldots,K$, from the accepted fixed-depth solution

    \If{$K=K_{\max}$ or another outer stopping criterion is satisfied}
        \State \textbf{break}
    \EndIf

    \State
    $\mathcal{M}_K
    \gets
    \operatorname{Mark}
    \bigl(
    \mathcal{T}_K,
    \{\eta_k\}_{k=1}^{K}
    \bigr)$

    \State Enforce
    $|\mathcal{M}_K|\leq K_{\max}-K$

    \If{$\mathcal{M}_K=\varnothing$}
        \State \textbf{break}
    \EndIf

    \State
    $K^+\gets K+|\mathcal{M}_K|$

    \State
    $\mathcal{T}_{K^+}
    \gets
    \operatorname{Refine}
    \bigl(
    \mathcal{T}_K,
    \mathcal{M}_K
    \bigr)$

    \State
    $\Theta_{\tau^+}^{(0)}
    \gets
    \mathcal{P}_{\tau\to\tau^+}
    \Theta_\tau^\star$

    \If{Adam state transfer is enabled and the selected iterate
    was produced by Adam}
        \State Prolong the first and second Adam moments by the
        same nodal interpolation
    \Else
        \State Restart the Adam moments on the refined
        discretization $\tau^+$
    \EndIf

    \State Set
    $K\gets K^+$,
    $\mathcal{T}_K\gets\mathcal{T}_{K^+}$,
    $\tau\gets\tau^+$, and
    $\Theta_\tau^{(0)}\gets\Theta_{\tau^+}^{(0)}$
\EndWhile

\State \Return the final partition, optimized coefficient matrix,
and phase history
\end{algorithmic}
\end{algorithm}

\section{Numerical examples}
\label{sec 8: Numerical examples}

In this section, we present two numerical examples, one for binary classification and the other for multiclass classification; see also \cite{haber2017stable}. The examples address different questions. The Swiss
roll provides a controlled setting in which to examine the mechanics of the local indicators, their response to layer insertion and re-optimization, and the influence of the regularization parameter. The more challenging five-class
Peaks problem is then used for a systematic ablation of the marking rule, comparing indicator-driven, random, and D\"orfler layer insertion under a common depth budget.

\begin{figure}[h!]
    \centering
    \includegraphics[width=1\linewidth]{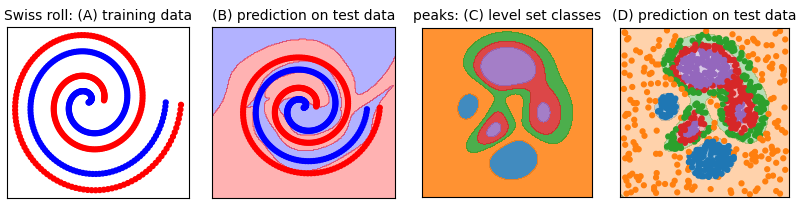}
    \caption{(\textbf{A}): Swiss roll training data with two colors indicating two different classes. (\textbf{B}): neural network prediction on the Swiss roll test data. (\textbf{C}): Peaks function level set classes. (\textbf{D}): neural network prediction on the Peaks test data.}
    \label{fig:swiss roll}
\end{figure}

\subsection{Binary classification}
We generate a two-class dataset, known as the Swiss roll, consisting of two concentric spirals in $\mathbb{R}^2$: 
\begin{align*}
f_{\text{blue}}(r, \varphi) = r (\cos(\varphi), \sin(\varphi))^{\top}, \quad f_{\text{red}}(r, \varphi) = (r +0.2) (\cos(\varphi), \sin(\varphi))^{\top},
\end{align*}
for $r \in [0,1]$ and $\varphi \in [0, 4\pi]$. We generate $1200$
samples from each class and associate with every input $x_{0}^{i}\in
\mathbb{R}^{2}$ a label $y^{i}\in\{0,1\}$ identifying the spiral from
which it was sampled, see Fig.~\ref{fig:swiss roll}(A). The complete data set is divided by a stratified 60/40 training--validation split. For this binary classification task, the neural ODE yields a score \eqref{output mapping}, where $\hat{y}^i \in (0,1) $ and $q_{\mathrm{out}}(s) = (1 + e^{-s})^{-1}$. The empirical binary cross-entropy 
\begin{align}\label{validation loss}
l(\boldsymbol{x}(T)):= -\frac{1}{m} \sum_{i=1}^m \big[ 
y^i \log \hat{y}^i + (1 - y^i)\log(1 - \hat{y}^i)\big]
\end{align}
is used to define the objective \eqref{Cross-entropy function}. Evaluating \eqref{validation loss} on the training and validation data yields the training and validation loss values. The neural ODE has hidden dimension $d=4$, terminal time $T=2.5$, and $\tanh$ activation. As an outer stopping criteria in Algorithm \ref{alg:stationarity-driven-layer-insertion}, we use the maximal layer budget: the refinement starts from $K=2$ intervals and continues to $K_{\max}=12$, with fixed-depth stationarity tolerance $\varepsilon_{\mathrm{stat}}=10^{-3}$ and switching parameter $\kappa_{\mathrm{s}}=5.0$ in Algorithm \ref{alg:fixed-depth-optimality-system}. We consider $\lambda \in \{10^{-2},10^{-3},10^{-4},10^{-5},10^{-6}\}$ to numerically study the influence of regularization parameter. 

\begin{figure}[h!]
    \centering
    \includegraphics[width=0.95\linewidth]{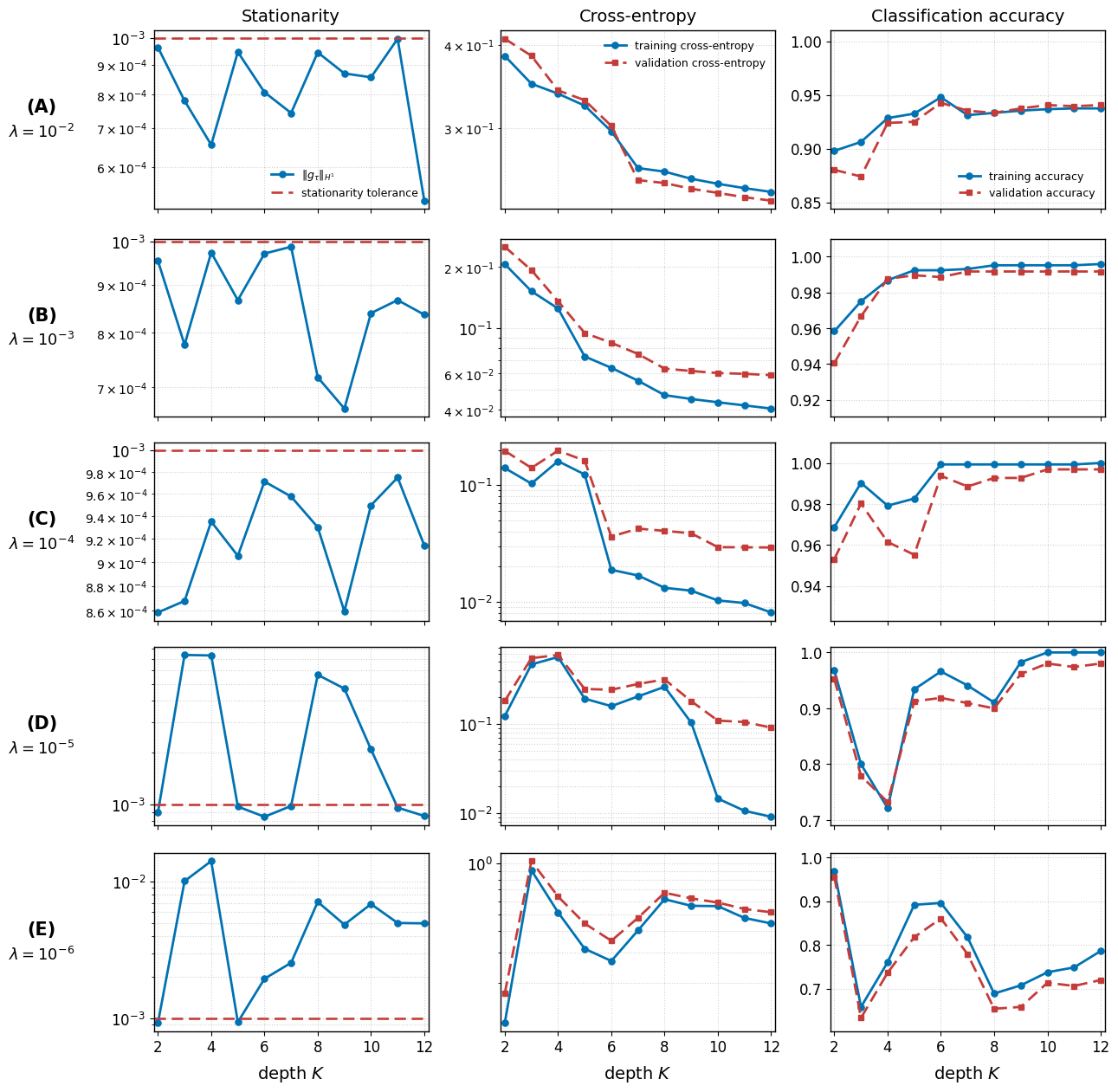}
    \caption{Swiss roll adaptive training: rows correspond to different regularization parameters; columns show the $H^1$ stationarity norm, training/validation cross-entropy loss, and training/validation accuracy.}
    \label{fig:SR learning}
\end{figure}

Table~\ref{tab:combined-swiss-new} summarizes the Swiss roll ablation, with selected learning curves shown in Figure~\ref{fig:SR learning}. We observe that the regularization parameter controls the stability of the adaptive hierarchy. For $\lambda=10^{-2}$ and $\lambda=10^{-3}$, every fixed-depth phase reaches stationarity, and the loss and accuracy improve monotonically along the hierarchy. As $\lambda$ decreases, stationarity becomes progressively harder to reach and the learning curves turn non-monotone, with the stationarity norm stagnating above the tolerance. This is reflected in the final quality: for $\lambda=10^{-6}$ the conservative run deteriorates to an accuracy worse than the strongly regularized $\lambda=10^{-2}$ case, now caused by a failure to optimize rather than by limited expressivity coming from the regularization. The weakly regularized runs also require substantially more optimization work measured by wall-clock time. Indeed, in such cases, the prolongated coarse solution may lie farther from a stationary point of the refined problem, resulting in greater optimization effort, as measured by both the number of Adam steps and the BFGS iteration count reported in Table~\ref{tab:combined-swiss-new}.

\begin{figure}[h!]
    \centering
    \includegraphics[width=0.85\linewidth]{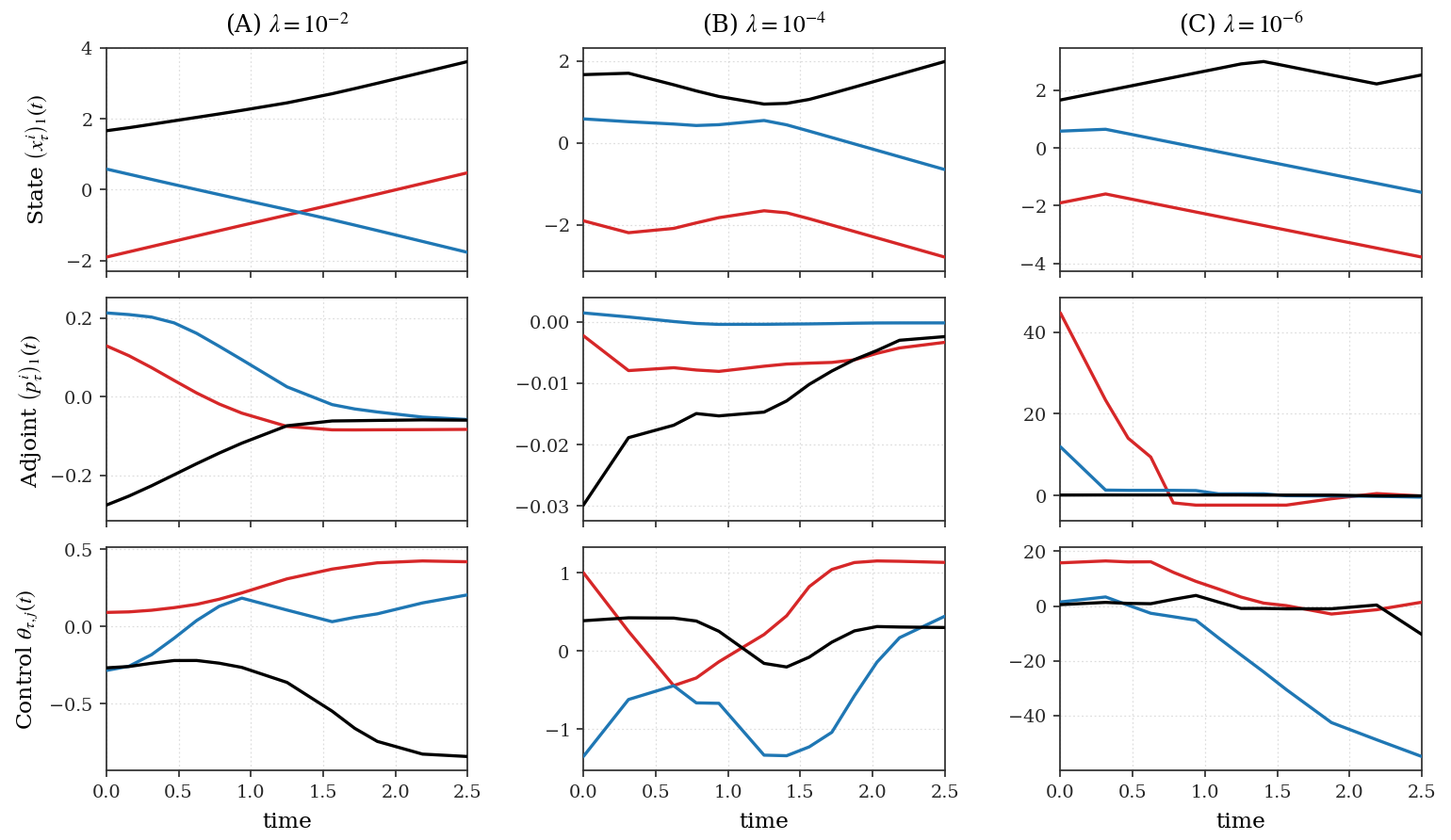}
    \caption{Representative trajectories \(t\mapsto\bx_{\tau}^{i}(t)\) for the state, \(t\mapsto\bp_{\tau}^{i}(t)\) for the adjoint, and \(t\mapsto\theta_{\tau}^{i}(t)\) for the control over the adaptive hierarchy for \textbf{(A)} \(\lambda=10^{-2}\), \textbf{(B)} \(\lambda=10^{-3}\), and \textbf{(C)} \(\lambda=10^{-4}\). Stronger regularization yields smoother state--adjoint--control trajectories with less variability.}
    \label{fig:trajectories}
\end{figure}

\begin{table}[!ht]
\centering
\scriptsize
\setlength{\tabcolsep}{3.0pt}
\renewcommand{\arraystretch}{1.08}
\caption{Swiss roll: final quality, reliability, and adaptive computational
work for the conservative and D\"orfler marking strategies. The stationarity column reports the percentage of fixed-depth phases satisfying
$\mathfrak{s}_{\tau}=\|g_\tau\|_{H^1}
\leq\varepsilon_{\mathrm{stat}}$. The Adam and BFGS columns report the total
iteration counts over the complete adaptive hierarchy.}
\label{tab:combined-swiss-new}
\resizebox{\textwidth}{!}{%
\begin{tabular}{lllcccccrrr}
\toprule
\multicolumn{3}{c}{Parameters}
&
\multicolumn{3}{c}{Final quality}
&
\multicolumn{2}{c}{Reliability}
&
\multicolumn{3}{c}{Adaptive work}
\\
\cmidrule(lr){1-3}
\cmidrule(lr){4-6}
\cmidrule(lr){7-8}
\cmidrule(lr){9-11}
$\lambda$
& Family
& Rule / $\vartheta_{\mathcal D}$
& final $\J_\lambda$
& val.\ loss
& val.\ accuracy [\%]
& stationarity [\%]
& \shortstack{all\\$r_{\mathrm{child}}<1$}
& \shortstack{Adam\\steps}
& \shortstack{BFGS\\iters.}
& time [s]
\\
\midrule

\multirow{4}{*}{$10^{-2}$}
& Conservative & Maximum
& $0.434270$ & $0.229760$ & $95.52$
& $100$ & no & $16\,000$ & $422$ & $47.8$
\\
& D\"orfler & $0.25$
& $0.435028$ & $0.232932$ & $95.42$
& $100$ & no & $15\,300$ & $400$ & $36.4$
\\
& D\"orfler & $0.50$
& $0.435027$ & $0.233103$ & $95.42$
& $100$ & no & $13\,100$ & $258$ & $27.0$
\\
& D\"orfler & $0.75$
& $0.490028$ & $0.322896$ & $90.31$
& $100$ & no & $11\,400$ & $333$ & $24.0$
\\

\addlinespace[3pt]

\multirow{4}{*}{$10^{-3}$}
& Conservative & Maximum
& $0.127845$ & $0.060734$ & $99.38$
& $100$ & no & $23\,500$ & $2\,728$ & $59.8$
\\
& D\"orfler & $0.25$
& $0.128824$ & $0.061539$ & $99.27$
& $100$ & no & $22\,600$ & $2\,624$ & $53.7$
\\
& D\"orfler & $0.50$
& $0.124228$ & $0.058747$ & $99.38$
& $100$ & no & $22\,200$ & $1\,109$ & $38.1$
\\
& D\"orfler & $0.75$
& $0.113583$ & $0.061611$ & $99.27$
& $100$ & no & $16\,300$ & $863$ & $27.3$
\\

\addlinespace[3pt]

\multirow{4}{*}{$10^{-4}$}
& Conservative & Maximum
& $0.195254$ & $0.210864$ & $85.94$
& $81.8$ & no & $50\,600$ & $9\,016$ & $107.0$
\\
& D\"orfler & $0.25$
& $0.212131$ & $0.232548$ & $88.75$
& $88.9$ & no & $40\,400$ & $9\,410$ & $94.9$
\\
&  D\"orfler &  $0.50$
&  $0.009880$ &  $0.001618$ &  $100.00$
&  $100$ &  yes &  $22\,300$ &  $7\,214$ &  $69.4$
\\
& D\"orfler & $0.75$
& $0.014234$ & $0.002888$ & $100.00$
& $66.7$ & no & $38\,000$ & $10\,464$ & $95.9$
\\

\addlinespace[3pt]

\multirow{4}{*}{$10^{-5}$}
& Conservative & Maximum
& $0.016436$ & $0.067812$ & $98.44$
& $90.9$ & no & $42\,400$ & $29\,440$ & $196.3$
\\
& D\"orfler & $0.25$
& $0.006523$ & $0.001114$ & $100.00$
& $80.0$ & no & $44\,000$ & $31\,130$ & $201.2$
\\
& D\"orfler & $0.50$
& $0.080139$ & $0.123051$ & $97.71$
& $71.4$ & no & $30\,500$ & $17\,390$ & $125.7$
\\
& D\"orfler & $0.75$
& $0.009206$ & $0.030139$ & $99.69$
& $80.0$ & no & $19\,400$ & $24\,590$ & $144.3$
\\

\addlinespace[3pt]

\multirow{4}{*}{$10^{-6}$}
& Conservative & Maximum
& $0.286500$ & $0.334552$ & $86.15$
& $54.5$ & no & $60\,000$ & $26\,757$ & $211.6$
\\
& D\"orfler & $0.25$
& $0.471307$ & $0.475620$ & $76.04$
& $44.4$ & no & $55\,100$ & $22\,835$ & $175.9$
\\
& D\"orfler & $0.50$
& $0.315465$ & $0.397232$ & $84.58$
& $28.6$ & no & $40\,300$ & $15\,415$ & $131.0$
\\
& D\"orfler & $0.75$
& $0.357055$ & $0.405849$ & $86.25$
& $16.7$ & no & $30\,500$ & $20\,367$ & $147.1$
\\

\bottomrule
\end{tabular}}
\end{table}
\noindent
Stronger regularization stabilizes the hierarchy at the price of expressivity, as seen directly in the state-adjoint-control triplet dynamics: for $\lambda=10^{-2}$ the state trajectories in Figure~\ref{fig:trajectories} are smooth and slowly varying, whereas decreasing $\lambda$ yields markedly more oscillatory trajectories, i.e., a richer feature transformation. The strongly regularized runs therefore cannot reduce the loss below a moderate floor, while the best approximation quality is attained in the intermediate regime, where the hierarchy remains stable but sufficiently expressive. Overall, $\lambda=10^{-3}$ offers the most robust compromise, combining full stationarity with over $99\%$ accuracy across all marking strategies, whereas $\lambda\leq10^{-5}$ trades a chance of higher accuracy against a loss of reliability.

\begin{figure}[h!]
\centering
\includegraphics[width=0.7\linewidth]{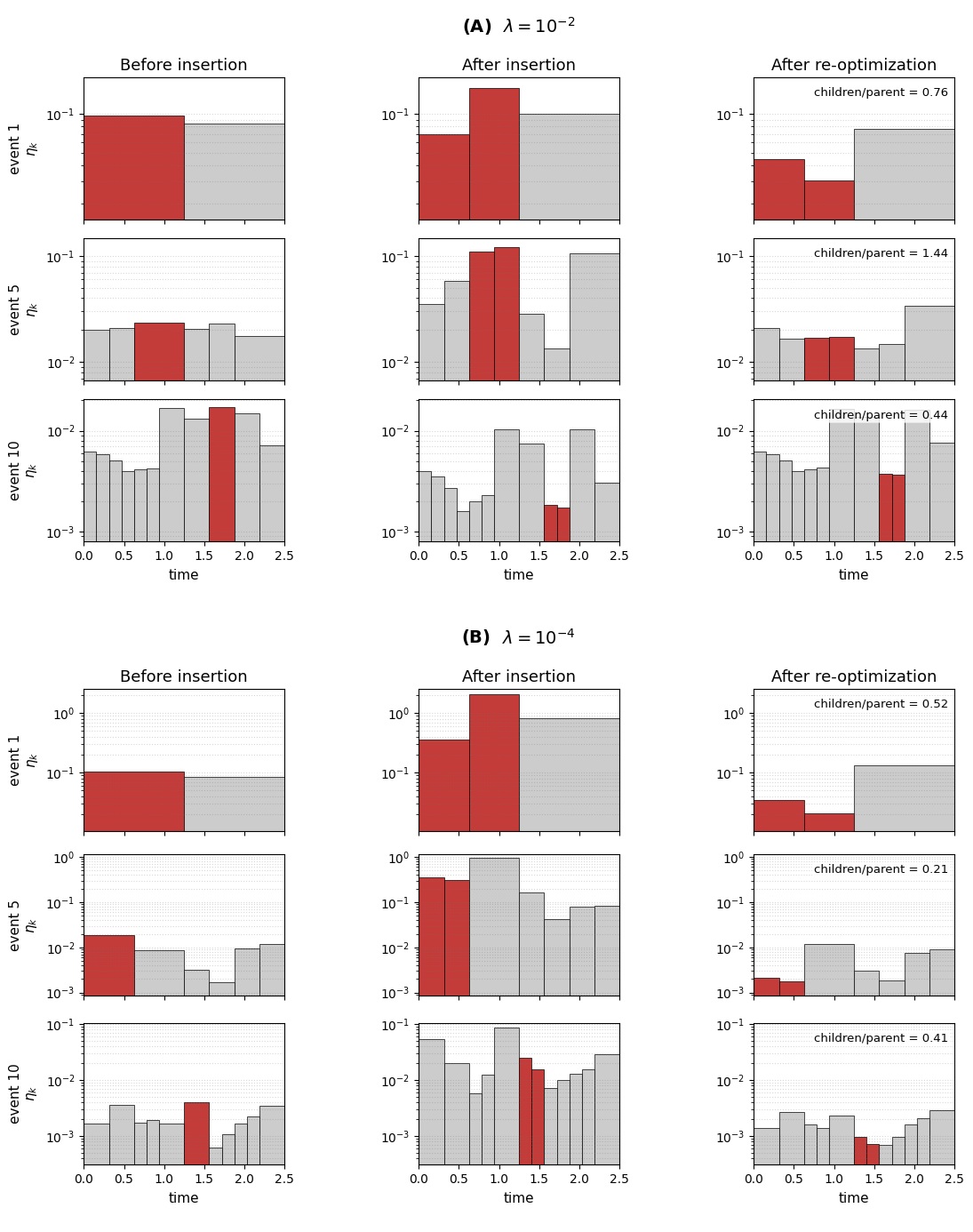}
\caption{Selected conservative refinement events
for \textbf{(A)}: $\lambda=10^{-2}$ and  \textbf{(B)}: $\lambda = 10^{-4}$. For each interval, local indicators are shown before insertion, immediately afterward, and after re-optimization. Red bars denote the marked parent before insertion and its two children afterward. The reported children/parent ratio is defined in \eqref{eq:children parent ratio}.}
\label{fig:conservative layer insertion}
\end{figure}

To quantify the effect of re-optimization on the indicator and on the overall adaptive process, for a parent interval $I_{k^\star}$ selected at depth $K$, we define its children-to-parent ratio by
\begin{align}\label{eq:children parent ratio}
r_{\mathrm{child},K}
=
\frac{
\eta_{K+1,\mathrm{child},1}^{\mathrm{post}}
+
\eta_{K+1,\mathrm{child},2}^{\mathrm{post}}
}{
\eta_{K,k^\star}^{\mathrm{before}}
}.
\end{align}
Here, ``post'' refers to the two child contributions after global
re-optimization on the refined grid. The condition
$r_{\mathrm{child},K}<1$ means that refinement and re-optimization reduce the combined contribution of the marked parent. However, we want to point out that this is only a local consistency diagnostic and not a guarantee of a monotonic decrease in the global indicator along the adaptive hierarchy. Indeed, the local indicators associated with unmarked intervals may also change during re-optimization, although in practice they tend to decrease rather than increase. Observe in Figure~\ref{fig:conservative layer insertion} that immediate layer insertion may increase the local indicators, whereas global re-optimization reduces the combined child contribution in the displayed events. The local indicators change only mildly for $\lambda=10^{-2}$ between the ``After insertion'' and ``After re-optimization'' stages. At the same time, weaker regularization $\lambda=10^{-4}$ produces a larger increase in the indicators after layer insertion and requires more re-optimization work because of the greater variability in the control--state--adjoint triplet between consecutive designs. This further supports the observation that larger values of $\lambda$ stabilize the adaptive sequence while limiting variability (and hence expressivity) of the triplet, whereas smaller values may yield more expressive networks but lead to less stable adaptive dynamics. We remark that we do not observe a strong general pattern in the unmarked intervals adjacent to the refined ones. Their indicator values do not tend to increase significantly and instead tend to decrease after re-optimization.

\begin{figure}[h!]
    \centering
    \includegraphics[width=0.95\linewidth]{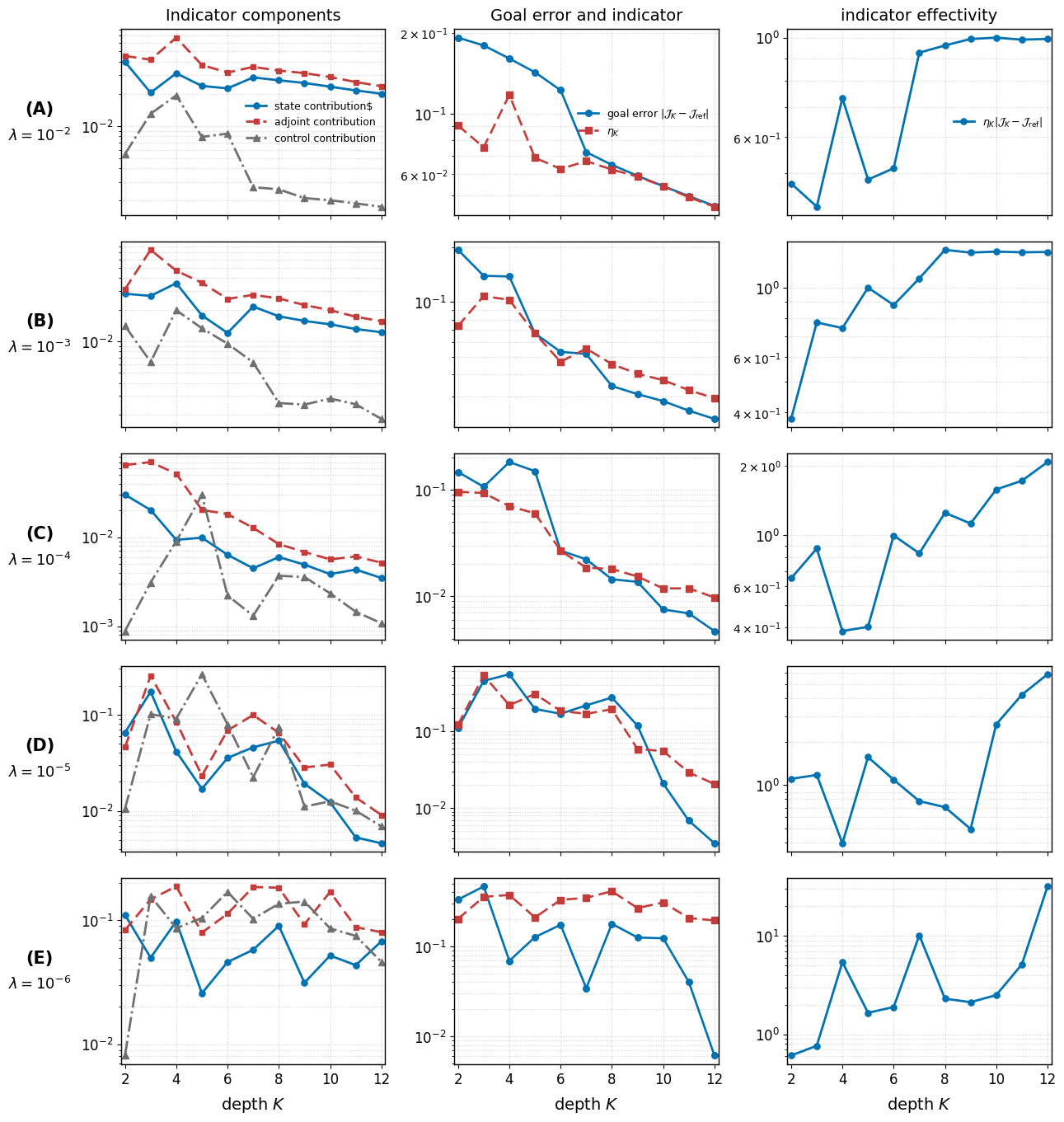}
    \caption{Indicator components, goal error vs global indicator $\eta_{K}$ for depth $K$, overall, and indicator effectivities.}
    \label{fig:indicator picture}
\end{figure}

To assess the global indicator, we build a nested reference hierarchy $\mathcal T^{\,(0)}_{12}\to\mathcal T^{\,(1)}_{24}\to\mathcal T^{\,(2)}_{48}$ by uniformly bisecting the final adaptive grid and prolonging and re-optimizing the control at each level. We take $\mathfrak J_{\tau^{(2)}}(\Theta_{\tau^{(2)}}^\star)$ as the reference value. Figure~\ref{fig:indicator picture} reports the three indicator components, the global indicator $\eta_K$ against the goal error $|\mathfrak J(\Theta)-\mathfrak J_{\mathrm{ref}}(\Theta_{\tau^{(2)}}^\star)|$, and the effectivity. For $\lambda\in\{10^{-2},10^{-3},10^{-4}\}$, the indicator tracks the goal error closely, with both decreasing with $K$ and with effectivity of order one. Observe that, at the early stages, the indicator may slightly underestimate the goal error. This is possibly related to neglecting the higher-order Taylor remainder, which may be non-negligible at these stages. However, for larger $K$, the indicator is reasonably calibrated.  For $\lambda\leq10^{-5}$, the effectivity drifts from one and grows with $K$, consistently with the loss of stationarity in Figure~\ref{fig:SR learning} (D), which may result in a larger optimization error dominating the discretization error in the objective measured by the indicator. Observe also that the three components of the indicator in Figure~\ref{fig:indicator picture} (A) are quite balanced in most cases, thus playing an equally important role in the adaptive process. Since $r_{\mathrm{child}}<1$ fails at least once in most cases according to Table~\ref{tab:combined-swiss-new}, re-optimization need not reduce the marked contribution in some events, and monotonicity of the global indicator cannot be expected. The $\lambda=10^{-6}$ case admits no meaningful sharpness assessment, making sufficient regularization and fixed-depth stationarity prerequisites for a clean adaptive interpretation.

\begin{figure}[h!]
    \centering
    \includegraphics[width=0.7\linewidth]{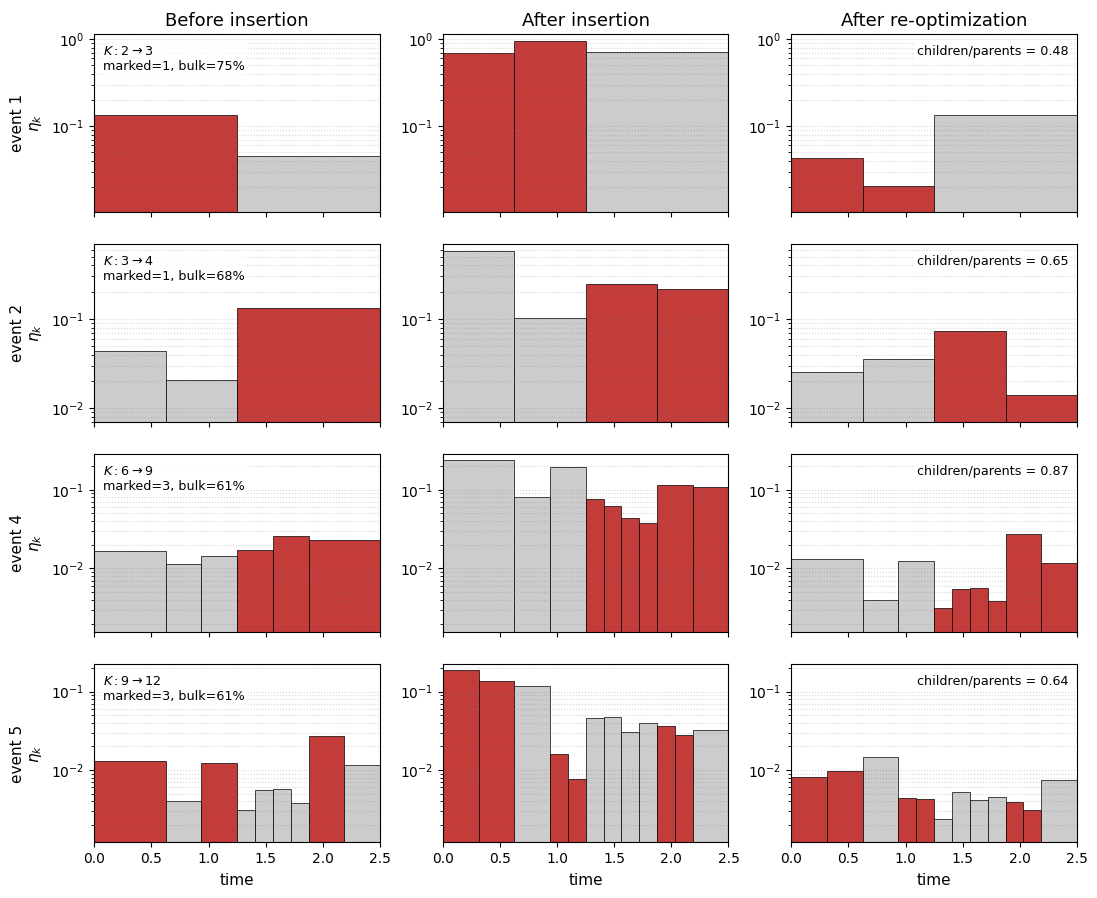}
    \caption{Selected D\"orfler refinement events
for $\lambda=10^{-3}$ and $\vartheta_{\mathcal{D}}=0.5$. For each interval, local indicators are shown before insertion, immediately afterward, and after re-optimization. Red bars denote the marked parents before insertion and their children afterward. The achieved bulk fraction for each event is computed according to \eqref{empirical bulk}.}
\label{fig:Dorfler refinement events}
\end{figure}

Results for D\"orfler marking are also reported in Table~\ref{tab:combined-swiss-new}. In the stable regimes $\lambda\in\{10^{-2},10^{-3}\}$, moderate marking with $\vartheta_{\mathcal D}\in{0.25,0.50}$ reproduces the validation performance of conservative refinement while reducing runtime. In particular, $\vartheta_{\mathcal D}=0.50$ reduces the conservative wall time by $44\%$ at $\lambda=10^{-2}$ and by $36\%$ at $\lambda=10^{-3}$. In the intermediate regime, the marking fraction becomes decisive. At $\lambda=10^{-4}$, only $\vartheta_{\mathcal D}=0.50$ reaches full stationarity and attains the best overall result, with a wall time $35\%$ shorter than that of the conservative run. By contrast, the conservative and $\vartheta_{\mathcal D}=0.25$ runs stall, while the $\vartheta_{\mathcal D}=0.75$ run fails to reach stationarity. For $\lambda\leq10^{-5}$, the outcome is highly sensitive to $\vartheta_{\mathcal D}$. Figure~\ref{fig:Dorfler refinement events} illustrates the mechanics of simultaneous layer placement. Overall, moderate D\"orfler marking reduces wall time by roughly one-third to almost one-half relative to conservative refinement in the stable regimes and can additionally stabilize the intermediate regime. By contrast, aggressive marking and weak regularization make the outcome increasingly dependent on the optimization trajectory.

\begin{figure}[h!]
    \centering
    \includegraphics[width=1\linewidth]{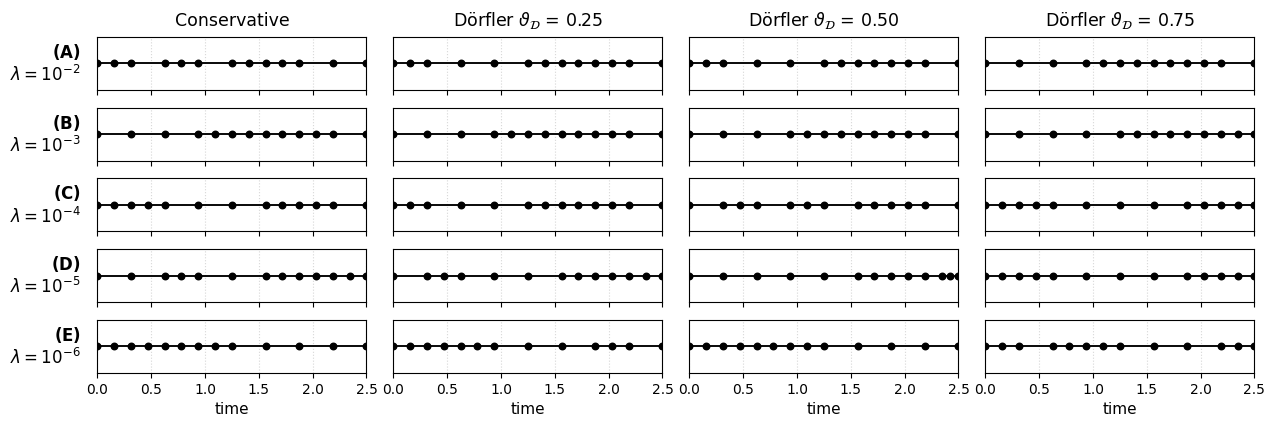}
    \caption{Grids produced by conservative marking (maximum-indicator insertion), random layer insertion, and minimum-indicator insertion across several random seeds.}
    \label{fig:SR grids}
\end{figure}

Finally, Figure~\ref{fig:SR grids} shows the grids produced by conservative maximum-indicator insertion and D\"orfler marking. We observe that the conservative strategy and D\"orfler marking with moderate marking fractions produce very similar grids under moderate regularization, whereas the differences between the grids become more visible for smaller regularization values or more aggressive marking. In view of the accuracy results displayed in Table \ref{tab:combined-swiss-new}, this suggests that D\"orfler marking can provide a plausible and less expensive alternative to conservative marking for stable adaptive hierarchies.

\subsection{Multiclass classification}

We next isolate the role of the marking rule on the five-class Peaks problem of
\cite{haber2017stable}, using the direct reconstruction-based indicator
\eqref{error indicator approximation} throughout. The data are generated from the
Peaks function
\begin{align*}
f_{\mathrm{ps}}(x)
={}&
3(1-x_{1})^2
\exp\!\left(-x_{1}^{2}-(x_{2}+1)^2\right)
-10\left(\frac{x_{1}}{5}-x_{1}^{3}-x_{2}^{5}\right)
\exp\!\left(-x_{1}^{2}-x_{2}^{2}\right) \\
&-\tfrac{1}{3}
\exp\!\left(-(x_{1}+1)^2-x_{2}^{2}\right),
\qquad x\in[-3,3]^2 ,
\end{align*}
discretized on a regular $256\times256$ grid, with classes
$\mathcal C_i=\{x:\ c_{i-1}\leq f_{\mathrm{ps}}(x)<c_i\}$ for $i=1,\dots,4$ and
$\mathcal C_5=\{x:\ f_{\mathrm{ps}}(x)\geq c_4\}$, where $c_0=\min f_{\mathrm{ps}}$
and $(c_1,c_2,c_3,c_4)=(-2.2,0.55,1.75,3.2)$, see also Figure \ref{fig:swiss roll} (C). We sample $1500$ points per class
and use a fixed stratified $60/40$ split, giving $4500$ training and $3000$
validation samples. With one-hot labels $y^i\in\mathbb R^{d_{\mathrm{out}}}$,
$d_{\mathrm{out}}=5$, we use the softmax output map and multiclass cross-entropy
\begin{align}\label{multiclass validation loss}
l(\boldsymbol{x}(T)) = -\frac{1}{m} \sum_{i=1}^{m}\sum_{j=1}^{d_{\mathrm{out}}} y_j^i\log\widehat y_j^i, \qquad q_{\mathrm{out}}(s^j) = \frac{\exp(s^j)}{\sum_{k=1}^{d_{\mathrm{out}}}\exp(s^k)} .
\end{align}
The neural ODE has width $d=5$, $\tanh$ activation, and terminal time $T=5$. Every adaptive hierarchy starts from $K_0=2$ and is refined to $K_{\max}=15$. We consider $\lambda \in \{10^{-2},10^{-3}\}$ and ten model initialization seeds $0,5,\ldots,45$, each fixing an initial guess shared by all marking rules. For each seed and $\lambda$, the conservative ablation builds one maximum-indicator, one minimum-indicator, and ten random-insertion hierarchies (averaged within each seed before aggregation), for $10\times2\times(1+1+10)=240$ runs; the D\"orfler ablation adds one hierarchy per $\vartheta_{\mathcal D}\in\{0.25,0.50,0.75\}$, for a further $10\times2\times3=60$ runs. Minimum-indicator insertion serves as a deliberately opposing strategy, refining the intervals identified as least important by the proposed indicator. In addition, to separate the effect of adaptive continuation from direct optimization at the final depth, we train a nonadaptive uniform-grid baseline with $K_{\max}=15$ from initialization. In Algorithm \ref{alg:fixed-depth-optimality-system}, we set $\varepsilon_{\mathrm{stat}}=5\times10^{-3}$ and $\kappa_{\mathrm{s}}=1.5$.

\begin{figure}[h!]
    \centering
    \includegraphics[width=1\linewidth]{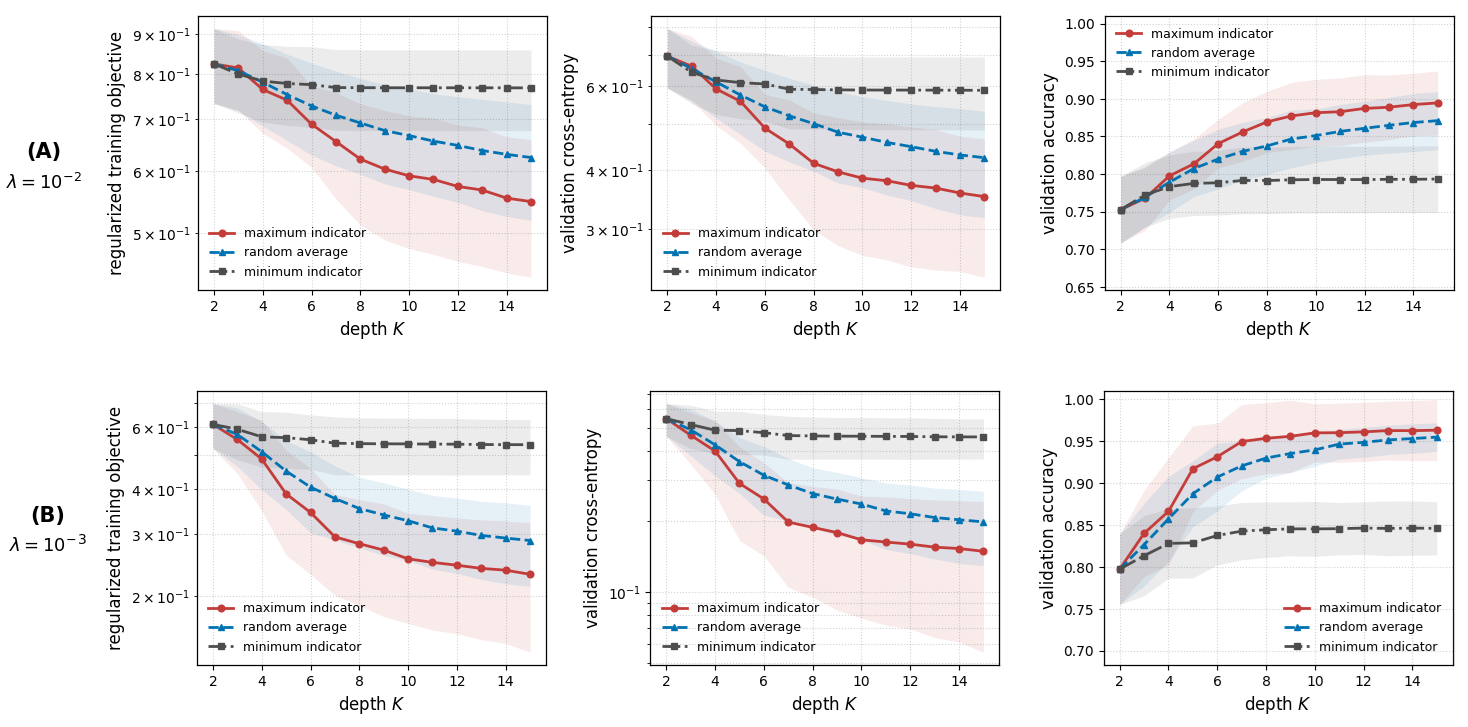}
    \caption{Peaks conservative marking ablation for \textbf{(A)}: $\lambda=10^{-2}$ and
    \textbf{(B)}: $\lambda=10^{-3}$.}
    \label{fig:peaks ablation}
\end{figure}

Table~\ref{tab:peaks table} and Figure~\ref{fig:peaks ablation} show that the ordering of the local indicators is informative for layer placement. Relative to random insertion, maximum-indicator refinement reduces the validation loss by $21.2\%$ for $\lambda=10^{-2}$ and $24.8\%$ for $\lambda=10^{-3}$; relative to minimum-indicator insertion, the reductions are $42.0\%$ and $61.7\%$. The corresponding accuracy gains are $3.03$ and $1.53$ percentage points over random insertion and $9.77$ and $10.33$ percentage points over minimum-indicator insertion. Compared with the uniform full-depth baseline, maximum-indicator refinement further reduces the validation loss by $12.9\%$ and $9.8\%$, with accuracy gains of $1.89$ and $1.22$ percentage points, respectively. The uniform grid solve reaches the prescribed stationarity tolerance for only $6/10$ and $1/10$ seeds, whereas maximum-indicator refinement satisfies it in $95.0\%$ and $87.1\%$ of its fixed-depth phases. Thus, although direct full-depth training requires less optimization work, adaptive continuation is substantially more reliable and yields better mean final quality. Taken together, these comparisons indicate that the gain is not explained by depth alone, but by the combination of indicator-informed layer placement and coarse-to-fine adaptive optimization.

\begin{table}[!ht]
\centering
\scriptsize
\setlength{\tabcolsep}{3.0pt}
\renewcommand{\arraystretch}{1.08}
\caption{Peaks: final quality, reliability, and optimization work for the
adaptive marking ablations and the nonadaptive uniform-grid baseline. Values
are reported as mean $\pm$ standard deviation over ten model initializations.
For random marking, the results are first averaged over ten random insertion paths for each initialization. The uniform baseline involves one fixed-depth optimization per initialization: its stationarity percentage reports the fraction of model seeds for which the final solution satisfies $\mathfrak{s}_{\tau}\leq\varepsilon_{\mathrm{stat}}$.}
\label{tab:peaks table}
\resizebox{\textwidth}{!}{%
\begin{tabular}{lllccccrrr}
\toprule
\multicolumn{3}{c}{Parameters}
&
\multicolumn{3}{c}{Final quality}
&
\multicolumn{1}{c}{Reliability}
&
\multicolumn{3}{c}{Optimization work}
\\
\cmidrule(lr){1-3}
\cmidrule(lr){4-6}
\cmidrule(lr){7-7}
\cmidrule(lr){8-10}
$\lambda$
& Family
& Rule / grid / $\vartheta_{\mathcal D}$
& final $\J_\lambda$
& val.\ loss
& val.\ accuracy [\%]
& \shortstack{stationarity [\%]}
& \shortstack{Adam\\steps}
& \shortstack{BFGS\\iters.}
& time [s]
\\
\midrule

\multirow{7}{*}{$10^{-2}$}
& Conservative & Maximum
& $0.5359\pmc0.0951$
& $0.3346\pmc0.0787$
& $90.13\pmc2.98$
& $95.0\pmc11.2$
& $27\,420\pmc12\,622$
& $3\,059\pmc3\,124$
& $126.0\pmc28.2$
\\
& Conservative & Random
& $0.6252\pmc0.1058$
& $0.4246\pmc0.1073$
& $87.10\pmc3.85$
& $96.7\pmc4.4$
& $17\,619\pmc4\,936$
& $1\,379\pmc907$
& $110.2\pmc12.9$
\\
& Conservative & Minimum
& $0.7646\pmc0.0816$
& $0.5771\pmc0.0949$
& $80.36\pmc4.50$
& $99.3\pmc2.3$
& $11\,000\pmc2\,953$
& $197\pmc194$
& $111.9\pmc17.6$
\\
& D\"orfler & $0.25$
& $0.5382\pmc0.0941$
& $0.3353\pmc0.0774$
& $90.18\pmc3.09$
& $94.0\pmc9.7$
& $24\,230\pmc10\,253$
& $3\,131\pmc3\,254$
& $96.7\pmc22.2$
\\
&  D\"orfler
&  $0.50$
&  $0.5591\pmc0.0989$
&  $0.3595\pmc0.0903$
&  $89.28\pmc3.43$
&  $96.3\pmc8.4$
&  $21\,250\pmc8\,015$
&  $1\,900\pmc1\,092$
&  $70.8\pmc16.1$
\\
& D\"orfler & $0.75$
& $0.5842\pmc0.1197$
& $0.3873\pmc0.1406$
& $88.59\pmc5.33$
& $79.3\pmc24.8$
& $17\,300\pmc5\,486$
& $1\,011\pmc678$
& $49.5\pmc8.8$
\\
\addlinespace[1pt]
& Nonadaptive & Uniform $K_{\max}$
& $0.5696\pmc0.1530$
& $0.3840\pmc0.1725$
& $88.24\pmc7.34$
& $60.0\pmc51.6$
& $5\,310\pmc5\,307$
& $694\pmc1\,006$
& $13.9\pmc8.7$
\\

\addlinespace[3pt]

\multirow{7}{*}{$10^{-3}$}
& Conservative & Maximum
& $0.2301\pmc0.0820$
& $0.1497\pmc0.0786$
& $97.04\pmc2.03$
& $87.1\pmc15.7$
& $46\,250\pmc15\,386$
& $14\,494\pmc3\,728$
& $226.2\pmc34.9$
\\
& Conservative & Random
& $0.2867\pmc0.0737$
& $0.1991\pmc0.0696$
& $95.51\pmc1.70$
& $92.3\pmc6.8$
& $30\,906\pmc7\,064$
& $8\,596\pmc4\,025$
& $187.6\pmc35.5$
\\
& Conservative & Minimum
& $0.4836\pmc0.1118$
& $0.3908\pmc0.1201$
& $86.71\pmc5.11$
& $97.1\pmc9.0$
& $15\,560\pmc5\,011$
& $975\pmc826$
& $156.4\pmc30.2$
\\
& D\"orfler & $0.25$
& $0.2259\pmc0.0884$
& $0.1473\pmc0.0816$
& $96.67\pmc2.85$
& $83.0\pmc12.6$
& $43\,330\pmc10\,198$
& $13\,082\pmc5\,113$
& $183.8\pmc35.3$
\\
&  D\"orfler
&  $0.50$
&  $0.2337\pmc0.0610$
&  $0.1488\pmc0.0536$
& $97.11\pmc1.51$
&  $72.5\pmc29.7$
& $38\,510\pmc8\,551$
&  $7\,422\pmc4\,096$
&  $120.4\pmc24.1$
\\
& D\"orfler & $0.75$
& $0.2436\pmc0.0692$
& $0.1575\pmc0.0580$
& $96.72\pmc1.49$
& $66.0\pmc23.2$
& $31\,030\pmc7\,533$
& $7\,740\pmc2\,632$
& $101.2\pmc10.7$
\\
\addlinespace[1pt]
& Nonadaptive & Uniform $K_{\max}$
& $0.2387\pmc0.0751$
& $0.1659\pmc0.0740$
& $95.82\pmc1.87$
& $10.0\pmc31.6$
& $6\,410\pmc3\,114$
& $1\,888\pmc2\,459$
& $20.6\pmc9.4$
\\
\bottomrule
\end{tabular}}
\end{table}

\begin{figure}[!h]
    \centering
    \includegraphics[width=0.9\linewidth]{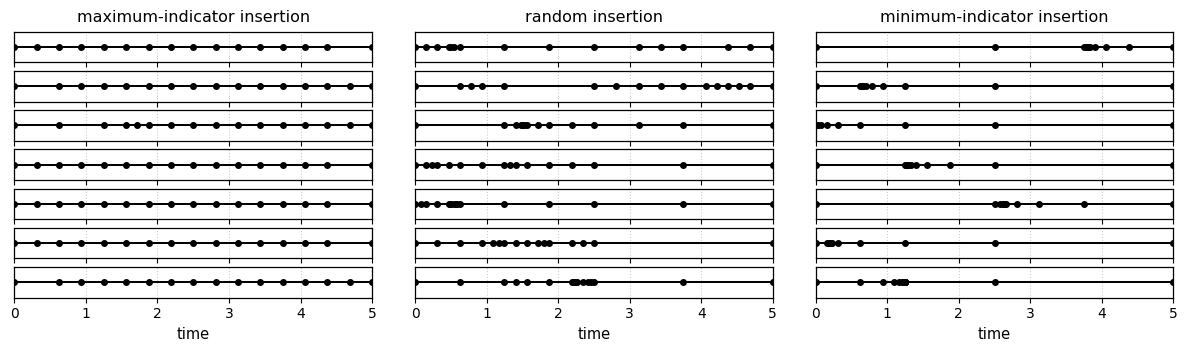}
    \caption{Final grids obtained for $\lambda=10^{-2}$ by conservative
    maximum-indicator insertion, random insertion, and minimum-indicator
    insertion across selected model initializations. Black circles denote the
    temporal-grid nodes.}
    \label{fig:peaks conservative grids 1e-2}
\end{figure}

\begin{figure}[!h]
    \centering
    \includegraphics[width=0.9\linewidth]{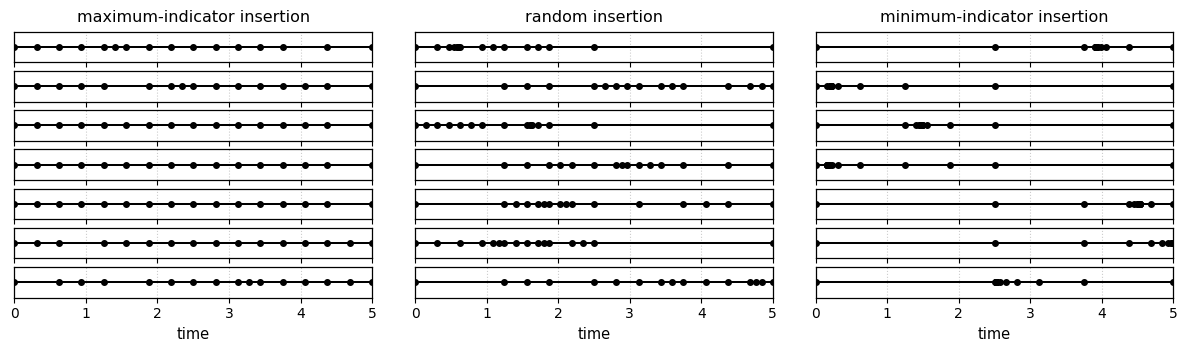}
    \caption{Final grids obtained for $\lambda=10^{-3}$ by conservative
    maximum-indicator insertion, random insertion, and minimum-indicator
    insertion across selected model initializations. Black circles denote the
    temporal-grid nodes.}
    \label{fig:peaks conservative grids 1e-3}
\end{figure}

The corresponding grids are shown in Figure~\ref{fig:peaks conservative grids 1e-2} for $\lambda=10^{-2}$, and for $\lambda=10^{-3}$ in Figure~\ref{fig:peaks conservative grids 1e-3}. Maximum-indicator insertion produces a consistent distribution of nodes over the horizon, whereas random insertion varies strongly across realizations and minimum-indicator insertion repeatedly leaves large temporal regions coarse. Observe that, for strong regularization with $\lambda=10^{-2}$, the resulting adaptive grids appear rather uniform. This is expected because the temporal variability of the control and, consequently, of the state and adjoint variables is limited in this case. In general, this supports the interpretation that the indicators identify temporally important parts of the dynamics rather than merely favoring additional depth.

Table~\ref{tab:peaks table} shows that D\"orfler marking reduces the cost of repeated global re-optimization while retaining the benefits of the conservative strategy. For both $\lambda$, $\vartheta_{\mathcal D}=0.25$ and $0.50$ match the conservative validation quality while cutting the mean adaptive time by up to $43.8\%$ ($\lambda=10^{-2}$) and $46.8\%$ ($\lambda=10^{-3}$) at $\vartheta_{\mathcal D}=0.50$. The aggressive choice $\vartheta_{\mathcal D}=0.75$ saves a further margin ($60.7\%$ and $55.3\%$ relative to conservative) but is less reliable in reaching stationarity, with a visible deterioration in validation quality. Across both regimes, $\vartheta_{\mathcal D}=0.50$ offers the most balanced compromise between quality, reliability, and cost, supporting the use of the  local indicators to guide layer refinement.

%\begin{figure}[!h]
%    \centering
%    \includegraphics[width=1\linewidth]{Figures/Dofler grids 1e-2.png}
%    \caption{Grids produced by Dörfler marking with $\lambda = 10^{-2}$ and with fraction parameters $\vartheta_{\mathcal{D}}=0.25$, $\vartheta_{\mathcal{D}}=0.5$ and $\vartheta_{\mathcal{D}}=0.75$ across several random seeds}
%    \label{fig:peaks Dorfler}
%\end{figure}

\section{Conclusion and outlook}
The identification of an “optimal''  neural network topology is a very important, yet very challenging task in machine learning. In this work we have presented a time-continuous viewpoint on this problem, by introducing neural ODEs and employing an adaptive (time) meshing technique. Our dual-weighted goal-oriented approach aims at optimally resolving a prescribed target quantity and, under moderate regularization, appears to produce a stable adaptive hierarchy that improves the validation metrics of interest. This indeed encourages to develop and refine the method further.

To this end, one may identify several avenues of research: (i) In our approach, we inserted layers of a pre-defined “constant'' size. This may be made more flexible by admitting variable layer-sizes. Borrowing once again terminology from adaptive finite element methods, the number of neurons per layer may be related to a spatial discretization of a (within a layer) continuous setting for $\boldsymbol{x}$. One may then attempt to carry over space-time adaptation concepts from adaptive finite element discretizations. This may also include more complex target quantities in order to stabilize the adaptation process. (ii) The produced hierarchy of meshes induces a hierarchy of objective functions and constraint realizations associated with the continuous neural ODE constrained optimization problem. Upon the identification of suitable transfer operators within this hierarchy, one may develop (stochastic) multilevel minimization algorithms for further computational speed-up. (iii) We close by mentioning another possible direction which is connected to carrying our approach over to more complex network structures as, for instance, in convolutional neural networks (CNNs), where specific sub-blocks of a network are dedicated to specific tasks. Incorporating the latter in an error estimator / indicator based technique would be beneficial when trying to keep the overall number of networks parameters as small as possible. Clearly, one has to overcome challenges such as the non-locality of convolutional layers, pooling layers, local connectivity etc.

\begin{funding}
M. Hintermüller acknowledges the support of the Deutsche Forschungsgemeinschaft
(DFG, German Research Foundation) under Germany’s Excellence Strategy– The Berlin Mathematics Research Center MATH+ (EXC-2046/1, project ID: 390685689). M. Hintermüller and D. Korolev additionally acknowledge the support of the Federal Ministry of Education and Research, Germany (funding reference: 01IS24081) under project HybridSolver. M. Hinze acknowledges funding of the Federal Ministry of Education and Research, Germany (funding reference: 16DHBKI039) under project IH-evrsKI.
\end{funding}

%\section*{CRediT Authorship Contribution Statement}

%M. Hintermüller: , M. Hinze: , D. Korolev: 

\appendix
\section{Discontinuous Galerkin formulation of neural ODE}
\label{sec:App A}

We derive the variational formulation \eqref{VF state} which is suitable for our discontinuous Galerkin approach. Starting from the form \eqref{state form}, we split it over the time intervals as follows:
\begin{align}
\mathcal{F}(\boldsymbol{x} , \theta; \varphi) : = \sum_{k=1}^{K}\int_{I_{k}}\big(\dot{\boldsymbol{x} } - N_{F}(\boldsymbol{x} , \theta), \ \varphi_{1} \big) \ dt  + \big(\boldsymbol{x}(0) - \boldsymbol{x} _{\text{in}}, \varphi_{0}\big).
\end{align}
We integrate by parts: 
\begin{align*}
\int_{I_{k}}\big(\dot{\boldsymbol{x} } - N_{F}(\boldsymbol{x} , \theta),  \varphi_{1} \big) \, dt = \int_{I_{k}} -\big(\boldsymbol{x}, \dot{\varphi}_{1} \big) - \big(N_{F}(\boldsymbol{x}, \theta), \varphi_{1}  \big) \  dt   + \big(\boldsymbol{x}(t_{k}^{-}),  \varphi_{1}(t_{k}^{-}) \big) - \big(\boldsymbol{x}(t_{k-1}^{+}),  \varphi_{1}(t_{k-1}^{+})\big). 
\end{align*}
We perform a downwind approximation by replacing $ \varphi_{1}(t_{k}^{-})$ with $ \varphi_{1}(t_{k}^{+})$, and integrate by parts once more  to obtain
\begin{align*}
\int_{I_{k}}\big(\dot{\boldsymbol{x} } - N_{F}(\boldsymbol{x} , \theta),  \varphi_{1} \big) \, dt  + (\boldsymbol{x}(t_{k}^{+}), \varphi_{1}(t_{k}^{-})) - (\boldsymbol{x}(t_{k}^{-}), \varphi_{1}(t_{k}^{-})).
\end{align*}
Recalling the definition of the jump, we define the following form:
\begin{equation*}
\begin{aligned}
\mathcal{F}_{\mathrm{DG}}(\boldsymbol{x}, \theta; \varphi): = \sum_{k=1}^{K} \bigg[\int_{I_{k}}\big(\dot{\boldsymbol{x} } - N_{F}(\boldsymbol{x} , \theta),  \varphi_{1} \big) \, dt  + \big(\jump{\boldsymbol{x}}^{k}, \varphi_{1}(t_{k}^{-})\big) \bigg]  + \big(\boldsymbol{x}(0) - \boldsymbol{x}_{\text{in}},   \varphi_{0}\big),
\end{aligned}
\end{equation*}
which is the form required in \eqref{VF state}.  

The variational formulation of the adjoint equation which is suitable for our discontinuous Galerkin approach is derived as follows. First, we compute the derivative $ D_{1}\mathcal{F}_{\text{DG}}(\boldsymbol{x}_{\tau}, \theta_{\tau}; \boldsymbol{p}_{\tau}, \varphi_x)$ in the direction of $\varphi_{x} \in \mathcal{W}_{\tau}$, and obtain
\begin{equation}\label{adjoint der 1}
\begin{aligned}
D_{1}\mathcal{F}_{\text{DG}}(\boldsymbol{x}_{\tau}, \theta_{\tau}; \boldsymbol{p}_{\tau}, \varphi_x) & = \sum_{k=1}^{K} \int_{I_{k}} \big( \boldsymbol{p}_{\tau}, \dot{\varphi}_{x} \big) - \big( D_{1}N_{F}(\boldsymbol{x}_{\tau}, \theta_{\tau})^{\ast} \boldsymbol{p}_{\tau}, \varphi_{x}\big) \, dt  \\ & + \sum_{k=1}^{K} \big( \jump{\varphi_{x}}^{k},  \boldsymbol{p}_{\tau}(t_{k}^{-}) \big) + \big(\varphi_{x}(0), \boldsymbol{p}_{\tau}(0) \big).
\end{aligned}
\end{equation}
Integrating the first term in \eqref{adjoint der 1} by parts, we get 
\begin{align*}
\int_{I_{k}} \big(\boldsymbol{p}_{\tau}, \dot{\varphi}_{x} \big) \, dt = \big( \boldsymbol{p}_{\tau}(t_{k}^{-}),  \varphi_{x}(t_{k}^{-}) \big) - \big( \boldsymbol{p}_{\tau}(t_{k-1}^{+}),  \varphi_{x}(t_{k-1}^{+}) \big) - \int_{I_{k}} \big(\dot{\boldsymbol{p}}_{\tau}, \varphi_{x} \big) \, dt.
\end{align*}
According to our jump notation, we have $ \jump{\varphi_{x}}^{k} = \varphi_{x}(t_{k}^{+}) - \varphi_{x}(t_{k}^{-})$. Therefore, we get
\begin{equation}\label{adjoint der 2}
\begin{aligned}
\int_{I_{k}} \big(\boldsymbol{p}_{\tau}, \dot{\varphi}_{x} \big) \, dt  +  \big( \jump{\varphi_{x}}^{k},  \boldsymbol{p}_{\tau}(t_{k}^{-}) \big) & =   \int_{I_{k}} \big(-\dot{\boldsymbol{p}}_{\tau}, \varphi_{x} \big) \, dt +\big(\boldsymbol{p}_{\tau}(t_{k}^{-}), \varphi_{x}(t_{k}^{+}) \big) - \big(\boldsymbol{p}_{\tau}(t_{k-1}^{+}), \varphi_{x}(t_{k-1}^{+}) \big).
\end{aligned}
\end{equation}
Observe further that
\begin{equation}\label{adjoint der3}
\begin{aligned}
 \sum_{k=1}^K \big(\boldsymbol{p}_{\tau}(t_{k}^{-}), \varphi_{x}(t_{k}^{+}) \big) &= \sum_{k=1}^{K} \big( \boldsymbol{p}_{\tau}(t_{k-1}^{-}), \varphi_{x}(t_{k-1}^{+})\big)  - \big( \boldsymbol{p}_{\tau}(0), \varphi_{x}(0)\big) +   \big( \boldsymbol{p}_{\tau}(T), \varphi_{x}(T)\big).
\end{aligned}
\end{equation}
Using the identity \eqref{adjoint der 2}  together with \eqref{adjoint der 1} and \eqref{adjoint der3} yields the first equation in \eqref{adjoint equation}. Computing $D_{1}\mathcal{J}(\boldsymbol{x}_{\tau}, \theta_\tau; \varphi_\theta)$ as in \eqref{adjoint RHS} gives the second equation in \eqref{adjoint equation}.

\section{Time-marching interpretation of the DG(0) approximation scheme.}
\label{sec:App B}
Here, we derive the forward Euler interpretation of the discrete state equation and the backward Euler interpretation of the discrete adjoint equation in \eqref{opt system discrete}. For $r=0$ in $\mathcal{W}_{\tau}$, the ansatz $\boldsymbol{x}_{\tau} \in \mathcal{W}_{\tau}$ is given by
\begin{align*}
\boldsymbol{x}_{\tau}(t) = \sum_{k=1}^{K} \boldsymbol{x}_{\tau}^{k-1} \psi^{k-1}(t) + \boldsymbol{x}_{\tau}^{K} \chi^{K}(t),
\end{align*}
where $\boldsymbol{x}_{\tau}^{k} \in \mathbb{R}^{md}$ for all $k=0,\ldots,K$, and 
\begin{align*}
\psi^{k-1}(t) = \begin{cases}
1, \ \ t \in [t_{k-1}, t_k), \\
0, \ \ \text{elsewhere},  
\end{cases}, \quad \chi^{K}(t) = \begin{cases}
1, \ \ t=t_{K}, \\ 
0, \ \ \text{elsewhere.}
\end{cases}
\end{align*}
To describe the test functions from $\mathcal{V}_{\tau}$ with $r=0$, we use $\boldsymbol{v}_{\tau}^{0}\chi^{0}(t)$, and $\boldsymbol{v}_{\tau}^{k} \ \xi^{k}(t)$ for $k \in \{1,\ldots,K\}$, where $\boldsymbol{v}_{\tau}^{k} \in \mathbb{R}^{md}$ for each $k \in \{0,\ldots,K\}$, and 
\begin{align*}
\xi^{k}(t) = \begin{cases}
1, \ \ t \in (t_{k-1}, t_k], \\
0, \ \ \text{elsewhere},  
\end{cases}, \quad \chi^{0}(t) = \begin{cases}
1, \ \ t=0, \\ 
0, \ \ \text{elsewhere.}
\end{cases}
\end{align*}
Since the trial functions are constant in time on $I_{k}$, we get $\dot{\boldsymbol{x}}_{\tau}=0$, $\boldsymbol{x}_{\tau}(t_{k}^{+})=\boldsymbol{x}^{k}_{\tau}$ and $\boldsymbol{x}_{\tau}(t_{k}^{-})=\boldsymbol{x}^{k-1}_{\tau}$. Thus, from the discontinuous formulation \eqref{VF state}, we get 
\begin{equation}
\begin{aligned}
\mathcal{F}_{\text{DG}}(\boldsymbol{x}_\tau, \theta_{\tau}; \boldsymbol{v}_{\tau})& = \sum_{k=1}^{K} \int_{I_{k}} -\big(N_{F}(\boldsymbol{x}_\tau, \theta_{\tau}),  \boldsymbol{v}_{\tau}^{k} \big) \, dt  + \big(\boldsymbol{x}_{\tau}^{k}- \boldsymbol{x}_{\tau}^{k-1}, \boldsymbol{v}_{\tau}^{k}\big)  + \big(\boldsymbol{x}_{\tau}^{0} - \boldsymbol{x}_{\text{in}},   \boldsymbol{v}_{\tau}^{0}\big),
\end{aligned}
\end{equation}
which yields the following expression after rearrangement:
\begin{align}
\big( \boldsymbol{x}_{\tau}^{k} -  \boldsymbol{x}_{\tau}^{k-1}, \boldsymbol{v}_{\tau}^{k} \big)  = \int_{I_{k}}\big(N_{F}(\boldsymbol{x}_\tau, \theta_{\tau}),  \boldsymbol{v}_{\tau}^{k} \big) \, dt, \quad \boldsymbol{x}_{\tau}^{0} = \boldsymbol{x}_{\text{in}}.
\end{align}
Discretizing the right-hand side using the midpoint quadrature rule and with the standard canonical basis in $\mathbb{R}^{md}$, we obtain an explicit ResNet-type time-stepping scheme:
\begin{align}\label{resnet}
\boldsymbol{x}_{\tau}^{0} = \boldsymbol{x}_{\text{in}}, \quad \boldsymbol{x}_{\tau}^{k} = \boldsymbol{x}_{\tau}^{k-1} + h_{k}F(\boldsymbol{x}_\tau^{k-1}, \theta_{\tau}(t_{k-1/2})), 
\end{align}
where $k \in \{1, \ldots, K\}$ and $t_{k-1/2}:=\tfrac{t_{k-1} + t_{k}}{2}$ is the midpoint of the interval $I_{k}$.

For the discretization of the adjoint equation \eqref{adjoint equation}, we use the ansatz $\boldsymbol{p}_{\tau} \in \mathcal{V}_{\tau}$, which is given by
\begin{align*}
\boldsymbol{p}_{\tau}(t) = \boldsymbol{p}_{\tau}^{0} \chi^{0}(t) + \sum_{k=1}^{K}\boldsymbol{p}_{\tau}^{k} \ \xi^{k}(t),
\end{align*}
and apply the test functions from $\mathcal{W}_{\tau}$. For the latter, we have $\boldsymbol{w}_{\tau}^{k-1} \psi^{k-1}(t)$ for $k\in \{1,\ldots, K\}$, and  $\boldsymbol{w}_{\tau}^{K} \chi^{K}(t)$, where $\boldsymbol{w}_{\tau}^{k} \in \mathbb{R}^{md}$ for all $k\in \{0,\ldots, K\}$. For the trial functions from $\mathcal{V}_{\tau}$, we get $\dot{\boldsymbol{p}}_{\tau} = 0$, $\boldsymbol{p}_{\tau}(t_{k-1}^+) = \boldsymbol{p}_{\tau}^k$ and $\boldsymbol{p}_{\tau}(t_{k-1}^-) = \boldsymbol{p}_{\tau}^{k-1}$ in  \eqref{adjoint equation}, which yields
\begin{equation}
\begin{aligned}
\sum_{k=1}^{K} \int_{I_{k}} \big (- D_{1}N_{F}(\boldsymbol{x}_{\tau}, \theta_{\tau})^{\ast}\boldsymbol{p}_{\tau}, \boldsymbol{w}_{\tau}^{k-1}) \ dt   &- \sum_{k=1}^{K}\big (\boldsymbol{p}_{\tau}^{k}- \boldsymbol{p}_{\tau}^{k-1} , \boldsymbol{w}_{\tau}^{k-1} \big) = 0 \\
 \big( \boldsymbol{p}_{\tau}(T), \boldsymbol{w}_{\tau}^{K} \big) &=  \big(l^{\prime}(\boldsymbol{x}_{\tau}(T)), \boldsymbol{w}_{\tau}^{K} \big).
\end{aligned}
\end{equation}
Rearranging, we get 
\begin{align*}
\big(\boldsymbol{p}_{\tau}^{k-1}, \boldsymbol{w}_{\tau}^{k-1} \big)  &=  \big(\boldsymbol{p}_{\tau}^{k}, \boldsymbol{w}_{\tau}^{k-1} \big)  + \int_{I_{k}} \big(D_{1}N_{F}(\boldsymbol{x}_{\tau}, \theta_{\tau})^{\ast}\boldsymbol{p}_{\tau}, \boldsymbol{w}_{\tau}^{k-1} \big) \ dt, \\ 
( \boldsymbol{p}_{\tau}^{K}, \boldsymbol{w}_{\tau}^{K} \big) &= \big(l^{\prime}(\boldsymbol{x}^{K}_{\tau}), \boldsymbol{w}_{\tau}^{K}\big).
\end{align*}
By testing with the standard canonical basis in $\mathbb{R}^{md}$, and approximating the integral with the midpoint quadrature rule, we get the time-marching scheme
\begin{equation}\label{discrete adjoint TM}
\begin{aligned}
\boldsymbol{p}_{\tau}^{k-1} &= \boldsymbol{p}_{\tau}^{k}  + h_{k} D_{1}F\big(\boldsymbol{x}_{\tau}^{k-1}, \theta_{\tau}(t_{k-1/2})\big)^{\ast}\boldsymbol{p}_{\tau}^{k}, \quad k=K, \dots , 1,\\
\boldsymbol{p}^{K}_{\tau} & =  l^{\prime}(\boldsymbol{x}^{K}_{\tau}).
\end{aligned}
\end{equation}
The right-hand side, involving the product of the Jacobian and $\boldsymbol{p}_{\tau}^{k}$ in \eqref{discrete adjoint TM}, is typically produced via backward-mode automatic differentiation without explicitly forming the Jacobian; see, e.g., \cite{chen2018neural}.

\printbibliography

\end{document}